%% file: 1DiagrCalc-2_newd.tex
\documentclass[a4paper,11pt,reqno]{amsart}

\usepackage{graphicx}
\usepackage{amsmath}
\usepackage{amssymb}
\usepackage{amsfonts}
\usepackage{epsfig}
\usepackage{epigraph}
\usepackage{comment}
\usepackage{color}

\usepackage{mathrsfs}

\usepackage{graphics}
\usepackage{float}

\graphicspath{{Linkages/}
              {Conj_Classes/}
              {Loctets/}}

\def\colorcom#1#2#3{
 \def#1##1 {
  \begin{color}{#3}
      \bf[\,#2: ##1\,]
  \end{color}}
} \colorcom\RS{RS}{red}

\makeatletter

\newtheorem{theorem}{Theorem}[section]

\@addtoreset{equation}{section}

\def\botL{\stackrel{}{L}}
\def\botLv{\stackrel{}{L^{\vee}}}

\newtheorem{corollary}[theorem]{Corollary}
\newtheorem{remark}[theorem]{Remark}
\newtheorem{proposition}[theorem]{Proposition}
\newtheorem{lemma}[theorem]{Lemma}

\def\PerfProof{{\it Proof.\ }}

\begin{document}

\title[Root systems and diagram calculus. II]
{\qquad\qquad Root systems and diagram calculus. \newline
  II. Quadratic forms for the Carter diagrams}
         \author{Rafael Stekolshchik}

\date{}


\begin{abstract}
   For any Carter diagram $\Gamma$ containing $4$-cycle, we introduce the partial Cartan matrix $B_L$, which is similar to the Cartan matrix associated with a Dynkin diagram. A linkage diagram is obtained from $\Gamma$ by adding one root together with its bonds such that the resulting subset of roots is linearly independent. The linkage diagrams connected under the action of dual partial Weyl group (associated with $B_L$) constitute the linkage system, which is similar to the weight system arising in the representation theory of the semisimple Lie algebras.  For Carter diagrams $E_6(a_i)$ and $E_6$ (resp. $E_7(a_i)$ and  $E_7$; resp. $D_n(a_i)$ and $D_n$),  the linkage system has, respectively, $2$, $1$, $1$ components,  each of which contains, respectively, $27$, $56$, $2n$ elements. Numbers $27$, $56$ and $2n$ are well-known dimensions of the smallest fundamental representations of semisimple Lie algebras, respectively, for $E_6$, $E_7$ and $D_n$.  The $8$-cell {\lq\lq}spindle-like{\rq\rq} linkage subsystems called loctets play the essential role in describing the linkage systems. It turns that weight systems also can be described by means of loctets.
\end{abstract}

\maketitle

\begin{center}
  \vspace*{2cm}
   \dedicatory{\it In memory of Gurii Fedorovich Kushner}
   \vspace*{3cm}
\end{center}

\tableofcontents

\newpage
~\\
~\\

\setlength{\epigraphwidth}{100mm}

\epigraph{
It seemed completely mad.
It seemed so mad, in fact, that Killing was rather upset that the exceptional
groups existed, and for a time he hoped they were a mistake that he
could eradicate. They spoiled the elegance of his classification. But they
were {\it there}, and we are finally beginning to understand {\it why} they are there.
In many ways, the five exceptional Lie groups now look much more interesting
than the four infinite families. They seem to be important in particle
physics, as we will see; they are definitely important in mathematics. And
they have a secret unity, not yet fully uncovered...
 }{Ian Stewart, \\ Why beauty is truth: a history of symmetry, \cite[p. 170]{S07}, 2007}

\section{\sc\bf Introduction}

\subsection{The linkage diagrams and linkage labels}
  \label{sec_linkage_diagr}
  We consider a class of diagrams called {\it linkage diagrams} that
  constitute the subclass of the class of connection diagrams introduced in \cite{St10}
  and generalize the Carter diagrams (= admissible diagrams) introduced
  by R.~Carter in \cite{Ca72} for the classification of conjugacy classes in a finite Weyl group $W$.
  The linkage diagram is obtained from a Carter diagram $\Gamma$ by adding one extra root $\gamma$  with its bonds such that the roots corresponding to vertices of $\Gamma$ together with $\gamma$ form some
  linearly independent root subset. The extra root $\gamma$ added to the Carter diagram $\Gamma$ is called a {\it linkage}, see Section \ref{sec_linkage}.  Any linkage diagram constructed by this way may be also a Carter diagram but this is not necessarily so.  The following inclusions hold:

 \begin{equation*}
   \fbox{$\begin{array}{c}
      \text{ Dynkin } \\
      \text{ diagrams of CCl\footnotemark[1] }
   \end{array}$} \quad\subset\quad
   \fbox{$\begin{array}{c}
      \text{ Carter } \\
      \text{ diagrams }
   \end{array}$} \quad\subset\quad
   \fbox{$\begin{array}{c}
      \text{ Linkage } \\
      \text{ diagrams }
   \end{array}$} \quad\subset\quad
   \fbox{$\begin{array}{c}
      \text{ Connection } \\
      \text{ diagrams }
   \end{array}$}
 \end{equation*}

 \footnotetext[1]{The Dynkin diagrams in this article appear in two ways:
 (1) associated with some Weyl group (customary use);
 (2) representing some conjugacy class (CCl), i.e, the Carter diagram
 which looked like a Dynkin diagram.
 In a few cases Dynkin diagrams represent two (and even three!) conjugacy classes, see Remark \ref{two_class_Al}.}

  With every linkage diagram we associate the {\it linkage labels vector}, or, for short,
 {\it linkage labels}. The linkage labels are similar to the Dynkin labels,
 see \cite{Sl81}, which are the {\lq\lq}numerical labels{\rq\rq} introduced by Dynkin in \cite{Dy50} for the study of irreducible linear representations of the semisimple Lie algebras, \cite{GOV90}, \cite{KOV95},
 \cite{Ch84}.

\begin{figure}[h]
\centering
\includegraphics[scale=0.7]{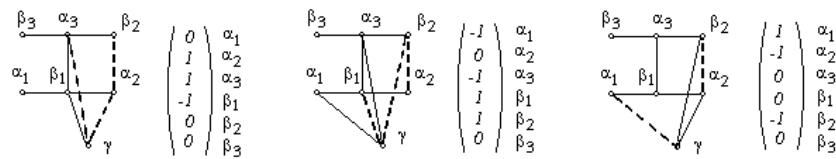}
 \caption{\hspace{3mm}Examples of linkage diagrams and linkage labels vectors for $E_6(a_1)$.}
\label{E6a1_exam_linkages}
\end{figure}

 For any simply-laced Carter diagram $\Gamma$, a linkage label
 takes one of three values $\{-1, 0, 1 \}$. There is one-to-one correspondence between
 the linkage diagrams obtained from the given simply-laced Carter diagram $\Gamma$ and the linkage label vectors taking coordinates from the set $\{-1, 0, 1 \}$.
 For this reason, we can use terms \underline{linkage labels} and  \underline{linkage diagrams} as
 convertible. Some linkage diagrams and their linkage labels for the Carter diagram $E_6(a_1)$  are depicted in Fig. \ref{E6a1_exam_linkages}.
  The linkage diagrams are the focus of this paper. We give a complete description of linkage diagrams constructed for every simply-laced Carter diagram containing $4$-cycle or a branch point, i.e., containing
  $D_4(a_1)$ or $D_4$. By abuse of notation, this description is essentially based on the following issue:
 ~\\

  {\it What vectors can be added to the irreducible linearly independent root subset,
  so that the resulting set would also be some irreducible linearly independent root subset?}
 ~\\
 ~\\
 It turns out that the answer to this question is very simple within the framework of the quadratic
 form associated with the Cartan matrix, see Section \ref{sec_criterion}.

\subsection{The Cartan matrix for a conjugacy class}
 \label{sec_part_Cartan}
  We consider two classes of simply-laced connected Carter diagrams:
  denote by $\mathsf{C4}$ the class of diagrams containing $4$-cycle $D_4(a_1)$,
  and by $\mathsf{DE4}$ the class of diagrams without cycles and containing $D_4$ as a subdiagram,
  i.e., $\mathsf{DE4}$ is the class consisting of Dynkin diagrams  $E_6$, $E_7$, $E_8$ and $D_l$ for $l \geq 4$.
  Conjugate elements in the Weyl group $W$ are associated with the same Carter diagram
 $\Gamma$. The converse is not true, the Carter diagram $\Gamma$ does not determine a
 single conjugacy class in $W$, \cite[Lemma 27]{Ca72}. Nevertheless, the converse
 statement takes place for $\mathsf{C4}$ and $\mathsf{DE4}$.

 \begin{theorem}
  \label{th_4.1.}
   Let $\Gamma$ be the diagram belonging to $\mathsf{C4}$ or to $\mathsf{DE4}$.

   Then $\Gamma$  determines only one conjugacy class. Two root subsets
 \begin{equation}
   \label{canon_dec_1}
    S = \{ \tau_1, \dots, \tau_{l} \}   \text{ and  }
    S'= \{ \tau'_1, \dots, \tau'_{l} \}
 \end{equation}
  corresponding to the same diagram $\Gamma$ are equivalent, i.e., there exists the element $U \in W$ such that
 \begin{equation}
   \label{uniq_B}
   \begin{split}
    & U\tau_i = \tau'_i, \text{ where }  i = 1, \dots, l,  \text{ and } \\
    & (\tau_i, \tau_j) = (\tau'_i, \tau'_j).
   \end{split}
 \end{equation}
 \end{theorem}
   For $\Gamma \in \mathsf{C4}$  the theorem follows from  \cite{St10}, Theorem 4.1 and Section 4.1.

   For $\Gamma \in \mathsf{DE4}$  the theorem follows from Proposition \ref{diagr_El_Dl}.

   For $A_n$ the theorem does not hold, see Remark \ref{two_class_Al}.
 \qed
~\\

 Let $L$ be the linear space spanned by the roots associated with $\Gamma$,
 $\gamma_L$ be the projection of the linkage $\gamma$ on $L$. The linkage labels vector
 is the element of the dual linear space $L^{\vee}$, see Section \ref{sec_linkage}. We denote
 the linkage labels vector by $\gamma^{\vee}$.
  For any Carter diagram $\Gamma \in \mathsf{C4}$,
  and consequently, for the conjugacy class associated with $\Gamma$, we introduce the {\it partial
  Cartan matrix} $B_L$ which is similar to the Cartan matrix ${\bf B}$ associated with a Dynkin diagram.  Thanks to Theorem \ref{th_4.1.} the matrix $B_L$ is well-defined, see Section \ref{sec_partial_B}.
  The matrix $B_L$ maps $\gamma_L \in L$ to the linkage labels $\gamma^{\vee} \in  L^{\vee}$ as follows:
 \begin{equation*}
   \label{B_pairing}
      \gamma^{\vee} = B_L \gamma_L, \quad B_L^{-1}\gamma^{\vee} =  \gamma_L,
 \end{equation*}
 see Proposition \ref{prop_link_diagr_conn}. In the dual space $L^{\vee}$ we take the quadratic form
 $\mathscr{B}^{\vee}_L$ associated with the inverse matrix $B_L^{-1}$. The quadratic form $\mathscr{B}^{\vee}_L$  provides the easily verifiable criterion that the vector $u^{\vee}$ is the linkage labels vector for a certain linkage $\gamma \not\in L$. This criterion
 (Theorem \ref{th_B_less_2}) is the following inequality:
  \begin{equation*}
      \mathscr{B}^{\vee}_L(\gamma^{\vee}) < 2.
  \end{equation*}

  \subsection{The linkage systems and loctets}

  A certain group $W^{\vee}_L$ named the dual partial Weyl group acts in
  the dual space $L^{\vee}$. This group acts on the linkage label vectors, i.e.,
  on the set of linkage diagrams:
 \begin{equation*}
      (w\gamma)^{\vee} = w^{*}\gamma^{\vee},
 \end{equation*}
 where $w^{*} \in W^{\vee}_L$, see Proposition \ref{prop_link_diagr_conn}. The set of
 linkage diagrams (=linkage labels) under action of $W^{\vee}_L$ constitute the
 diagram called the {\it linkage system} similarly to the weight system\footnotemark[1]
 in the theory of representations  of semisimple Lie algebras, \cite[p. 30]{Sl81}.
 We denote by $\mathscr{L}(\Gamma)$ the linkage system associated with the Carter diagram $\Gamma$.
 \footnotetext[1]{Frequently, in the literature (see, for example, \cite{Va00}), the term
  \underline{weight diagram} is used instead of the term \underline{weight system}. However, the term {\lq\lq}diagram{\rq\rq} is heavily overloaded in our context.}

\begin{figure}[H]
\centering
\includegraphics[scale=1.2]{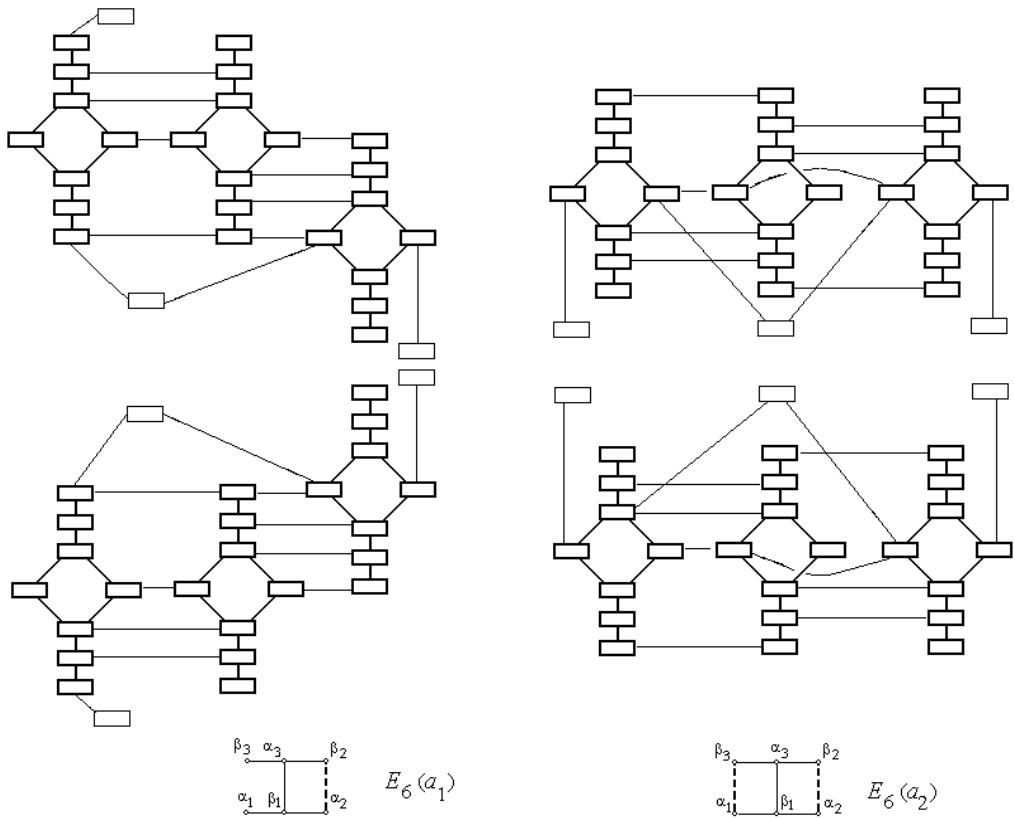}
\caption[\hspace{3mm}Linkage systems for $E_6(a_1)$ and $E_6(a_2)$. The $8$-cell bold subdiagrams are {\bf loctets}]{\hspace{3mm}Linkage systems for $E_6(a_1)$ and $E_6(a_2)$, see Fig. \ref{E6a1_linkages} and Fig. \ref{E6a2_linkages}. The $8$-cell bold subdiagrams are {\bf loctets}, see Fig. \ref{fig_loctets}}
\label{linkage_systems_example}
\end{figure}

 The linkage systems for the Carter diagrams $E_6(a_1)$ and $E_6(a_2)$ depicted in Fig. \ref{linkage_systems_example}. The linkage systems with linkage labels vectors
 for all Carter diagrams are presented in Fig. \ref{D4a1_linkages}-\ref{E7a4_linkages},
  \ref{Dk_al_wind_rose},   \ref{E6pure_loctets_norm_gam8},  \ref{E7pure_linkage_system}, \ref{D5pure_loctets},
  \ref{D6pure_loctets}, \ref{D7a1_linkages_cdef}-\ref{D7a1_D7a2_D7pu_loctets_comp3}.
 Every linkage diagram containing at least one non-zero $\alpha$-label (see Section \ref{sec_Carter})
 belongs to a certain $8$-cell "spindle-like" linkage subsystem called
 \underline{loctet} (= linkage octet). The loctets are the main construction blocks in every linkage system. If all $\alpha$-labels (resp. $\beta$-labels) of the linkage diagram $\gamma^{\vee}$ are zero, we call $\gamma^{\vee}$ the {\it $\beta$-unicolored} (resp. {\it $\alpha$-unicolored}) linkage diagram. Every linkage system is the union of several loctets and several $\beta$-unicolored linkage diagrams, see Section \ref{sec_enum_loctets}. There are exactly $6$ loctets in any linkage system of type $E_6(a_i)$ for $i = 1, 2$ or type $E_7(a_i)$ for $i = 1, 2, 3, 4$.  For $E_6(a_1)$ and $E_6(a_2)$,  the linkage system is the union of $6$ loctets and $6$ $\beta$-unicolored linkage diagrams, in sum
 $54 = 2 \times 27 = 6 \times 8 + 6$ linkage diagrams, see Fig. \ref{linkage_systems_example}.

  The knowledge of the structure of linkage systems is the essential part in our proof of the Carter theorem
  on decomposition of every element in the Weyl group into the product of two involutions, see \cite{St11}.

\subsection{The numbers $27$, $56$ and $2l$}
 We observe that numbers of linkage diagrams for Carter diagrams $E_6(a_i)$,
 $E_7(a_i)$ and $D_l(a_i)$ are, respectively, $27$, $56$ and $2l$ that coincide with the dimensions of the smallest fundamental representations of semisimple Lie algebras, respectively, $E_6$,  $E_7$, and $D_l$. Moreover, the number of components in the linkage systems $E_6(a_i)$, $E_7(a_i)$ and $D_l(a_i)$ are, respectively, equal to $2$, $1$ and $1$ that coincide with the number of different smallest fundamental representations of  semisimple Lie algebras $E_6$, $E_7$ and $D_l$, see \cite[Ch. VIII, Tabl. 2]{Bo05}, \cite[Tabl. 30]{Dy52}. Of course, these facts require {\it a priori} reasoning. It seems that ideas of Section \ref{sec_proj} concerning the projection of linkage systems and Theorem \ref{th_map_linkage_systems} may help to explain these facts.

\subsection{The Carter diagrams, connection diagrams and linkage diagrams}
  \label{sect_diagram}
  Several types of diagrams are considered in this paper.
  The {\it Carter diagram} introduced by R.~Carter \cite{Ca72}
  describes a bicolored decomposition of some element $w \in W$, see Section \ref{sec_Carter}.
  The {\it connection diagrams} introduced in \cite[\S1.1]{St10}
  generalize the Carter diagram; the connection diagram describes a decomposition
  of certain element $w \in W$,
  and this diagram is supplied with an order of reflections $\Omega$,  see Section \ref{sec_connection}.
  In both cases all reflections are associated with roots which are not necessary simple.
  The {\it linkage diagram} is a particular case of the connection diagram obtained from a certain Carter diagram by adding one extra vertex with its bonds, see Section \ref{sec_linkage}.
  The linkage diagrams are the focus of this paper.

\subsubsection{The Dynkin diagrams}
  \label{sec_Dynkin}
  Let $\Gamma$ be a certain Dynkin diagram, $\varPhi$ be the associated root system,
  $\Pi$ be the set of all simple roots in $\varPhi$,
  $E$ be the linear space spanned by all roots,
  $W$ be the finite Weyl group associated with $\Gamma$ and acting in the linear space $E$.
  Let ${\bf B}$ be the corresponding Cartan matrix, $(\cdot, \cdot)$ be the corresponding symmetric bilinear form, and $\mathscr{B}$ be the quadratic Tits
  form associated with ${\bf B}$,  \cite[Ch. 2]{St08}.
  We suppose that the diagonal elements of ${\bf B}$ are $2$, see Remark \ref{rem_classic}. The following relation
  is the well-known property connecting roots and the quadratic Tits forms\footnotemark[1]:
  \begin{equation}
     \label{eq_Kac}
      \mathscr{B}(\alpha) = 2 \Longleftrightarrow \alpha \in \varPhi.
  \end{equation}
  For two non-orthogonal simple roots $\alpha$, $\beta$, we have
  \begin{equation}
     \label{eq_Kac_2}
      (\alpha, \beta) = \Arrowvert \alpha \Arrowvert \Arrowvert \beta \Arrowvert
      \cos(\widehat{\alpha, \beta}) = \sqrt{2} \cdot \sqrt{2} (-\frac{1}{2}) = -1.
  \end{equation}

  \footnotetext[1]{In order to obtain the values of the linkage labels (see Section \ref{sec_linkage})
  by integers as in \eqref{eq_Kac_2}, we choose the diagonal elements equal $2$. Frequently, diagonal elements are chosen equal $1$, and \eqref{eq_Kac} looks as follows: $\mathscr{B}(\alpha) = 1 \Longleftrightarrow \alpha \in \varPhi$,
  see \cite{Kac80}.}

\subsubsection{The Carter diagrams}
  \label{sec_Carter}

  The Carter diagram (= admissible diagram) \cite[\S4]{Ca72}
  is the diagram $\Gamma$  satisfying two conditions:
  \vspace{3mm}

  (a) The nodes of $\Gamma$ correspond to a set of linearly independent roots.

  (b) Each subgraph of $\Gamma$ which is a cycle contains even number of vertices.
   \vspace{3mm}

Let $w = w_1 w_2$ be the decomposition of $w$ into the product of
two involutions. By \cite[Lemma 5]{Ca72} each of $w_1$ and $w_2$ can
be expressed as products of reflections corresponding to mutually orthogonal roots as
follows:
\begin{equation}
   \label{two_invol}
       w = w_1 w_2, \quad
       w_1 = s_{\alpha_1} s_{\alpha_2} \dots s_{\alpha_k}, \quad
       w_2 = s_{\beta_1} s_{\beta_2} \dots s_{\beta_h}, \quad
      \text{where} \quad  k + h = l_C(w).
 \end{equation}
For details, see \cite[\S4]{Ca72}, \cite[\S1.1]{St10}.  
We denote by {\it $\alpha$-set} (resp. {\it $\beta$-set}) the subset of
 roots corresponding to $w_1$ (resp. $w_2$):
\begin{equation}
   \label{two_sets}
       \alpha\text{-set} = \{ \alpha_1, \alpha_2, \dots, \alpha_k \}, \quad
       \beta\text{-set} =  \{ \beta_1, \beta_2, \dots, \beta_h \}.
 \end{equation}
 Any coordinate from $\alpha$-set (resp. $\beta$-set) of the linkage labels vector
 we call {\it $\alpha$-label} (resp. {\it $\beta$-label}).
 We call the decomposition \eqref{two_invol} the {\it bicolored decomposition}.
 Let $L \subset E$ be the linear subspace spanned by root subsets \eqref{two_sets},
 $L^{\vee}$ be the dual linear space. The corresponding root basis which vectors
 are not necessarily simple roots, we denote by $\Pi_w$:
\begin{equation}
   \label{root_subset_L}
       \Pi_w = \{ \alpha_1, \alpha_2, \dots, \alpha_k, \beta_1, \beta_2, \dots, \beta_h \}.
 \end{equation}

\subsubsection{The connection diagrams}
  \label{sec_connection}
  Each element $w \in W$ can be expressed in the from
 \begin{equation}
   \label{conn_decomp}
    w = s_{\alpha_1}s_{\alpha_2}\dots_{\alpha_k}, \text{ where } \alpha_i \in \Phi,
 \end{equation}
 $\Phi$ is the root system associated with the Weyl group $W$, $s_{\alpha_i}$ are reflections in
 $W$ corresponding to not necessarily simple roots $\alpha_i \in \Phi$.
  The {\it connection diagram} is the pair
  $(\Gamma, \Omega)$, where $\Gamma$ is the diagram describing
  connections between roots as it is described by
  the Dynkin diagrams or by the Carter diagrams, and $\Omega$ is the
  order of elements in the (not necessarily bicolored) decomposition \eqref{conn_decomp},
  see \cite{St10}.

  For the Dynkin diagrams, a number of bonds for non-orthogonal roots
  describes the angle between roots, and the ratio of lengths of two roots.
  For the Carter diagrams and connection diagrams,
  we add designation distinguishing acute and obtuse angles
  between roots. Recall, that for the Dynkin diagrams, all angles between simple roots
  are obtuse and a special designation is not necessary.
  A {\it solid edge} indicates an obtuse angle between roots exactly as for simple roots
  in the case of Dynkin diagrams. A {\it dotted edge} indicates an acute angle between the
  roots considered. For details, see \cite{St10}.

\subsubsection{Linkages and linkage diagrams}
  \label{sec_linkage}
 Let $w = w_1w_2$ be the bicolored decomposition of some element
 $w \in W$, where $w_1$, $w_2$ are two involutions
 associated, respectively, with $\alpha$-set $\{ \alpha_1, \dots, \alpha_k \}$ and $\beta$-set
 $\{ \beta_1, \dots, \beta_h \}$  of roots from the root system $\varPhi$, see
 \eqref{two_invol}, \eqref{two_sets}, and let $\Gamma$ be the Carter
 diagram associated with this bicolored decomposition.
 We consider the {\it extension} of the root basis $\Pi_w$ by means of the root $\gamma \in \varPhi$,
 such that the set of roots
 \begin{equation}
   \label{alpha_beta}
    \Pi_w(\gamma) = \{ \alpha_1, \dots, \alpha_k, \beta_1, \dots, \beta_h, \gamma \}
 \end{equation}
 is linearly independent. Let us multiply $w$ on the right by the reflection
 $s_{\gamma}$ corresponding to $\gamma$ and consider the diagram
 $\Gamma' = \Gamma \cup \gamma$ together with new edges.
By \eqref{eq_Kac_2}, these edges are
 ~ \\
\begin{minipage}{10.8cm}
 \begin{equation*}
   \begin{cases}
  \text{{\it solid}, for $(\gamma, \tau) = -1$}, \\
  \text{{\it dotted}, for $(\gamma, \tau) = 1$}, \\
   \end{cases}
 \end{equation*}
 where $\tau$ one of elements \eqref{alpha_beta}. We call the diagram $\Gamma'$
 a {\it linkage diagram}, and the root $\gamma$ we call a {\it linkage} or a {\it $\gamma$-linkage}. The roots $\tau$ corresponding to the new edges ($(\gamma, \tau) \neq 0$)
 we call {\it endpoints} of the linkage diagram. Endpoints lying in $\alpha$-set (resp. $\beta$-set)
 we call {\it $\alpha$-endpoints} (resp. {\it $\beta$-endpoints}).
 Consider vectors $\gamma^{\vee}$ belonging to the dual space $L^{\vee}$ and defined by
 \eqref{dual_gamma}.
 We call vectors \eqref{dual_gamma} {\it linkage labels vectors} or, for brevity, {\it linkage labels}.
 \end{minipage}
\begin{minipage}{5.5cm}
  \quad
 \begin{equation}
   \label{dual_gamma}
   \gamma^{\vee} :=
   \left (
    \begin{array}{c}
    (\gamma, \alpha_1) \\
    \dots,     \\
    (\gamma, \alpha_k)  \\
    (\gamma, \beta_1)  \\
    \dots,  \\
    (\gamma, \beta_h) \\
    \end{array}
   \right )
 \end{equation}

\end{minipage}
~\\
 There is, clearly, the one-to-one correspondence between linkage labels vectors $\gamma^{\vee}$
 (with labels $\gamma^{\vee}_i \in \{0, -1, 1\}$) and simply-laced linkage diagrams
 (i.e., such linkage diagrams that
  $(\gamma, \tau) \in \{0, -1, 1\}$).

\input 2results.tex

\input 2a_more4cycles.tex
\input 3dualWeyl.tex

\input 4inverse.tex

\input 5loctets.tex

\input 6homogen.tex

\input 7Dynkin.tex

\input 5a_extendable.tex

\begin{appendix}
\input A4matricesB.tex

\input A5solutions.tex

\input A6linkages.tex
\end{appendix}
\newpage
\listoffigures
\newpage
\input notation.tex
\newpage
\input biblio-2.tex

\end{document}

%% file: 2results.tex
\subsection{The main results}

 \subsubsection{The partial Cartan matrix and dual partial Weyl group}
    \label{sec_criterion}
 In Section \ref{sec_Cartan_matr_ccl} we introduce the partial Cartan matrix $B_L$ associated with
 the linear subspace $L \subset E$ in such a way that $B_L$ coincide with the Cartan matrix ${\bf B}$
 restricted on $L$. The matrix $B_L$ is positive definite (Proposition \ref{restr_forms_coincide}).
 We introduce the partial Weyl group $W_L$ generated by reflections
  $\{s_{\tau_1}, \dots, s_{\tau_l} \}$, where $l = \dim L$, and the dual partial Weyl group $W^{\vee}_L$ generated by dual reflections  $\{s^*_{\tau_1}, \dots, s^*_{\tau_l} \}$,   where $s^*_{\tau_i}$ are the dual reflections associated with not necessarily simple roots $\tau_i$.
 The linkage diagrams $\gamma^{\vee}$ and $(w\gamma)^{\vee}$ are related as follows:
 $(w\gamma)^{\vee} = w^{*}\gamma^{\vee}$, where $w^{*} \in W^{\vee}_L$ (Proposition \ref{prop_link_diagr_conn}).
 The quadratic Tits form $\mathscr{B}_L$ takes a certain constant value for all elements
 $w\gamma$, where $w$ runs over  $W_L$ (Proposition \ref{prop_uniq_proj}).
 Let $\gamma_L$ is the projection of the root $\gamma$ on $L$.
 The main result of Section \ref{sec_Cartan_matr_ccl} is the following theorem that verifies whether or not
 a given vector is a linkage labels vector:
  ~\\

 {\bf Theorem} (Theorem \ref{th_B_less_2}).
     {\it A vector $u^{\vee} \in L^{\vee}$ is the linkage labels vector
     corresponding to a certain root $\gamma \in \varPhi$, $\gamma \not\in L$
     (i.e., $u^{\vee} = B_{L}\gamma_L$) if and only if
     \begin{equation}
       \label{eq_p_less_2_loc}
          \mathscr{B}^{\vee}_L(u^{\vee}) < 2.
     \end{equation}
   }

 \subsubsection{Three loctet types}

\begin{figure}[h]
\centering
\includegraphics[scale=1.3]{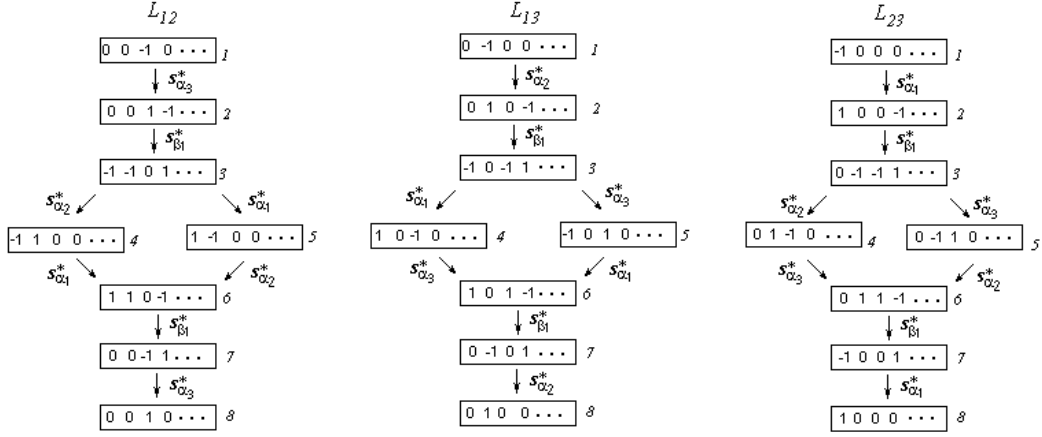}
 \caption{\hspace{3mm}The loctet types $L_{12}$, $L_{13}$ and $L_{23}$.}
\label{fig_loctets}
\end{figure}

 In Section \ref{sec_enum_loctets}, by means of inequality \eqref{eq_p_less_2_loc},
 we obtain a complete description of linkage diagrams for all linkage systems.
 We introduce the $8$-cell linkage subsystem called
 {\it loctet (= linkage octet)} as is depicted in Fig. \ref{fig_loctets}.
 Loctets are the main construction blocks in the structure of the linkage systems.
  Consider roots $\gamma^{\vee}_{ij}(n)$ depicted in Fig. \ref{fig_loctets},
  where $\{ij\}$ is associated with type $L_{ij}$ and $n$ is the order number of the linkage diagrams
  in the vertical numbering in Fig. \ref{fig_loctets}.
  The octuple of linkages depicted in every connected component
  in Fig \ref{fig_loctets} we call the {\it loctet } of type $L_{12}$
  (resp. $L_{13}$, resp. $L_{23}$).

 ~\\
  {\bf Corollary} (on the structure of loctets and linkage diagrams (Corollary \ref{corol_loctet})).
  {\it
    1) Any linkage diagram containing non-zero $\alpha$-label
    belongs to one of the loctets of the linkage system.

    2) Any linkage diagram of the loctet uniquely determines the whole loctet.

    3) If two loctets have one common linkage diagram they coincide.

    4) Every linkage diagram from the linkage system either belongs to one of the loctets
    or is $\beta$-unicolored.
   }

   \begin{remark} {\rm 1)
     \label{rem_l_less_8}
    The Carter diagrams $E_8(a_i)$ for $i = 1,\dots,4$ and $E_8$  do not represent any conjugacy classes in  $W(D_n)$:

     (a) For Carter diagrams from $\mathsf{C4}$, only  $D_l(a_k)$ (for some $l, k$) are conjugacy classes in  $W(D_n)$,
     see \cite[p. 13]{Ca72}, \cite[\S3]{St10}, i.e., $E_8(a_i)$ can not be mapped into  $W(D_n)$.

     (b) For Carter diagrams from $\mathsf{DE4}$, only  $D_l$ (for some $l$) are conjugacy classes in  $W(D_n)$,
     see Proposition \ref{prop_E6_Dk}, i.e., $E_i$ can not be mapped into  $W(D_n)$.

    Hence, there is no such a Weyl group containing linearly independent $9$-element root subset that embraces the root subset associated with $E_8(a_i)$ or $E_8$, i.e., \underline{there are no linkages} for any Carter diagram of type
    $E_8(a_i)$, where $i = 1,\dots,4$, and $E_8$, see \cite[\S2.3.3 and Table 2.3]{St10}.  For this reason, among simply-laced Carter diagrams of $E$-type we can consider only diagrams with a number of vertices $l < 8$.

   2) For the union of two Carter diagrams $\Gamma = \Gamma_1 \coprod \Gamma_2$ we have
      the direct sums $L = L_1 \oplus L_2$ and $L^{\vee} = L_1^{\vee} \oplus L_2^{\vee}$, and the linkage diagram $\gamma^{\vee} \in L^{\vee}$ splits into two summands $\gamma^{\vee} = \gamma^{\vee}_1 + \gamma^{\vee}_2$, where $\gamma^{\vee}_i \in L_i^{\vee}$. Accordingly, the partial Cartan matrix $B_L$, and its inverse matrix  $B_L^{-1}$ are decomposed into a direct sum, and dual partial Weyl group $W_L^{\vee} = W_{L_1}^{\vee} \times W_{L_2}^{\vee}$. Thus, the linkage system of $\Gamma$ is the direct product of linkage systems for $\Gamma_1$ and  $\Gamma_2$.
      \qed
   }
   \end{remark}

  In Section \ref{sec_calc_gamma_8} the calculation technique for loctet diagrams  $\gamma^{\vee}_{ij}(8)$ (= $8$th linkage diagram of the loctet) is explained. According to Corollary \ref{corol_loctet}, the whole loctet is uniquely determined
  from $\gamma^{\vee}_{ij}(8)$. By Tables \ref{sol_inequal_1}-\ref{sol_inequal_6} one can recover the calculation of $\gamma^{\vee}_{ij}(8)$ for Carter diagrams $\Gamma \in \mathsf{C4} \coprod \mathsf{DE4}$. 
  Similarly, in Section \ref{sec_calc_homog}
  the calculation technique for $\beta$-unicolored linkage diagrams is explained.
  By Tables \ref{homog_inequal_1}-\ref{homog_inequal_3} one can recover the calculation of $\beta$-unicolored linkage diagrams for Carter diagrams $\Gamma \in \mathsf{C4} \coprod \mathsf{DE4}$.
  In Section \ref{sec_diagr_per_comp} the loctets per components
  for all linkage systems are listed in Table \ref{tab_seed_linkages_6}.
  In Section \ref{sec_inv_matr} the partial Cartan matrix $B_L$ and the matrix $B^{-1}_L$, the inverse of the partial Cartan matrix $B_L$ for all simply-laced Carter diagrams containing $4$-cycle and Dynkin diagrams $E_6$, $E_7$ are listed in Tables \ref{tab_partial Cartan_1}-\ref{tab_partial Cartan_3}.
  In Section \ref{sect_linkage_diagr} all linkage systems are depicted in Fig. \ref{D4a1_linkages}-\ref{Dk_al_wind_rose}. The description of all linkage systems is presented
  in Theorem \ref{th_full_descr}.

\subsubsection{Linkage systems and weight systems}

  The linkage and weight systems for $E_6$ coincide,
  see Fig. \ref{E6pure_loctets_norm_gam8} and Fig. \ref{27_weight_diagr_E6__2comp}.
  It becomes obvious after recognizing loctets in both diagrams. The comparative figure containing
  both the linkage systems and the weight systems together with all their loctets can be seen in Fig. \ref{27_weights}.
  Similarly, the linkage and weight systems for $E_7$ coincide,
  see Fig. \ref{E7pure_linkage_system} and Fig. \ref{56_weight_diagr_E7}.

\begin{remark}
  \label{rem_weight_system}
  {\rm
  1)  Let us $\overline{\omega_l}$ be highest weight corresponding to the dimension of the smallest fundamental representation of semisimple Lie algebra associated with the Dynkin diagram $\Gamma$,  see \cite[Ch. 8, Table 2]{Bo05}, \cite[Ch. 6.1.10]{Bo02}  and $\gamma_l$ be the linkage connected only to the
  endpoint vertex $l$. This is not a surprise that the linkage system for $\Gamma$ and the weight system for  $\overline{\omega_l}$ coincide. The reason for this fact is the following: the dual partial Weyl group  $W_L$  for the Dynkin diagram of CCl (see footnote in Section \ref{sec_linkage_diagr}) coincides with the Weyl group $W$ of this type, the Dynkin labels of $\overline{\omega_l}$ coincide with linkage labels $\gamma^{\vee}_l$, so we have that
  \begin{equation}
    \label{start_vector}
     \gamma^{\vee}_l = B_L{\gamma_l} = {\bf B}\overline{\omega_l} =
       \left ( \begin{array}{c}
                1\\
                0\\
                \dots\\
                0 \\
              \end{array}
              \right ).
  \end{equation}
  Note that the vector \eqref{start_vector} is the vector from which the group $W^{\vee} = W_L^{\vee}$ starts to act (up to permutations of the top and bottom of the figure and permutations of appropriate coordinates), see Fig. \ref{E6pure_loctets_norm_gam8}, Fig. \ref{27_weight_diagr_E6__2comp},
  (resp. Fig. \ref{E7pure_linkage_system}, Fig. \ref{56_weight_diagr_E7}) for $E_6$ (resp. $E_7$).
  Further, the orbit of the vector ${\bf B}\overline{\omega_l}$
  (under action of the dual Weyl group $W^{\vee}$)  is the weight system,
  and the orbit of $\gamma^{\vee}_l$
  (under action of the partial dual Weyl group $W^{\vee}_L$, which coincides with $W^{\vee}$ in this case)
  is the linkage system.

  2) The weight system corresponding to the highest weight $\overline{\omega_1}$ for type $D_l$
     is taken from \cite[Fig. 4]{PSV98}, see Fig. \ref{Dl_weight_system}.  The number of weights is
     $2l$.  According to heading 3),  the linkage system has the same diagram. Compare with
     the linkage systems for $D_l(a_k)$, see Fig. \ref{Dk_al_linkages}, \ref{Dk_al_wind_rose},
     \ref{D5pure_loctets}, \ref{D6pure_loctets}.
\begin{figure}[H]
\centering
\includegraphics[scale=1.0]{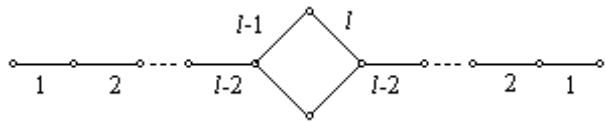}
\vspace{3mm}
\caption{\hspace{3mm}The weight system $(D_l, \overline{\omega_1})$.}
\label{Dl_weight_system}
\end{figure}
   Note that the linkage system $D_l$ (= weight system $(D_l, \overline{\omega_1})$) in Fig. \ref{Dl_weight_system} has exactly the same shape as the Carter diagram $D_l(a_k)$ despite the
   fact that the vertices in these diagrams are of the different nature.

  3)  There are two non-conjugate conjugacy classes $A'_l$ and $A{''}_l$  in $W(D_n)$,
 see Remark \ref{two_class_Al}.
 Nevertheless, the linkage system for $A_l$ can be constructed. According to heading 3),
 this linkage system coincides with the weight system for $A_l$. Figures of the weight system
 $(A_l, \overline{\omega_1})$ can be found in \cite[Fig. 1]{PSV98}. The shape of this weight
 system coincides with the Dynkin diagram $A_{l+1}$.
  } \qed
\end{remark}

\begin{theorem}
  \label{th_full_descr}
     In Table \ref{val_Bmin1},
     the values $\mathscr{B}^{\vee}_L(\gamma^{\vee})$, the number of components and number of
     linkage diagrams are collected for Carter diagrams from $\mathsf{C4}$ and $\mathsf{DE4}$.
 \end{theorem}

  \PerfProof The number of linkage diagrams is obtained from the enumeration of loctets and
  $\beta$-unicolored linkage diagrams for every Carter diagram from Table \ref{val_Bmin1}, see Section \ref{sec_enum_loctets}. The number of components is obtained from the shape of linkage systems.
  For the Carter diagrams from $\mathsf{C4}$, see Fig. \ref{D5a1_linkages} - \ref{E7a4_linkages},
  Fig. \ref{D7a1_linkages_cdef} - \ref{D7a2_linkages_comp2}, Fig. \ref{D7a1_D7a2_D7pu_loctets_comp3}.
  For the Carter diagrams from $\mathsf{DE4}$,
  see Fig. \ref{E6pure_loctets_norm_gam8},  Fig. \ref{E7pure_linkage_system}, Fig. \ref{D5pure_loctets},
  Fig. \ref{D6pure_loctets},  Fig. \ref{D7pu_loctets_comp1}-\ref{D7a1_D7a2_D7pu_loctets_comp3}.
  For $D_l(a_k)$ and $D_l$, where $l \geq 8$, the statement is proved in Proposition \ref{prop_Dl_ak}.  \qed

  \begin{table}[H]
  \centering
  \renewcommand{\arraystretch}{1.5}
  \begin{tabular} {||c|c|c|c|c|c||}
  \hline \hline
      The Carter &  Number of  &
      \multicolumn{2}{c|}{$E$-type components}   & \multicolumn{2}{c||}{$D$-type components} \cr
      \cline{3-6}
        diagram  &  components   &
      $\mathscr{B}^{\vee}_L(\gamma^{\vee})$ &  Number of &
      $\mathscr{B}^{\vee}_L(\gamma^{\vee})$ &  Number of  \cr
        &   &  &  linkage diagrams & & linkage diagrams \\
    \hline \hline
       $D_4(a_1)$ &   $3$  & $1$   & $16 (=8\times2)$ & $1$ & $8$ \\
    \hline
       $D_4$      &   $3$  & $1$   & $16 (=8\times2)$ & $1$ & $8$ \\
    \hline \hline
       $D_5(a_1)$ &   $3$  & $\frac{5}{4}$   & $32 (=16\times2)$  & $1$ & $10$ \\
    \hline
       $D_5$      &   $3$  & $\frac{5}{4}$   & $32 (=16\times2)$  & $1$ & $10$ \\
    \hline \hline
       $E_6(a_1)$ &   $2$  &  $\frac{4}{3}$  & $54 (=27\times2)$  & - & -  \\
    \hline
       $E_6(a_2)$ &   $2$  &   $\frac{4}{3}$  & $54 (=27\times2)$ &  - & - \\
    \hline
       $E_6$      &   $2$  &   $\frac{4}{3}$  & $54 (=27\times2)$ &  - & - \\
    \hline \hline
       $D_6(a_1)$ &   $3$  &  $\frac{3}{2}$  & $64 (=32\times2)$  & $1$ &  $12$ \\
    \hline
       $D_6(a_2)$ &   $3$  &  $\frac{3}{2}$ & $64 (=32\times2)$ & $1$ & $12$ \\
    \hline
       $D_6$      &   $3$  &  $\frac{3}{2}$ & $64 (=32\times2)$ & $1$ & $12$ \\
    \hline \hline
       $E_7(a_1)$ &   $1$  &   $\frac{3}{2}$  & $56$ &  - & - \\
    \hline
       $E_7(a_2)$ &   $1$  &   $\frac{3}{2}$  & $56$  &  - & - \\
    \hline
       $E_7(a_3)$ &   $1$  &  $\frac{3}{2}$  &  $56$ & - &  - \\
    \hline
       $E_7(a_4)$ &   $1$  &  $\frac{3}{2}$  &  $56$ & - &  - \\
    \hline
       $E_7$      &   $1$  &  $\frac{3}{2}$  &  $56$ & - &  - \\
    \hline \hline
       $D_7(a_1)$ &   $3$  &   $\frac{7}{4}$  & $128 (=64\times2)$ &  $1$  & $14$ \\
    \hline
       $D_7(a_2)$ &   $3$  &   $\frac{7}{4}$  & $128 (=64\times2)$ &  $1$ &  $14$  \\
    \hline
       $D_7$      &   $3$  &   $\frac{7}{4}$  & $128 (=64\times2)$  &  $1$ &  $14$  \\
    \hline \hline
       $D_l(a_k)$, $l > 7$ &   $1$  &   -  & - & $1$  & $2l$ \\
    \hline
       $D_l$, $l >  7$ &   $1$  &   -  & - & $1$  & $2l$ \\
   \hline  \hline  
\end{tabular}
  \vspace{2mm}
  \caption{\hspace{3mm}Values of $\mathscr{B}^{\vee}_L(\gamma^{\vee})$ and the number of
   linkage diagrams for the Carter diagrams.}
  \label{val_Bmin1}
  \end{table}

\begin{remark}[Additions to Table \ref{val_Bmin1}]
  \label{rem_DE_except}
{\rm
   1) For $D_4(a_1)$, two $E$-type components of the linkage system
   are also $D$-type components, see Fig. \ref{D4a1_linkages}, components (II) and (III).
   This is reflected in fact that
   $\mathscr{B}^{\vee}_L(\gamma^{\vee}) = 1$ for all $3$ components.
   The third component is only of the $D$-type, see Fig. \ref{D4a1_linkages}, component (I).

   2) For $D_4$, there are $3$ components in the linkage system, each of which is
   the $D$-type component and also $E$-type component.
   To obtain these components
   one can take only $4$ first coordinates for any linkage diagrams in Fig. \ref{fig_loctets},
   see Fig. \ref{D4_loctets}.
   Every component is exactly the \underline{loctet}. These components coincide with $3$
   weight systems of $3$ fundamental representations of semisimple Lie algebra $D_4$:
    $(D_4, \overline{\omega}_1)$, $(D_4, \overline{\omega}_3)$, $(D_4, \overline{\omega}_4)$,
    see \cite[Fig. 10]{PSV98}.
}
\end{remark}

\begin{figure}[H]
\centering
\includegraphics[scale=1.0]{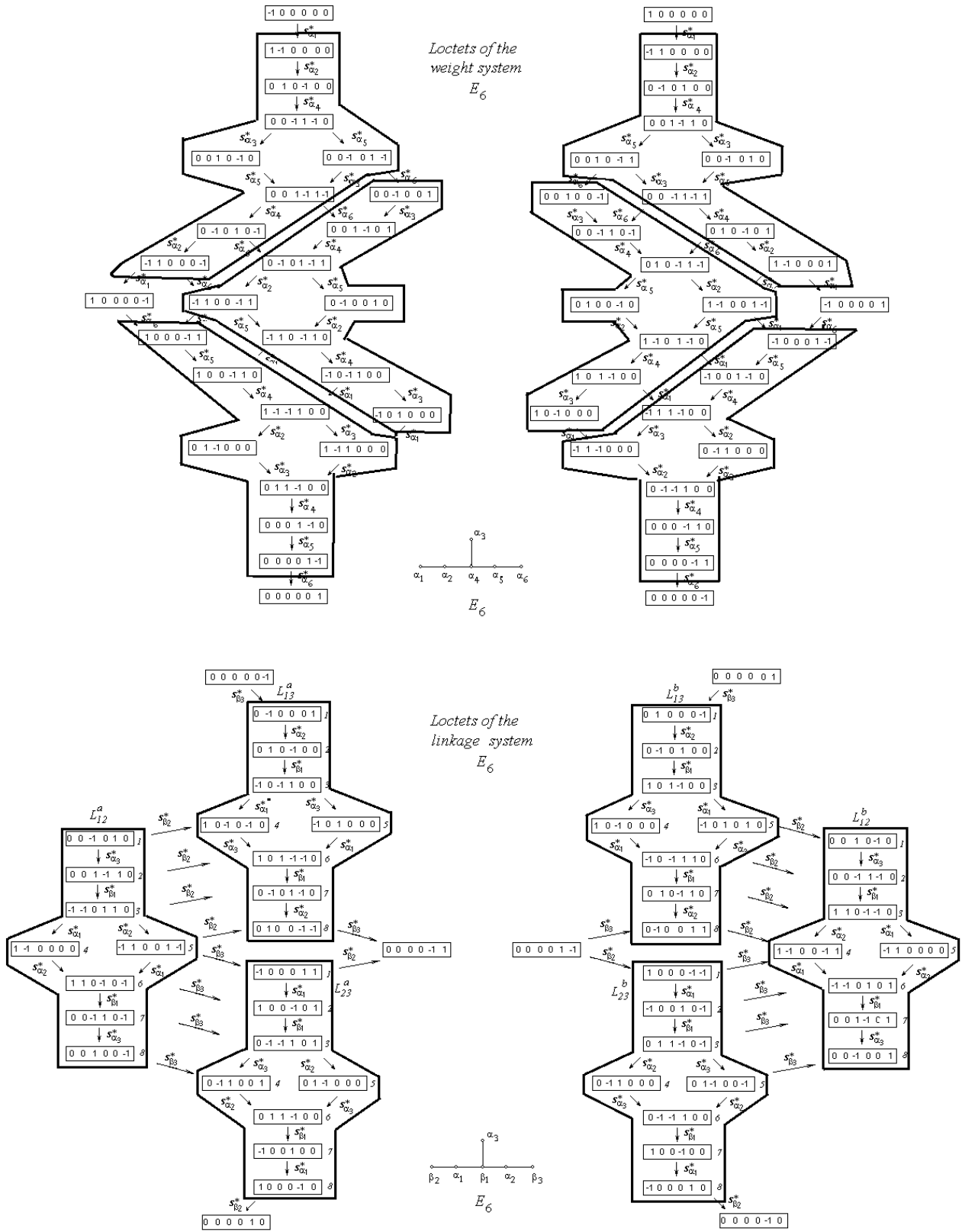}
\caption{\hspace{3mm}Loctets in the weight system and in the linkage system $E_6$}
\label{27_weights}
\end{figure}

 \subsubsection{Projection of linkage diagrams and loctets}

  Section \ref{sec_proj} is devoted to the projection of linkage systems, see Fig. \ref{inclusions_E6_E7},
  Fig. \ref{inclusions_D6a_D5a1}.
  Primarily, we consider the simply extendable Carter diagrams in Section \ref{sect_extendable}.
  The Carter diagram $\Gamma$ is called {\it simply extendable in the vertex $\tau_p$} if
  the new diagram obtained by the extra vertex $\tau_{l+1}$ together with the additional connection edge
  $\{ \tau_p, \tau_{l+1} \}$  is also the Carter diagram.
  In the case, where the Carter diagram $\widetilde\Gamma$ is an
  simple extension of other Carter diagram $\Gamma$, we construct the projection of
  the linkage diagrams of $\widetilde\Gamma$ to the linkage diagrams of $\Gamma$.
  The main result of this section is the following theorem:

\begin{figure}[H]
\centering
\includegraphics[scale=0.6]{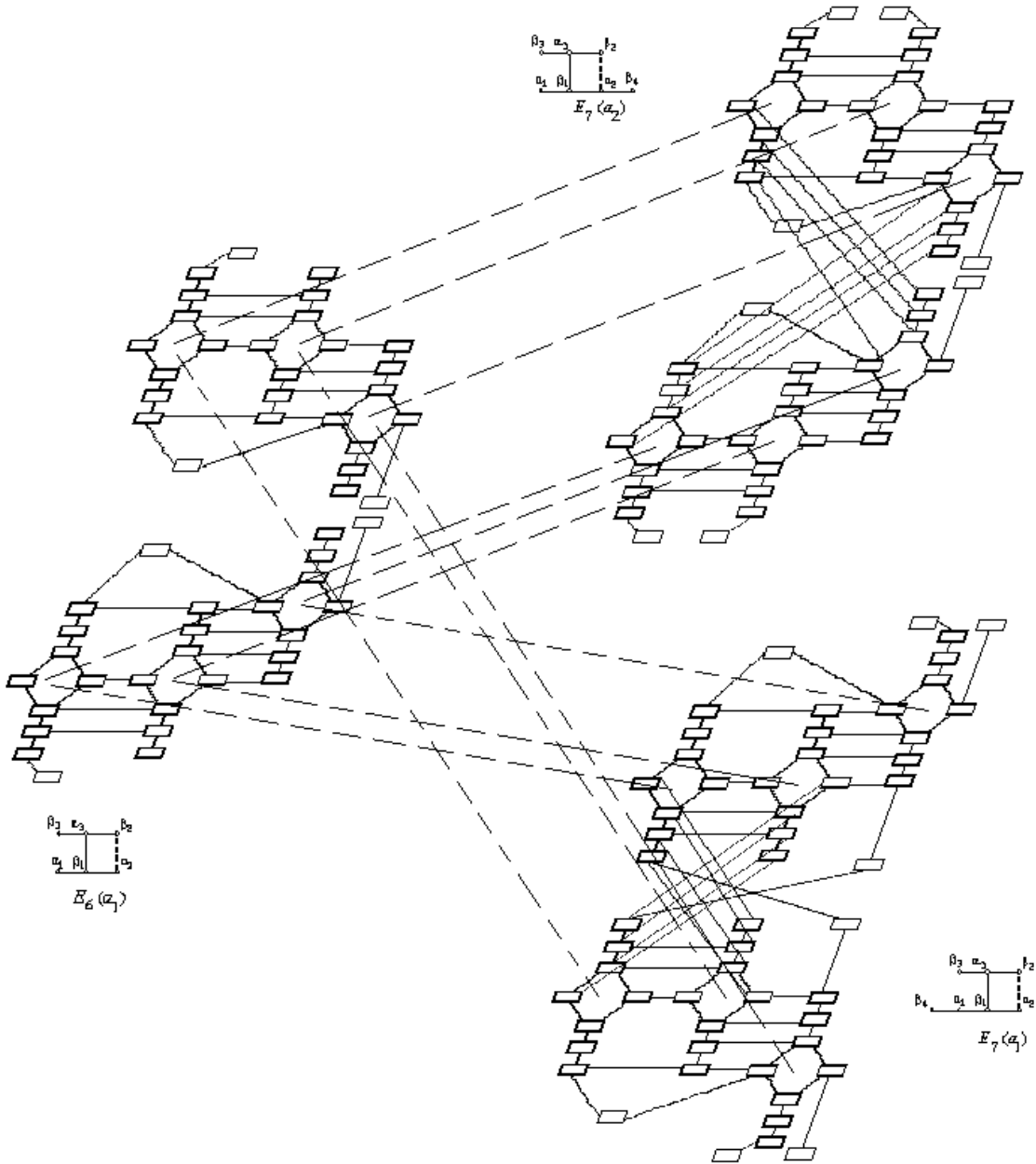}
\vspace{3mm}
\caption[\hspace{3mm}Projection of linkage systems $E_7(a_1)$, $E_7(a_2)$ onto the linkage system $E_6(a_1)$]{\hspace{3mm}Projection of linkage systems $E_7(a_1)$, $E_7(a_2)$ onto the linkage system $E_6(a_1)$
 (depicted in Fig. \ref{E6a1_linkages}, Fig. \ref{E7a1_linkages}, Fig. \ref{E7a2_linkages}).}
\label{inclusions_E6_E7}
\end{figure}

  {\bf Theorem} (Theorem \ref{th_map_linkage_systems}).
  {\it  Let the Carter diagram $\widetilde\Gamma$ be the simple extension of the Carter diagram $\Gamma$ in the vertex $\tau_p$, such that the vertex $\tau_p$ is connected to the vertex $\tau_{l+1}$.
  Let $\widetilde\gamma^{\vee}$ be a certain linkage diagram for $\widetilde\Gamma$, and $\gamma^{\vee}$
  be the vector obtained from $\widetilde\gamma^{\vee}$ by removing the coordinate $\tau_{l+1}$.
  If $\gamma^{\vee} \neq 0$ then $\gamma^{\vee}$ is the linkage diagram for $\Gamma$.}

  This theorem is used for describing linkage systems $D_l(a_k)$, $D_l$, where $l > 7$, in Proposition \ref{prop_Dl_ak}.
  According to Remark \ref{rem_mapping}, we obtain the projection of loctets associated with Carter
  diagrams $\widetilde\Gamma$ and $\Gamma$. This projection is presented by means of dotted lines connecting
  centers of loctets of $\Gamma$ and $\widetilde\Gamma$, see Fig. \ref{inclusions_E6_E7}.

~\\
{\bf Dedication.} This paper is dedicated to the memory {\it Gurii Fedorovich Kushner}, \cite{K70, K72, K79}, who in the early $80$'s first told me about the highest weight representations of semisimple Lie algebras.

%% file: 2a_more4cycles.tex
\newpage
\section{\sc\bf The Cartan matrix associated with conjugacy class}
  \label{sec_Cartan_matr_ccl}

\subsection{More on $4$-cycles}

\subsubsection{How many endpoints may be in a linkage diagram?}

 In that follows, we show that the number of endpoints in any linkage diagram is not more
 than $6$,  and in some cases this number is not more than $4$.

\begin{proposition}
   \label{prop_numb_ep}
    Let $w = w_\alpha w_\beta$ be the bicolored decomposition
    of $w$ into the product of two involutions, and $\Gamma$ be the Carter
    diagram corresponding to this decomposition.
    Let $\gamma^{\vee}$ be any linkage diagram obtained from $\Gamma$.

    1) The linkage diagram $\gamma^{\vee}$ does not have more than $3$ $\alpha$-endpoints and $3$ $\beta$-endpoints.
~\\
\begin{minipage}{12.5cm}
    \quad 2) If $\alpha$-set (resp. $\beta$-set) contains $3$ points
    connecting to one branch point $\beta$ (resp. $\alpha$),
    i.e., $\{ \beta \cup \alpha\text{-set} \}$ (resp. $\{ \alpha \cup \beta\text{-set} \}$) forms
    the diagram $D_4$ then there are no more
    than two $\alpha$-endpoints (resp. $\beta$-endpoints) in the
    given $\alpha$-set (resp. $\beta$-set). (Corollary 2.4, from \cite{St10}.)
\end{minipage}
\begin{minipage}{4.5cm}
  \qquad \psfig{file=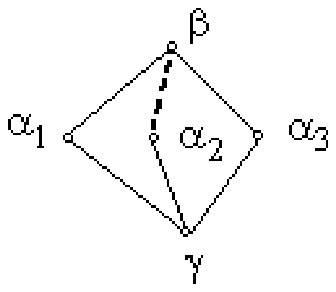, width=1.0in}
\end{minipage}
~\\

    3) Let $\{\alpha_1, \beta_1, \alpha_2, \beta_2 \}$ be
    the square in certain connected diagram. There does not exist
    a root $\gamma$ connected to all vertices of the square. (Corollary 2.4, from \cite{St10}.)
 \end{proposition}

\PerfProof
   1) If the linkage $\gamma^{\vee}$ has $4$ $\alpha$-endpoints then the connection diagram
      contains the diagram
      $\widetilde{D}_4 = \{ \gamma, \alpha_1, \alpha_2, \alpha_3, \alpha_4 \}$,
      then the vector\footnotemark[1]
   \footnotetext[1]{We denote by the same letters vertices and the corresponding roots.}
      $\varphi = 2\gamma + \sum\limits_{i=1}^4\alpha_i$
      has zero length, since
\begin{equation}
  (\varphi, \varphi) = 4(\gamma, \gamma) + \sum\limits_{i=1}^4(\alpha_i, \alpha_i)  +
        4\sum\limits_{i=1}^4(\gamma, \alpha_i) = 4\cdot{2} + 4\cdot{2} - 16\cdot{1} = 0.
\end{equation}
      Hence, $\varphi = 0$, contradicting the linear independence of roots
      $\{ \gamma, \alpha_1, \alpha_2, \alpha_3, \alpha_4 \}$.
 \qed

 \begin{corollary}
       a) For the Carter diagrams $D_5(a_1)$, $E_6(a_1)$,  $E_6(a_2)$, $D_6(a_2)$, any linkage diagram
        does not contain more than two endpoints in $\alpha$-set (resp. $\beta$-set).

       b) For the Carter diagrams $D_6(a_1)$, $E_7(a_1)$, $E_7(a_2)$, $E_7(a_3)$, $E_7(a_4)$,
          any linkage diagram does not contain more than two endpoints in $\alpha$-set.

       c) For the Carter diagrams $D_7(a_1)$, $D_7(a_2)$, there are linkage diagrams containing
          more than two endpoints in $\alpha$-set or in $\beta$-set.
 \end{corollary}
   See $\alpha$-set and $\beta$-set of Carter diagrams in
   Tables \ref{tab_partial Cartan_1}-\ref{tab_partial Cartan_2}.
\qed

\subsubsection{The diagonal in a square}

\begin{figure}[h]
\centering
\includegraphics[scale=0.8]{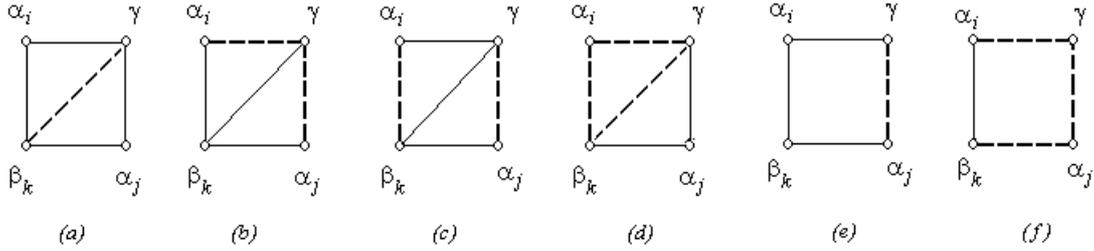}
 \caption{\hspace{3mm}Linkage diagrams containing a square.}
\label{square_n_diag}
\end{figure}

\begin{proposition}[On squares]
  \label{prop_diagonal}
   Let $\gamma$ form the linkage diagram containing
   the square $\{ \alpha_i, \beta_k, \alpha_j, \gamma \}$
   without the diagonal $\{ \alpha_i, \alpha_j \}$, i.e.,
   the roots $\{ \alpha_i, \beta_k, \alpha_j, \gamma \}$ are linearly independent, 
 \begin{equation*}
   (\alpha_i, \beta_k) \neq 0, \quad  (\beta_k, \alpha_j) \neq 0, \quad
   (\alpha_j, \gamma) \neq 0, \quad  (\gamma, \alpha_i) \neq 0, \quad \text{ and } \quad  (\alpha_i, \alpha_j) = 0.
 \end{equation*}
   If there is an even number of dotted edges in the square then
   there exists a diagonal in the square.
   If there is an odd number of dotted edges in the square then
   there is no any diagonal in the square.
   Namely:

    {\rm (a)} If there is no dotted edge in the square
    then \underline{there exists the dotted diagonal
   $\{ \gamma, \beta_k \}$}, i.e., $(\gamma, \beta_k) = 1$,
   see Fig. \ref{square_n_diag},(a).

   {\rm (b)} If there are two dotted edges $\{ \gamma, \alpha_i \}$ and
   $\{ \gamma, \alpha_j \}$ in the square, i.e.,
   $(\gamma, \alpha_i) = (\gamma, \alpha_j) = 1$,
   and remaining edges are  solid then \underline{there exists the solid diagonal
   $\{ \gamma, \beta_k \}$}, i.e., $(\gamma, \beta_k) = -1$,
   see Fig. \ref{square_n_diag},(b).

   {\rm (c)} If there are two dotted edges $\{ \gamma, \alpha_j \}$ and
   $\{ \beta_k, \alpha_i \}$, i.e.,
   $(\gamma, \alpha_j) = (\beta_k, \alpha_i) = 1$,
   and remaining edges are  solid then \underline{there exists the solid diagonal
   $\{ \gamma, \beta_k \}$}, i.e., $(\gamma, \beta_k) = -1$,
   see Fig. \ref{square_n_diag},(c).

   {\rm (d)} If there are two dotted edges $\{ \gamma, \alpha_i \}$ and
   $\{ \beta_k, \alpha_i \}$, i.e.,
   $(\gamma, \alpha_i) = (\beta_k, \alpha_i) = 1$,
   and remaining edges are solid then \underline{there exists the dotted diagonal
   $\{ \gamma, \beta_k \}$}, i.e., $(\gamma, \beta_k) = -1$,
   see Fig. \ref{square_n_diag},(d).

   {\rm (e)} If there is only one dotted edge $\{ \gamma, \alpha_j \}$, i.e.,
   $(\gamma, \alpha_j) = 1$,  and remaining edges are
   solid then \underline{there is no any diagonal}, i.e., $(\gamma, \beta_k) = 0$,
   see Fig. \ref{square_n_diag},(d).

   {\rm (f)} If there are three dotted edges
   $\{ \gamma, \alpha_i \}, \{ \gamma, \alpha_j \}, \{ \beta_k, \alpha_j \}$, i.e.,
   $(\gamma, \alpha_i) = (\gamma, \alpha_j) = (\beta_k, \alpha_j) = 1$,
   and remaining edge is
   solid then \underline{there is no any diagonal}, i.e., $(\gamma, \beta_k) = 0$,
   see Fig. \ref{square_n_diag},(e).
\end{proposition}

  \PerfProof If there is no diagonal $\{ \gamma, \beta_k \}$ for one of cases (a),(b), (c) or (d),
   we get the extended Dynkin diagram $\widetilde{A}_4$ by following changes:
  \begin{equation*}
     (a) \quad \text{no changes,}
     \qquad (b) \quad \gamma \longrightarrow -\gamma, 
     \qquad (c) \quad \gamma \longrightarrow -\gamma, \quad  \alpha_j \longrightarrow  -\alpha_j, 
     \qquad (d) \quad \alpha_j \longrightarrow  -\alpha_j, \\
  \end{equation*}
  contradicting Lemma A.1 from \cite{St10}.
  If the diagonal (dotted or solid) exists in cases (e) or (f)
  then one of obtained triangles (if necessary, after the change $\gamma \longrightarrow -\gamma$) is
   the extended Dynkin diagram $\widetilde{A}_3$, contradicting Lemma A.1 from \cite{St10}.
   \qed


%% file: 3dualWeyl.tex
 \subsection{The partial Cartan matrix}
    \label{sec_partial_B}
 Let $L \subset E$ (resp. $L' \subset E$) be the linear subspace spanned by root subset\footnotemark[1] $S$ (resp. $S'$), see  \eqref{canon_dec_1}.
 \footnotetext[1]{ The root subset $S$ (resp. $S'$) is not the root subsystem since
 roots of $S$ (resp. $S'$) are not necessarily simple.}
 From now on, if the bicolored decomposition does not matter, we use the notation $\tau_i$ for the roots
 instead of the bicolored notation $\alpha_i$ and $\beta_j$, i.e.,
 \begin{equation}
   \label{tau_noation}
   \begin{split}
    & S = \{ \tau_1, \dots,  \tau_{k+h} \},
    \text{ where } \tau_i = \alpha_i \text{ for } i = 1, \dots, k \text{ and } \tau_{j+k} = \beta_j
    \text{ for } j = 1, \dots, h, \\
    & S' = \{ \tau'_1, \dots, \tau'_{k+h} \},
    \text{ where } \tau'_i = \alpha'_i \text{ for } i = 1, \dots, k \text{ and } \tau'_{j+k} = \beta'_j
    \text{ for } j = 1, \dots, h. \\
   \end{split}
 \end{equation}

  Similarly to the Cartan matrix associated with Dynkin diagrams, we determine the Cartan matrix for each Carter diagram $\Gamma$ from $\mathsf{C4}$ or $\mathsf{DE4}$ as follows

 \begin{equation}
   \label{canon_dec_2}
   B_L :=
      \left (
        \begin{array}{ccc}
         (\alpha_i, \alpha_j) & \dots & (\alpha_i, \beta_r)  \\
         \dots                & \dots & \dots \\
         (\beta_q, \alpha_j) & \dots & (\beta_q, \beta_r )  \\
        \end{array}
      \right ) \quad
        \begin{array}{c}
            i,j = 1, \dots, k, \\
            \dots \\
            q,r = 1, \dots, h.
         \end{array}
 \end{equation}
 We call this matrix the {\it partial Cartan matrix}. According to \eqref{uniq_B},
 $B_L = B_{L'}$, i.e., the partial Cartan matrix $B_L$ is well-defined. Thus,
 for every Carter diagram $\Gamma \in \mathsf{C4} \coprod \mathsf{DE4}$, we associate the Cartan matrix also with the conjugacy
 class corresponding to $\Gamma$. The symmetric bilinear form associated with
 the partial Cartan matrix $B_L$ is denoted by $(\cdot, \cdot)_L$ and the
 corresponding quadratic form is denoted by $\mathscr{B}_L$.
 \begin{proposition}
   \label{restr_forms_coincide}
 1)  The restriction of the bilinear form associated with the Cartan matrix
  ${\bf B}$ on the subspace $L$ coincide with the bilinear form associated with the
  partial  Cartan matrix $B_L$, i.e., for any pair of vectors $v, u \in L$ we have
 \begin{equation}
   \label{restr_q}
       (v, u)_{\botL} = (v, u), \text{ and }
       \mathscr{B}_L(v) = \mathscr{B}(v).
 \end{equation}

 2) For every Carter diagram, the matrix $B_L$ is positive definite.

 \end{proposition}
\PerfProof
 1) From \eqref{canon_dec_2} we deduce:
 \begin{equation*}
     (v, u)_{\botL} =   (\sum\limits_i{t_i{\tau_i}}, \sum\limits_j{q_j{\tau_j}})_{\botL} =
      \sum\limits_{i,j}t_i{q}_j(\tau_i, \tau_j)_{\botL} = \sum\limits_{i,j}t_i{q}_j(\tau_i, \tau_j) =
      (v, u).
 \end{equation*}

 2) This follows from 1).
\qed

\begin{remark}[The classical case]
   \label{rem_classic}
 \rm{
  Recall that the $n\times{n}$ matrix $K$ satisfying the following properties
 \begin{equation*}
   \begin{split}
     (C1) & \quad k_{ii} = 2 \text{ for } i = 1,\dots, n, \\
     (C2) & \quad -k_{ij} \in \mathbb{Z} = \{0, 1, 2, \dots \} \text{ for } i \neq j, \\
     (C3) & \quad  k_{ij} = 0 \text{ implies } k_{ji} = 0 \text{ for } i, j = 1, \dots, n
    \end{split}
 \end{equation*}
 is called a {\it generalized Cartan matrix}, \cite{Kac80}, \cite[\S 2.1]{St08}.

  The condition (C2) is not valid for the partial Cartan matrix:
 a few values $k_{ij}$ associated with dotted edges are positive,
 see Tables \ref{tab_partial Cartan_1}, \ref{tab_partial Cartan_2}.

 If the Carter diagram does not contain any cycle,
 then the Carter diagram is the Dynkin diagram, the corresponding conjugacy class
 is the conjugacy class of the Coxeter element, and
 the partial Cartan matrix is the classical Cartan matrix, which is
 the subclass of generalized Cartan matrices.
 }
\end{remark}

 \begin{remark}[The predefined numbering]
  \label{two_predef_classes}
 \rm{
   Every Carter diagrams $\Gamma \in \mathsf{C4}$ with number of vertices $l \geq 5$ contains
   the subdiagram $D_5(a_1)$, see Tables \ref{tab_partial Cartan_1}, \ref{tab_partial Cartan_2}.
   We predefine the vertex numbering of the subdiagram $D_5(a_1)$ in any encompassing diagram
   as in Fig \ref{D5a1_D4_pattern}.
 \begin{figure}[H]
\centering
\includegraphics[scale=0.9]{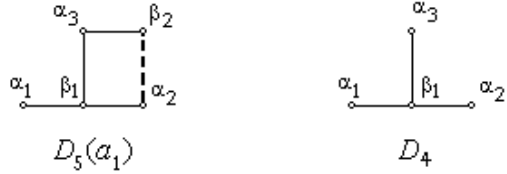}
\vspace{3mm}
\caption{\hspace{3mm}The predefined numbering in two pattern subdiagrams: $D_5(a_1)$ and $D_4$}
\label{D5a1_D4_pattern}
\end{figure}
   For Carter diagrams $\Gamma \in \mathsf{DE4}$, i.e., Dynkin diagrams $D_l$, where $l \geq 4$,
   $E_6$, $E_7$, $E_8$, the diagram
   $\Gamma$ contains subdiagram $D_4$. We predefine the vertex numbering of
   the subdiagram $D_4$ in any encompassing diagram as in Fig \ref{D5a1_D4_pattern}.
   The predefined numbering of vertices Fig. \ref{D5a1_D4_pattern} is presented, for example, in
   Fig. \ref{E7a1_gamma_ij_3}-\ref{E7a3_gamma_ij_4}.
 \qed
}
\end{remark}

 \subsection{The dual partial Weyl group associated with a conjugacy class}
   Let  $\Gamma \in \mathsf{C4} \coprod \mathsf{DE4}$, and $\Pi_w$ be
   the root basis associated with $\Gamma$. Up to conjugacy, the root basis $\Pi_w$
   is well-defined for any $\Gamma \in \mathsf{C4} \coprod \mathsf{DE4}$.
   Indeed, let $w$ and $w'$ have the same Carter diagram and be conjugate,
   i.e., $w' = \tilde{w}w\tilde{w}^{-1}$, for some $\tilde{w} \in W$.
   For any bicolored decomposition
   \begin{equation}
      w = s_{\alpha_1} \dots s_{\alpha_k} s_{\beta_1} \dots s_{\beta_h},
   \end{equation}
   the conjugate element $w'$  has also a bicolored decomposition associated with the same diagram $\Gamma$:
  \begin{equation}
      w' = s_{\tilde{w}\alpha_1} \dots s_{\tilde{w}\alpha_k} s_{\tilde{w}\beta_1} \dots s_{\tilde{w}\beta_h}.
   \end{equation}
   ~\\
   Let us define dual vectors $\tau_i^{\vee} \in L^{\vee}$ for any  $\tau_i \in \Pi_w$:
   \begin{equation}
      \label{eq_tau_vee}
       \tau_i^{\vee} := B_L{\tau_i}.
   \end{equation}
   Eq. \eqref{eq_tau_vee} is consistent with \eqref{dual_gamma}. The mapping ${}^{\vee} : \tau_i \longrightarrow \tau_i^{\vee}$ is expanded to the linear mapping $L \longrightarrow  L^{\vee}$,
   and
   \begin{equation}
      \label{eq_gamma_vee}
       u^{\vee} = B_L{u} =
       \left (
       \begin{array}{c}
         (u, \tau_1) \\
          \dots \\
         (u, \tau_l) \\
       \end{array}
       \right ) \text{ for any } u \in L.
   \end{equation}
   ~\\
   Lengths of vectors $\tau_i^{\vee}$ in the sense of symmetric bilinear form
   $(\cdot , \cdot)_{\botLv}$ associated with $B^{-1}_L$ are equal to $2$, since
    \begin{equation}
       (\tau_i^{\vee}, \tau_i^{\vee})_{\botLv}  =
       \langle B^{-1}_L \tau_i^{\vee}, \tau_i^{\vee} \rangle =
       \langle \tau_i, B_L\tau_i \rangle = (\tau_i, \tau_i)_{\botL} = (\tau_i, \tau_i) = 2.
   \end{equation}
   Let us consider the restriction of the reflection $s_{\tau_i}$ on the subspace $L$.
   For any $v \in L$, by Proposition \ref{restr_forms_coincide} we have:
   \begin{equation}
     \label{refl_tau}
       s_{\tau_i}v = v  - 2\frac{(\tau_i, v)}{(\tau_i, \tau_i)}\tau_i =
                      v  - (\tau_i, v)_{\botL}\tau_i =
                      v  - \langle B_L\tau_i, v \rangle \tau_i =
                      v - \langle \tau_i^{\vee}, v \rangle \tau_i.
   \end{equation}
   We define the reflection $s^{*}_{\tau_i}$ acting on $u \in L^{\vee}$ as follows:
   \begin{equation}
     \label{refl_tau_dual}
       s^{*}_{\tau_i}u := u - 2\frac{(u, \tau_i^{\vee})_{\botLv}}{(\tau_i^{\vee}, \tau_i^{\vee})_{\botLv}}\tau_i^{\vee} = u - (u, \tau_i^{\vee})_{\botLv}\tau_i^{\vee} =
       u - \langle u, B^{-1}_L\tau_i^{\vee} \rangle \tau_i^{\vee} =
       u - \langle u, \tau_i \rangle \tau_i^{\vee}.
   \end{equation}
   Let $W_L$ (resp. $W_L^{\vee}$) be the group generated by reflections
   $\{s_{\tau_i} \mid \tau_i \in \Pi_w \}$
   (resp. $\{s^{*}_{\tau_i} \mid \tau_i \in \Pi_w \})$.

   \begin{proposition}
     1) For any $\tau_i \in \Pi_w$, we have
   \begin{equation}
      \label{refl_transp_2}
        s^{*}_{\tau_i} = {}^t{s}_{\tau_i} = {}^t{s}_{\tau_i}^{-1}.
   \end{equation}
     2) The mapping
   \begin{equation}
        w \rightarrow {}^t{w}^{-1}
   \end{equation}
     determines an isomorphism of $W_L$ onto $W_L^{\vee}$.
   \end{proposition}
   \PerfProof
  1) By \eqref{refl_tau} and \eqref{refl_tau_dual}, for any $v \in L, u \in L^{\vee}$ we have:
   \begin{equation}
     \label{refl_transpon}
     \begin{split}
      & \langle  s^{*}_{\tau_i}u, v \rangle =
       \langle  u - \langle u, \tau_i \rangle \tau_i^{\vee}, v \rangle =
        \langle u, v \rangle - \langle u, \tau_i \rangle \langle v, \tau_i^{\vee} \rangle, \\
      &   \langle  u, s_{\tau_i}v \rangle =
        \langle  u , v - \langle \tau_i^{\vee}, v \rangle \tau_i \rangle =
        \langle u, v \rangle - \langle \tau_i^{\vee}, v \rangle \langle u, \tau_i \rangle.
     \end{split}
   \end{equation}
   Thus,
   \begin{equation}
      \langle  s^{*}_{\tau_i}u, v \rangle = \langle  u, s_{\tau_i}v \rangle,
      \quad \text{ for any } v \in L, u \in L^{\vee},
   \end{equation}
   and \eqref{refl_transp_2} holds.
  \qed

  One should note that $W_L$ (and, therefore, $W_L^{\vee}$) is not necessarily Weyl group, since the roots $\tau_i \in \Pi_w$ are not necessarily simple and they do not constitute a root subsystem.
  We call $W_L$ the partial Weyl group, and $W_L^{\vee}$ the {\it dual partial Weyl group
  associated with a conjugacy class}, or, for short, the {\it dual partial Weyl group}.
    Then
    \begin{equation}
       \label{dual_refl}
       (s^*_{\tau_i}\gamma^{\vee})_{\tau_k} =
        \begin{cases}
            -\gamma^{\vee}_{\tau_i}, & \text{ for } k = i, \\
            \gamma^{\vee}_{\tau_k} + \gamma^{\vee}_{\tau_i}, & \text{ if } \{\tau_k, \tau_i \}
            \text{ is the {\it solid} edge, i.e., } (\tau_k, \tau_i) = -1, \\
            \gamma^{\vee}_{\tau_k} - \gamma^{\vee}_{\tau_i}, & \text{ if }  \{\tau_k, \tau_i \}
            \text{ is the {\it dotted} edge, i.e., } (\tau_k, \tau_i) = 1, \\
            \gamma^{\vee}_{\tau_i}, & \text{ if } \tau_k \text{ and } \tau_i
             \text{ are not connected, i.e. } (\tau_k, \tau_i) = 0.
        \end{cases}
    \end{equation}
\qed
~\\
Let $\Pi$ be the set of the simple roots in the root system $\varPhi$ associated with
the Weyl group $W$, let $E$ be the linear space spanned by simple roots from $\Pi$.
Let $M$ be the orthogonal complement of $L$ in $E$ in the sense of the bilinear form $(\cdot, \cdot)$:
\begin{equation}
    E = L \oplus M,  \quad M \perp L.
\end{equation}
Any root $\gamma \in \varPhi$ is uniquely  decomposed into the following sum:
\begin{equation}
  \label{decomp_mu}
    \gamma = \gamma_L +  \mu,
\end{equation}
where $\gamma_L \in L$,  $\mu \in M$.

\begin{proposition}
  \label{prop_link_diagr_conn}
    1) For the linkage labels vector $\gamma^{\vee}$ and the vector $\gamma_L$ from  \eqref{decomp_mu},
    we have
  \begin{equation}
     \label{conn_L_1}
    \gamma^{\vee} = B_L\gamma_L.
  \end{equation}
  2) For $s^{*}_{\tau_i}$, the dual reflection  \eqref{refl_tau_dual},
  the following relations hold
  \begin{equation}
      \label{conn_dual_1}
        B_L s_{\tau_i} = s^*_{\tau_i} B_L,
     \end{equation}
  \begin{equation}
     \label{conn_L_2}
    (s_{\tau_i}\gamma)^{\vee} = s^{*}_{\tau_i}B_L\gamma_L = s^{*}_{\tau_i}\gamma^{\vee}.
  \end{equation}
  3) For $w^* \in W^{\vee}_L$ $(w^* := {}^t{w}^{-1})$, the dual element of $w \in W$, we have
    \begin{equation}
      \label{conn_L_3}
        {(w\gamma)}^{\vee} = w^*{\gamma^\vee} \quad(= {}^t{w}^{-1}{\gamma^\vee}).
     \end{equation}
  4) The following relations hold    
    \begin{equation}
      \label{conn_L_5}
        \begin{split}
        & \mathscr{B}_L (s_{\tau_i} v)  =  \mathscr{B}_L (v),  \text{ for any } v \in L, \\
        & \mathscr{B}^{\vee}_L(s^*_{\tau_i} u) =
           \mathscr{B}^{\vee}_L(u), \text{ for any } u \in  L^{\vee}.
        \end{split}
     \end{equation}
     \end{proposition}

\PerfProof

1) Since $(\tau_i, \mu) = 0$ for any $\tau_i \in S$, by \eqref{eq_gamma_vee} we have
 \begin{equation*}
    \begin{split}
    \gamma^{\vee} = &
      \left (
      \begin{array}{c}
         (\gamma, \alpha_1) \\
         \dots \\
         (\gamma, \beta_h) \\
      \end{array}
      \right ) =
      \left (
      \begin{array}{c}
         (\gamma_L + \mu, \alpha_1) \\
         \dots \\
         (\gamma_L + \mu, \beta_h) \\
      \end{array}
      \right ) =
      \left (
      \begin{array}{c}
         (\gamma_L, \alpha_1) \\
         \dots \\
         (\gamma_L, \beta_h) \\
      \end{array}
      \right ) =  \gamma_L^{\vee} =
    B_L\gamma_L.
     \end{split}
  \end{equation*}

2)  The equality \eqref{conn_dual_1} holds since the following is true for
   any $u, v \in L$:
 \begin{equation*}
    \begin{split}
     & (s_{\tau_i}{u}, v)_L = ({u}, s_{\tau_i}{v})_L, \text{ i.e., }
     \langle B_L s_{\tau_i}{u}, v \rangle = \langle B_L {u}, s_{\tau_i}{v} \rangle =
     \langle s^*_{\tau_i}{B}_L {u}, v \rangle \text{, and } \\
     & \langle (B_L s_{\tau_i} - s^*_{\tau_i}{B}_L) {u}, v \rangle = 0.
    \end{split}
 \end{equation*}

Let us consider eq. \eqref{conn_L_2}.
Since $(\tau_i, \mu) = 0$ for any $\tau_i \in S$, and $s_{\tau_i}\mu = \mu$, we have
 \begin{equation*}
  \begin{split}
    (s_{\tau_i}\gamma)^{\vee} =
      & \left (
      \begin{array}{c}
         (s_{\tau_i}\gamma, \alpha_1) \\
         \dots \\
         (s_{\tau_i}\gamma, \beta_h) \\
      \end{array}
      \right ) =
      \left (
      \begin{array}{c}
         (s_{\tau_i}\gamma_L + \mu, \alpha_1) \\
         \dots \\
         (s_{\tau_i}\gamma_L + \mu, \beta_h) \\
      \end{array}
      \right ) =
      \left (
      \begin{array}{c}
         (s_{\tau_i}\gamma_L, \alpha_1) \\
         \dots \\
         (s_{\tau_i}\gamma_L, \beta_h) \\
      \end{array}
      \right ) =
      \\
    &  \\
    &
    (s_{\tau_i}\gamma_L)^{\vee} = B_L{s}_{\tau_i}\gamma_L =
    s^{*}_{\tau_i}B_L\gamma_L = s^{*}_{\tau_i}\gamma^{\vee}.
     \end{split}
  \end{equation*}

3)  Let $w = s_{\tau_1}s_{\tau_2}\dots{s}_{\tau_m}$ be
the decomposition of $w \in W$. Since  $s^*_{\tau} = {}^ts^{-1}_{\tau} = {}^ts_{\tau}$,
 we deduce from \eqref{conn_L_2} the following:
\begin{equation*}
  \begin{split}
   & ({w\gamma})^{\vee} = ({s_{\tau_1}s_{\tau_2}\dots{s}_{\tau_m}\gamma})^{\vee} =
    s^*_{\tau_1}({s_{\tau_2}\dots{s}_{\tau_m}\gamma})^{\vee} =
    s^*_{\tau_1}s^*_{\tau_2}({{s}_{\tau_3}\dots{s}_{\tau_m}\gamma})^{\vee} = \dots = \\
   & s^*_{\tau_1}{s}^*_{\tau_2}\dots{s}^*_{\tau_m}{\gamma}^{\vee} =
     {}^t{(s_{\tau_m}\dots{s}_{\tau_2}{s}_{\tau_1})}{\gamma}^{\vee} =
     {}^t{({s}_{\tau_1}s_{\tau_2}\dots{s}_{\tau_m})}^{-1}{\gamma}^{\vee} =
     {w}^*{\gamma}^{\vee}
  \end{split}
\end{equation*}

 4)  
 Further, \eqref{conn_L_5} holds since
 \begin{equation*}
    \begin{split}
    & \mathscr{B}_L(s_{\tau_i}{v}) = \langle B_L s_{\tau_i}{v}, s_{\tau_i}{v} \rangle =
     \langle s^*_{\tau_i}{B}_L {v}, s_{\tau_i}{v} \rangle =
     \langle {B}_L {v}, {v} \rangle = \mathscr{B}_L(v). \\
    & \mathscr{B}_L^{\vee}(s^{*}_{\tau_i}{u}) = \langle B_L^{\vee} s^{*}_{\tau_i}{u}, s^{*}_{\tau_i}{u} \rangle =
     \langle s_{\tau_i}{B}_L^{\vee} {u}, s^{*}_{\tau_i}{u} \rangle =
     \langle {B}_L^{\vee} {u}, {u} \rangle = \mathscr{B}_L^{\vee}(u). \\
    \end{split}
 \end{equation*}
 \qed

 \begin{remark}{\rm
1) The  mapping ${}^{\vee}$ defined on $L$ preserves dimensions and
  coincides with $B_L$,  ${}^{\vee} : L \longrightarrow L^{\vee}$ ,
  see \eqref{eq_gamma_vee}. From Proposition \ref{prop_link_diagr_conn}, heading 1)
  we see that the mapping ${}^{\vee}$ is also defined on $E$:
  \begin{equation}
    \begin{split}
     & {}^{\vee} : E = L \oplus M \longrightarrow L \longrightarrow L^{\vee},  \\
     & {}^{\vee} : M \longrightarrow 0.
    \end{split}
  \end{equation}

2) Note that $M$ can be the zero space. For example, for the Carter diagram $E_6(a_1)$  corresponding
to the conjugacy class lying in the Weyl group $E_6$, we have $L = E$, $M = 0$, see \cite[Tab.9, p.50]{Ca72}.
In this case, $\gamma = \gamma_L$ in \eqref{decomp_mu}. From now on, we consider the case $\dim M = 1$, i.e.,
one-dimensional extension of the root subset corresponding to $\Gamma$, see Section \ref{sec_proj_len}.
}
\end{remark}

%% file: 4inverse.tex
\subsection{The inverse quadratic form $\mathscr{B}^{\vee}_L$}

\subsubsection{The length of the projection of the root}
  \label{sec_proj_len}
   Let $\Pi_w(\theta)$ be the root subset obtained as the extension of $\Pi_w$ by some root $\theta$ linearly independent of $\Pi_w$ (i.e., $\theta \not\in L$), let $W_L(\theta)$ be the subgroup generated by $W_L$ and reflection $s_\theta$,  and $\varPhi_w(\theta)$ be the root subset that is the orbit of the action
   of $W_L(\theta)$ on $\Pi_w(\theta)$. Let $\gamma$ be any root from $\varPhi_w(\theta)$, $\gamma \not\in L$,
   $\gamma = \widetilde{w}\theta$ for some $\widetilde{w} \in W_L(\theta)$, i.e.,
   $s_\gamma = \widetilde{w}s_\theta\widetilde{w}^{-1} \in W_L(\theta)$.
   It is clear  that
\begin{equation}
  \label{varPhi_indep}
   W_L(\theta) = W_L(\gamma), \quad
   \varPhi_w(\theta) = \varPhi_w(\gamma) \text{ for any } \gamma \in \varPhi_w(\theta).
\end{equation}
   We will show that the length
   of the projection $\gamma_L$ of $\gamma$ (given by \eqref{varPhi_indep}) is independent of $\gamma$.
   We will see that this length essentially depends on the root system
   $\varPhi$ encompassing $\varPhi_w(\gamma)$, see Section \ref{sec_rat_p}.

   Since $B_L$ is positive definite, eigenvalues of $B_L$ are positive. Then also eigenvalues of $B^{-1}_L$ are positive, and the matrix $B^{-1}_L$ is positive definite as well.
   The quadratic form corresponding to inverse matrix $B^{-1}_L$ we call the
   {\it inverse quadratic form $\mathscr{B}^{\vee}_L$}.  The form $\mathscr{B}^{\vee}_L$
   is  positive definite. We will consider subgroup $W_L$ generated by reflections $\{s_{\alpha_1}, \dots, s_{\alpha_k}, s_{\beta_1}, \dots, s_{\beta_h} \}$ and its dual partial Weyl group
   $W^{\vee}_L$ generated by reflections
   $\{s^*_{\alpha_1}, \dots, s^*_{\alpha_k}, s^*_{\beta_1}, \dots, s^*_{\beta_h} \}$.

 \begin{figure}[H]
\centering
\includegraphics[scale=0.4]{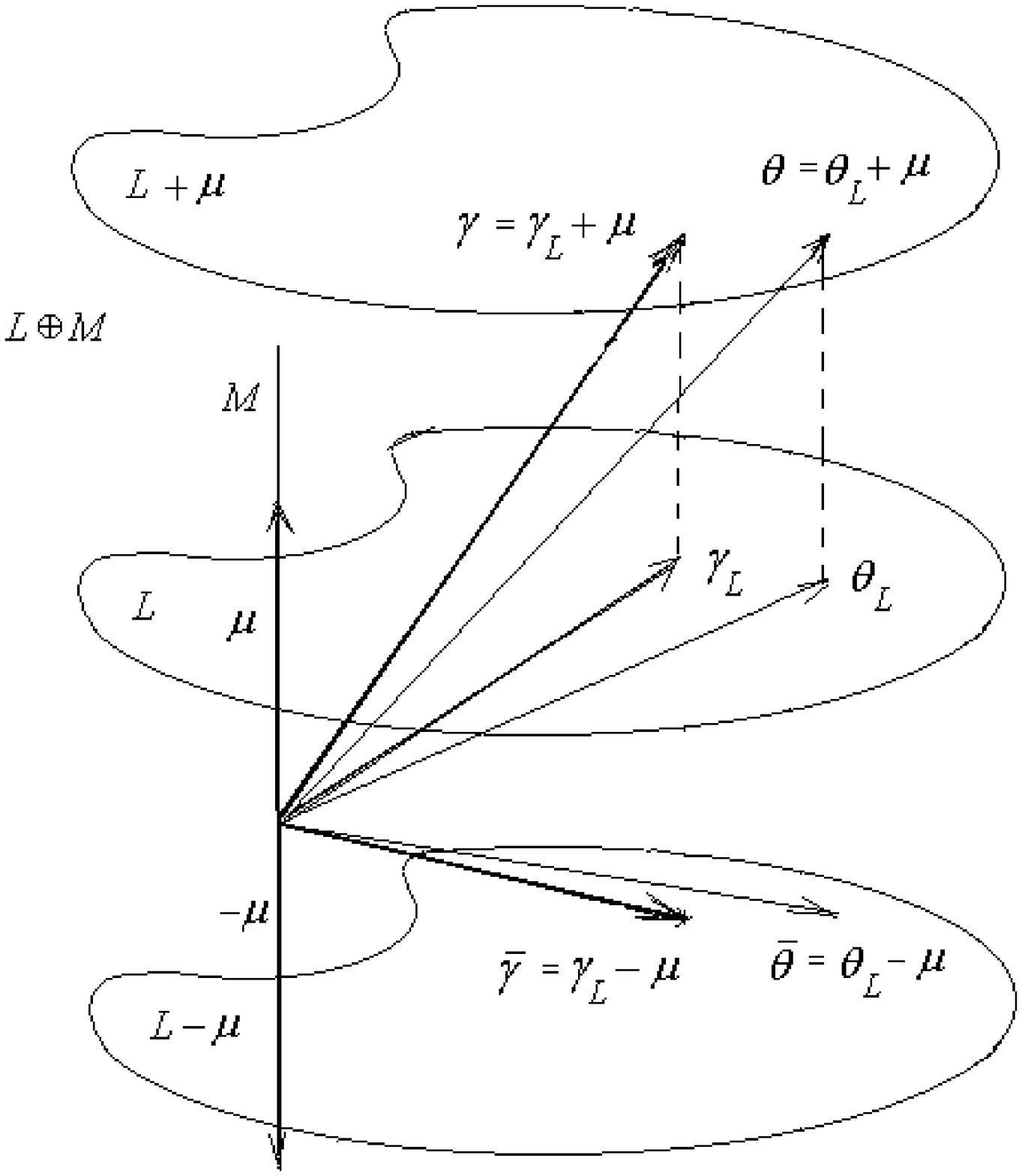}
 \caption{\hspace{3mm}The roots $\gamma = \gamma_L + \mu$ and $\overline{\gamma} = \gamma_L - \mu$.}
\label{L_and_mu}
\end{figure}

\begin{proposition}
   \label{prop_uniq_proj}
   Let $\gamma$ be any root belonging to the root subset $\varPhi_w(\theta)$
   and $\gamma \not\in L$.  Let us consider the decomposition
\begin{equation}
  \label{unique_gam_prim}
   \gamma = \gamma_L + \mu,  \text{ where }   \gamma_L \in L, ~ \mu \in M \text{ and } \mu \perp L,
\end{equation}
   see  Fig. \ref{L_and_mu}.

   1) The component $\mu$, up to sign, is a fixed vector for any root $\gamma \in \varPhi_w(\theta)$.

   2) The value $\mathscr{B}_L(\gamma_L)$ is constant for any root $\gamma \in \varPhi_w(\theta)$.

   3) The vector $\gamma = \gamma_L + \mu$ is a root in $\varPhi$ if and only if
   $\overline{\gamma} = \gamma_L - \mu$ is a root in
   $\varPhi$ (not necessarily both vectors belong to $\varPhi_w(\theta)$).

   4) If $\delta$ is the root form $\varPhi_w(\theta)$ such that $\delta^{\vee} = \gamma^{\vee}$
   then $\delta = \gamma_L + \mu$ or $\gamma_L - \mu$.
 \end{proposition}

\PerfProof
  1) Let $\delta \not\in L$ be another root from the root subset $\varPhi_w(\theta) = \varPhi_w(\gamma)$, i.e.,
  $\delta = w\gamma$ for some $w$ from ${W}_L(\gamma)$. We have $\delta = w\gamma_L + w\mu$.
 On the other hand, $\delta = \delta_L + t\mu$ for some $\delta_L \in L$ and some rational $t$, i.e.,
 $w\mu = t\mu$. Since $w$ preserves the length of $\mu$, we have $t = \pm{1}$, and
\begin{equation}
   \label{w_gam_plusmin}
      \delta = \delta_L  + \mu \text{ or } \delta = \delta_L - \mu  \text{ for some } \delta_L \in L.
\end{equation}

  2) By \eqref{unique_gam_prim} we have $\gamma_L \perp \mu$, and
  $\mathscr{B}(\gamma) = \mathscr{B}(\gamma_L) + \mathscr{B}(\mu)$.
  Here, $\mathscr{B}(\gamma) = 2$ since $\gamma$ is the root, and by 1), $\mathscr{B}(\mu)$ takes
  constant values. Therefore, $\mathscr{B}(\gamma_L)$ is constant. 
  By \eqref{restr_q} $\mathscr{B}_L(\gamma_L) = \mathscr{B}(\gamma_L)$,
  i.e.,  $\mathscr{B}_L(\gamma_L)$ is also constant for all $\gamma \in \varPhi_w(\theta)$.

  3) Let $\gamma$ be a root, i.e., $\mathscr{B}(\gamma) = \mathscr{B}(\gamma_L) + \mathscr{B}(\mu) = 2$.
  Then given $\overline{\gamma} = \gamma_L - \mu$, we have $\mathscr{B}(\overline{\gamma})
  = \mathscr{B}(\gamma_L) + \mathscr{B}(-\mu) = 2$, and $\overline{\gamma}$ is a root.

  4) From $\delta^{\vee} = \gamma^{\vee}$ and \eqref{conn_L_1} we have $B_L\delta_L = B_L\gamma_L$,
  i.e., $\delta_L = \gamma_L$. From \eqref{unique_gam_prim} and  heading 1), we have
  $\delta = \gamma_L + \mu$ or $\gamma_L - \mu$.
  \qed
  ~\\

  The group $W_L(\theta)$ acts on the linkages (= roots of $\varPhi_w(\theta)$)
  from the space $L(\theta)$, where $\dim L(\theta) = l + 1$.
  The group $W^{\vee}_L$ acts on the linkage diagrams, i.e., linkage labels vectors
  from $L^{\vee}$, where $\dim L^{\vee} = l$.
 \begin{corollary}
 Let $\gamma, \delta \in \varPhi_w(\theta)$, $w \in W_L(\theta)$, $w^{*} \in  W_L^{\vee}$.
 Then
  \begin{equation}
     \label{action_consist}
       \begin{split}
         & \delta = w \gamma \quad ( \in \mathbb{R}^{l+1}) \quad \Longrightarrow
         \quad \delta^{\vee}  = w^{*} \gamma^{\vee} \quad ( \in \mathbb{R}^{l}), \\
         & \delta^{\vee} = w^{*}\gamma^{\vee} \quad ( \in \mathbb{R}^{l}) \quad \Longrightarrow \quad
            w\gamma = \delta_L + \mu \quad \text{or} \quad w\gamma = \delta_L - \mu \quad ( \in \mathbb{R}^{l+1}). \\
       \end{split}
   \end{equation}
 \end{corollary}
  \qed

 \subsubsection{The length of the linkage labels vector}

 The following proposition checks whether or not the given vector $u^{\vee} \in L^{\vee}$
 is the linkage labels vector for a certain root $\gamma \in \varPhi_w(\theta)$.
\begin{proposition}
  \label{prop_unique_val}
   1) For any root $\gamma \in \varPhi_w(\theta)$, we have
\begin{equation}
    \label{length_dual_1}
     \mathscr{B}_L^{\vee}(\gamma^{\vee}) = \mathscr{B}_L(\gamma_L),
\end{equation}
~\\
     and, $\mathscr{B}^{\vee}_L(\gamma^{\vee})$ takes constant values for all roots  $\gamma \in \varPhi_w(\theta)$.

    2)  Let
\begin{equation}
   \label{length_dual_2}
     \mathscr{B}^{\vee}_L(u^{\vee}) =  \mathscr{B}^{\vee}_L(\gamma^{\vee}),
\end{equation}
~\\
where $\gamma^{\vee}$ is the linkage labels vector for some root $\gamma \in \varPhi_w(\theta)$,
and $u^{\vee}$ is a vector from $L^{\vee}$.
Then there exists the root $\delta \in \varPhi_w(\theta)$ such that $u^{\vee}$ coincides with $\delta^{\vee}$.
\end{proposition}
\PerfProof
  1) We have
 \begin{equation}
      \mathscr{B}^{\vee}_L(\gamma^{\vee}) =
      \langle B_L^{-1}\gamma^{\vee}, \gamma^{\vee} \rangle =
     \langle {\gamma_L},  B_L\gamma_L \rangle = \mathscr{B}_L(\gamma_L). \\
 \end{equation}
  By \eqref{conn_L_3}, for any $w^* \in W^{\vee}_L$, we have
  $w^*\gamma^{\vee} = (w\gamma)^{\vee}$,  and by \eqref{conn_dual_1}:
 ~\\
\begin{equation*}
      \mathscr{B}^{\vee}_L(w^*\gamma^{\vee}) =
      \langle B_L^{-1}w^*\gamma^{\vee}, w^*\gamma^{\vee} \rangle =
      \langle {w}B_L^{-1}\gamma^{\vee}, w^*\gamma^{\vee} \rangle =
      \langle {w}^{-1}{w}B_L^{-1}\gamma^{\vee}, \gamma^{\vee} \rangle =
     \mathscr{B}_L(\gamma_L). \\
 \end{equation*}
~\\
   2) Set $\delta = B_L^{-1}{u^{\vee}} + \mu$, where $\mu$ is
the fixed vector from \eqref{unique_gam_prim}.
  Then,
\begin{equation}
 \begin{split}
  \mathscr{B}(\delta) = & \mathscr{B}(B_L^{-1}u^{\vee}) +  \mathscr{B}(\mu) = \\
  & \mathscr{B}_L(B_L^{-1}u^{\vee}) +  \mathscr{B}(\mu) =
     \langle B_L(B_L^{-1}{u^{\vee}}), (B_L^{-1}u^{\vee}) \rangle + \mathscr{B}(\mu) =  \\
  & \langle u^{\vee}, B_L^{-1}{u^{\vee}} \rangle + \mathscr{B}(\mu) =
    \mathscr{B}^{\vee}_L(u^{\vee}) + \mathscr{B}(\mu) \stackrel{by \eqref{length_dual_2}}{=}  \mathscr{B}^{\vee}_L(\gamma^{\vee}) + \mathscr{B}(\mu) \stackrel{by \eqref{length_dual_1}}{=}   \\
  & \mathscr{B}_L({\gamma_L}) + \mathscr{B}(\mu) \stackrel{by \eqref{restr_q}}{=}
    \mathscr{B}({\gamma_L}) + \mathscr{B}(\mu) = \mathscr{B}(\gamma) = 2.
 \end{split}
\end{equation}
So $\mathscr{B}(\delta) = 2$, therefore $\delta$ is a root, see Section \ref{sec_Dynkin}.
Since $\delta = B_L^{-1}{u^{\vee}} + \mu$,
we get
\begin{equation*}
 \delta^{\vee} = B_L(B^{-1}_L{u^{\vee}} + \mu) = u^{\vee}. 
\end{equation*} 
 \qed

 Let us summarize:
\begin{theorem}
   \label{th_B_less_2}
     Let $\gamma^{\vee} \in L^{\vee}$ be the linkage labels vector
     corresponding to a certain root $\gamma \in \varPhi_w(\theta)$, i.e., $\gamma^{\vee} = B_{L}\gamma_L$.

     1) The root $\gamma \in \varPhi_w(\theta)$ is linearly independent of roots of $\varPhi_w$
     if and only if
     \begin{equation}
       \label{eq_p_less_2}
          \mathscr{B}^{\vee}_L(\gamma^{\vee}) < 2.
     \end{equation}

     2) We have $\mathscr{B}^{\vee}_L(\delta^{\vee}) = \mathscr{B}^{\vee}_L(\gamma^{\vee})$
     if and only if $\delta^{\vee} = B_{L}\delta_L$,  where $\delta_L$ is the projection on $L$ of some
     root $\delta$ belonging $\varPhi_w(\theta)$.

     Values $\mathscr{B}^{\vee}_L(\gamma^{\vee})$ are constant for any $\gamma \in \varPhi_w(\theta)$.

     (In a few cases $\varPhi_w$ can be extended to two different root subsets $\varPhi_w(\theta)$, see Section \ref{sec_rat_p})
\end{theorem}
  \PerfProof
    1) Let $\mathscr{B}^{\vee}_L(\gamma^{\vee}) = 2$. By Proposition \ref{prop_unique_val} we have  $\mathscr{B}_L(\gamma_L) = 2$, and  by Proposition \ref{restr_forms_coincide} we have also $\mathscr{B}(\gamma_L) = 2$. Since $\gamma \in \varPhi_w(\theta)$ is a root in $\varPhi$, then
    $\mathscr{B}(\gamma) = 2$ and $\mathscr{B}(\mu) = \mathscr{B}(\gamma) - \mathscr{B}(\gamma_L) = 0$.
  Therefore, $\mu = 0$,  and by \eqref{unique_gam_prim} $\gamma$ coincides with its projections on $\varPhi_w$, i.e., $\gamma$ is linearly depends on vectors of $\varPhi_w$.
     Conversely, let $\mathscr{B}^{\vee}_L(\gamma^{\vee}) < 2$, i.e., $\mathscr{B}_L(\gamma_L) < 2$.
     As above, we have  $\mathscr{B}(\gamma_L) < 2$
     and $\mathscr{B}(\mu) = \mathscr{B}(\gamma) - \mathscr{B}(\gamma_L) \neq 0$, i.e., $\mu \neq 0$ and $\gamma$ is linearly independent of roots of $\varPhi_w$.

    2)  If $\delta$ and $\gamma$ belong to the same root system $\varPhi_w(\theta)$ then
    $w\gamma = \delta$ for some $w \in W_L(\theta)$, and
    ${\delta}^{\vee} = (w\gamma)^{\vee} = w^{*}\gamma^{\vee}$.
    Thus, $\mathscr{B}^{\vee}_L$ takes the same values on $\gamma$ and $\delta$.
    The converse statement follows from Proposition \ref{prop_unique_val}, heading 2).
 \qed

\subsubsection{The rational number $p = \mathscr{B}^{\vee}_L(u^{\vee})$}
  \label{sec_rat_p}
  The rational number $p = \mathscr{B}^{\vee}_L(u^{\vee})$
  is the invariant characterizing the pair $\{ \Gamma, \varPhi_w(\theta) \}$,
  where $\Gamma$ is the Carter diagram  and  $\varPhi_w(\theta)$ is the root
  system extending the root subset $\varPhi_w$.

  By Theorem \ref{th_4.1.}, the conjugacy class with the Carter diagram $\Gamma$ belonging $\mathsf{C4}$ or $\mathsf{DE4}$ is uniquely characterized by the root subset $\varPhi_w$.
  Let us consider the set of linkage diagrams $S^{\vee} = \{ \gamma^{\vee} \}$ obtained from $\varPhi_w$ in such a way that  every $\gamma^{\vee} \in S^{\vee}$ is the linkage diagram associated  with a certain root $\gamma \in \varPhi_w(\theta)$, (i.e., $\gamma^{\vee} = B_{L}\gamma_L$), and the inverse quadratic form $\mathscr{B}^{\vee}_L$ takes the same rational value $p < 2$. We call this set {\it the extension set of the conjugacy class by p},
  or equivalently,  {\it the extension set of the Carter diagram $\Gamma$ by p}.
  According to Proposition \ref{prop_unique_val}, this set is well-defined and it is independent of
  the concrete root $\gamma \in \varPhi_w(\theta)$.  We denote this extension by $\{ \Gamma, p \}$.
  Note that  the extension $\{ \Gamma, p \}$ does depend on the choice of the root subsystem $\varPhi_w(\theta)$.
  Namely, we will see, that the Carter diagrams $D_5(a_1), D_6(a_1), D_6(a_2), D_7(a_1), D_7(a_2), D_5, D_6, D_7$ have two extensions -- any such extension can be either of type $D$ or of type $E$. Respectively, we call this extension either the {\it $D$-type extension} or the {\it $E$-type extension}. For example, $D_5(a_1)$ has two extension sets:  the $D$-type extension set $\{ D_5(a_1), 1 \}$ contains $10$ linkages, the $E$-type extension set $\{ D_5(a_1), \frac{5}{4} \}$ contains $32$ linkages, see Fig. \ref{D5a1_linkages} and Table \ref{val_Bmin1}.

%% file: 5loctets.tex
\subsection{Loctets and unicolored linkage diagrams}
  \label{sec_loctets_diagr}
 In this section we give the complete description of linkage diagrams for every linkage system.
 It turns that each linkage diagram containing at least one non-zero $\alpha$-label
 belongs to a certain $8$-cell linkage
 subsystem which we call loctet, see Fig. \ref{fig_loctets}.
 Every linkage system is the union of several loctets and several $\beta$-unicolored
 linkage diagrams. For the exact description, see Tables \ref{tab_seed_linkages_6},
  Fig. \ref{D4a1_linkages}-\ref{Dk_al_wind_rose}, Theorem \ref{th_full_descr}.

   \begin{proposition}
     Let $\Gamma$ be the Carter diagram from $\mathsf{C4}$, and $\gamma^{\vee} \in \mathscr{L}(\Gamma)$.
       \begin{equation}
          \gamma^{\vee} =
            \left (
            \begin{array}{c}
              a_1 \\
              a_2 \\
              a_3 \\
              \dots \\
              b_1 \\
              \dots \\
            \end{array}
            \right )
                =
            \left (
            \begin{array}{c}
              (\alpha_1, \gamma) \\
              (\alpha_2, \gamma) \\
              (\alpha_3, \gamma) \\
              \dots \\
              (\beta_1, \gamma) \\
              \dots \\
            \end{array}
            \right )
            \begin{array}{c}
              \alpha_1 \\
              \alpha_2 \\
              \alpha_3 \\
              \dots \\
              \beta_1 \\
              \dots \\
            \end{array}
       \end{equation}
       Among labels $a_i$ of the linkage diagram $\gamma^{\vee}$
       at least one label $a_i$ is equal $0$.
   \end{proposition}
 \PerfProof
   In the diagram $D_5(a_1)$, which is the part of every
   simply-laced connected Carter diagram containing $4$-cycle, the vertices
   $\{ \alpha_1, \alpha_2,  \alpha_3 \}$ connected to $\beta_1$, see Fig. \ref{D5a1_D4_pattern}.
   Thus, any root $\gamma$ can not be connected to all $\alpha_i$, where $i = 1,2,3$,
   otherwise we get the contradiction with the case of Corollary 2.4 from \cite{St10},
   see Proposition \ref{prop_numb_ep}, heading 2).
 \qed

  The following proposition explains relations between linkage diagrams depicted in Fig. \ref{fig_loctets} and shows that every linkage diagram containing at least one non-zero $\alpha$-label belongs to one of the loctets $L_{12}$,  $L_{13}$, $L_{23}$.

\begin{proposition}
  \label{prop_generating_diagr}
  1) The linkage labels $\gamma^{\vee}_{ij}(n)$ depicted in Fig. \ref{fig_loctets}
  are connected by means of dual reflections $s^{*}_{\alpha_i}$, where $i = 1,2,3$,
  and reflection $s^{*}_{\beta_1}$ as follows:
\begin{equation}
  \label{refl_dual_linkages}
  \begin{split}
   & s^{*}_{\alpha_k}\gamma^{\vee}_{ij}(8) = \gamma^{\vee}_{ij}(7), \quad
     s^{*}_{\alpha_k}\gamma^{\vee}_{ij}(1) = \gamma^{\vee}_{ij}(2), \\
   & s^{*}_{\beta_1}\gamma^{\vee}_{ij}(7) = \gamma^{\vee}_{ij}(6), \quad
     s^{*}_{\beta_1}\gamma^{\vee}_{ij}(2) = \gamma^{\vee}_{ij}(3), \\
    s^{*}_{\alpha_i}\gamma^{\vee}_{ij}(6) & =
     s^{*}_{\alpha_j}\gamma^{\vee}_{ij}(3) = \gamma^{\vee}_{ij}(4), \quad
     s^{*}_{\alpha_j}\gamma^{\vee}_{ij}(6) =
     s^{*}_{\alpha_i}\gamma^{\vee}_{ij}(3) = \gamma^{\vee}_{ij}(5),
  \end{split}
\end{equation}
where $\{ i,j, k \} = \{1, 2, 3 \}$.
Relations of the last line in \eqref{refl_dual_linkages} hold up to permutation of
indices $n = 4$ and $n = 5$ in  $\gamma^{\vee}_{ij}(n)$.

  2) If $\gamma^{\vee}$ contains exactly two non-zero labels $a_i$, $a_j$ (corresponding to
   coordinates $\alpha_i, \alpha_j$), then $\gamma^{\vee}$ is one of the following linkage diagrams:
   \begin{equation*}
      \gamma^{\vee}_{ij}(3), \quad \gamma^{\vee}_{ij}(4), \quad
      \gamma^{\vee}_{ij}(5), \quad \gamma^{\vee}_{ij}(6).
   \end{equation*}

   3) If $\gamma^{\vee}$ contains exactly one non-zero labels $a_i$ (corresponding to
   $\alpha_i$), then $\gamma^{\vee}$ is one of the following linkage diagrams:
   \begin{equation*}
      \gamma^{\vee}_{ij}(1), \quad \gamma^{\vee}_{ij}(2), \quad
      \gamma^{\vee}_{ij}(7), \quad \gamma^{\vee}_{ij}(8).
   \end{equation*}
\end{proposition}

\PerfProof
  1)   By \eqref{dual_refl} we have
  \begin{equation}
       \label{dual_refl_alpha_beta}
      \begin{split}
       (s^*_{\alpha_i}\gamma^{\vee})_{\alpha_k} =
        \begin{cases}
            -\gamma^{\vee}_{\alpha_i},  \text{ for } k = i, \\
            \hspace{3mm} \gamma^{\vee}_{\alpha_k},  \text{ for } k \neq i, \\
        \end{cases} \quad
        (s^*_{\alpha_i}\gamma^{\vee})_{\beta_k} =
        \begin{cases}
            \gamma^{\vee}_{\beta_k} + \gamma^{\vee}_{\alpha_i}, \text{ for } (\beta_k, \alpha_i) = -1, \\
            \gamma^{\vee}_{\beta_k} - \gamma^{\vee}_{\alpha_i},  \text{ for } (\beta_k, \alpha_i) = 1, \\
            \gamma^{\vee}_{\beta_k},  \text{ for } (\beta_k, \alpha_i) = 0.
        \end{cases} \\
          \\
       (s^*_{\beta_i}\gamma^{\vee})_{\beta_k} =
        \begin{cases}
            -\gamma^{\vee}_{\beta_i}, \text{ for } k = i, \\
            \hspace{3mm} \gamma^{\vee}_{\beta_k}, \text{ for } k \neq i, \\
        \end{cases} \quad
        (s^*_{\beta_i}\gamma^{\vee})_{\alpha_k} =
        \begin{cases}
            \gamma^{\vee}_{\alpha_k} + \gamma^{\vee}_{\beta_i}, \text{ for } (\beta_i, \alpha_k) = -1, \\
            \gamma^{\vee}_{\alpha_k} - \gamma^{\vee}_{\beta_i}, \text{ for } (\beta_i, \alpha_k) = 1, \\
            \gamma^{\vee}_{\alpha_k},  \text{ for } (\beta_i, \alpha_k) = 0.
        \end{cases} \\
       \end{split}
    \end{equation}
  We show \eqref{refl_dual_linkages} only for $\gamma^{\vee}_{ij}(n)$, where $n = 6, 7, 8$.
  One can get the remaining cases $n=1, 2, 3$ only by changing the sign of $\gamma^{\vee}_{ij}(n)$,
  see Fig. \ref{fig_loctets}.
  Applying $s^*_{\alpha_i}$ to $\gamma^{\vee}_{ij}(8)$,
  we have the first line of \eqref{refl_dual_linkages} as follows:
  \begin{equation*}
    \tiny
       \label{dual_refl_alpha_2}
    \begin{split}
     s^*_{\alpha_3}
      \left (
       \begin{array}{c}
          0 \\
          0 \\
          1 \\
          0 \\
       \dots \\
       \end{array} \right ) =
      \left (
       \begin{array}{c}
          0 \\
          0 \\
          -1 \\
          1 \\
       \dots \\
       \end{array} \right ), \quad
     s^*_{\alpha_2}
      \left (
       \begin{array}{c}
          0 \\
          1 \\
          0 \\
          0 \\
       \dots \\
       \end{array} \right ) =
      \left (
       \begin{array}{c}
          0 \\
          -1 \\
          0 \\
          1 \\
       \dots \\
       \end{array} \right ), \quad
       s^*_{\alpha_1}
      \left (
       \begin{array}{c}
          1 \\
          0 \\
          0 \\
          0 \\
       \dots \\
       \end{array} \right ) =
      \left (
       \begin{array}{c}
          -1 \\
          0 \\
          0 \\
          1 \\
       \dots \\
       \end{array} \right ).
     \end{split}
  \end{equation*}
  Applying $s^*_{\beta_1}$ to $\gamma^{\vee}_{ij}(7)$,
  we have the second line of \eqref{refl_dual_linkages}:
  \begin{equation*}
    \tiny
       \label{dual_refl_alpha_3}
    \begin{split}
     & s^*_{\beta_1}
      \left (
       \begin{array}{c}
          0 \\
          0 \\
          -1 \\
          1 \\
       \dots \\
       \end{array} \right ) =
      \left (
       \begin{array}{c}
          1 \\
          1 \\
          0 \\
          -1 \\
       \dots \\
       \end{array} \right ), \quad
      s^*_{\beta_1}
      \left (
       \begin{array}{c}
          0 \\
          -1 \\
          0 \\
          1 \\
       \dots \\
       \end{array} \right ) =
      \left (
       \begin{array}{c}
          1 \\
          0 \\
          1 \\
          -1 \\
       \dots \\
       \end{array} \right ), \quad
       s^*_{\beta_1}
      \left (
       \begin{array}{c}
          -1 \\
          0 \\
          0 \\
          1 \\
       \dots \\
       \end{array} \right ) =
      \left (
       \begin{array}{c}
          0 \\
          1 \\
          1 \\
          -1 \\
       \dots \\
       \end{array} \right ).
     \end{split}
  \end{equation*}
  Applying $s^*_{\alpha_i}$, $s^*_{\alpha_j}$ to $\gamma^{\vee}_{ij}(6)$
  we have the last line of \eqref{refl_dual_linkages}:
  \begin{equation*}
    \tiny
       \label{dual_refl_alpha_4}
    \begin{split}
     s^*_{\alpha_1}
      \left (
       \begin{array}{c}
          1 \\
          1 \\
          0 \\
          -1 \\
       \dots \\
       \end{array} \right ) =
      \left (
       \begin{array}{c}
          -1 \\
          1 \\
          0 \\
          0 \\
       \dots \\
       \end{array} \right ), \quad
     & s^*_{\alpha_2}
      \left (
       \begin{array}{c}
          1 \\
          1 \\
          0 \\
          -1 \\
       \dots \\
       \end{array} \right ) =
      \left (
       \begin{array}{c}
          1 \\
          -1 \\
          0 \\
          0 \\
       \dots \\
       \end{array} \right ), \quad
     s^*_{\alpha_1}
      \left (
       \begin{array}{c}
          1 \\
          0 \\
          1 \\
          -1 \\
       \dots \\
       \end{array} \right ) =
      \left (
       \begin{array}{c}
          -1 \\
          0 \\
          1 \\
          0 \\
       \dots \\
       \end{array} \right ), \quad
     s^*_{\alpha_3}
      \left (
       \begin{array}{c}
          1 \\
          0 \\
          1 \\
          -1 \\
       \dots \\
       \end{array} \right ) =
      \left (
       \begin{array}{c}
          1 \\
          0 \\
          -1 \\
          0 \\
       \dots \\
       \end{array} \right ), \\
     & \\
     & s^*_{\alpha_2}
      \left (
       \begin{array}{c}
          0 \\
          1 \\
          1 \\
          -1 \\
       \dots \\
       \end{array} \right ) =
      \left (
       \begin{array}{c}
          0 \\
          -1 \\
          1 \\
          0 \\
       \dots \\
       \end{array} \right ), \quad
     s^*_{\alpha_3}
      \left (
       \begin{array}{c}
          0 \\
          1 \\
          1 \\
          -1 \\
       \dots \\
       \end{array} \right ) =
      \left (
       \begin{array}{c}
          0 \\
          1 \\
          -1 \\
          0 \\
       \dots \\
       \end{array} \right ).
    \end{split}
  \end{equation*}

\begin{figure}[H]
\centering
\includegraphics[scale=0.8]{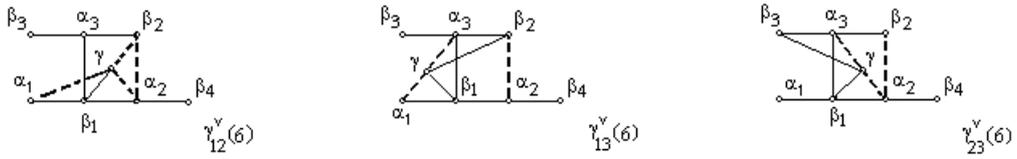}
 \caption{\hspace{3mm}The linkage diagrams $\gamma^{\vee}_{ij}(6)$ for $E_7(a_2)$, loctets $L_{ij}^b$.}
\label{E7a2_gamma_ij_6}
\end{figure}

\begin{figure}[H]
\centering
\includegraphics[scale=0.8]{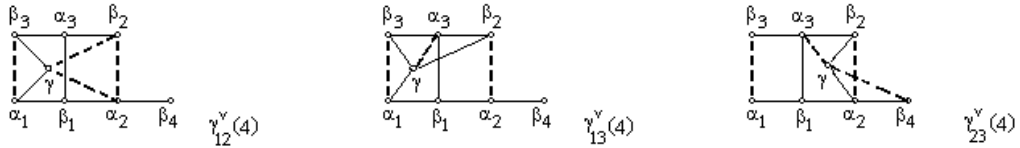}
 \caption{\hspace{3mm}The linkage diagrams $\gamma^{\vee}_{ij}(4)$ for $E_7(a_3)$, loctets $L_{ij}^b$}
\label{E7a3_gamma_ij_4}
\end{figure}

  2) Here, it suffices to prove that the label $b_1$ corresponding to the coordinate $\beta_1$
  is uniquely determined  by $\alpha_i, \alpha_j$. For the linkage diagram $\gamma^{\vee}_{ij}(3)$
  the statement follows from  Proposition \ref{prop_diagonal},(a), see Fig. \ref{square_n_diag},(a).
  For $E_7(a_1)$, the linkage diagrams $\gamma^{\vee}_{ij}(3)$, where $\{ij\} \in \{ \{ 12 \}, \{ 13 \}, \{ 23 \} \}$,  depicted
\begin{figure}[H]
\centering
\includegraphics[scale=0.8]{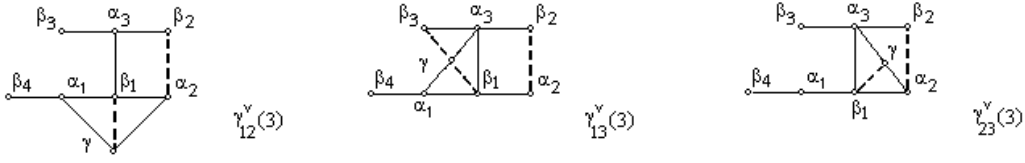}
 \caption{\hspace{3mm}The linkage diagrams $\gamma^{\vee}_{ij}(3)$ for $E_7(a_1)$, loctets $L_{ij}^b$.}
\label{E7a1_gamma_ij_3}
\end{figure}
~\\
  in Fig. \ref{E7a1_gamma_ij_3}, see the linkage system $E_7(a_1)$, loctets $L_{ij}^b$  in Fig. \ref{E7a1_linkages}.
  For the linkage diagram $\gamma^{\vee}_{ij}(6)$, the statement follows from
  Proposition \ref{prop_diagonal},(b), see Fig. \ref{square_n_diag},(b).
  For $E_7(a_2)$, the linkage diagrams $\gamma^{\vee}_{ij}(6)$, where $\{ij\} \in \{\{ 12 \}, \{ 13 \}, \{ 23 \}\}$,
  depicted in Fig. \ref{E7a2_gamma_ij_6}, see the linkage system $E_7(a_2)$, loctets $L_{ij}^b$
  in Fig. \ref{E7a2_linkages}. For the linkage diagram $\gamma^{\vee}_{ij}(4)$  and $\gamma^{\vee}_{ij}(5)$, the statement follows from Proposition
~\\
 \ref{prop_diagonal},(e), see Fig. \ref{square_n_diag},(e). For $E_7(a_3)$, the linkage diagrams $\gamma^{\vee}_{ij}(4)$, where $\{ij\} \in \{\{ 12 \}, \{ 13 \}, \{ 23 \}\}$,
  depicted in Fig. \ref{E7a3_gamma_ij_4}, see the linkage system $E_7(a_3)$, loctets $L_{ij}^b$
  in Fig. \ref{E7a3_linkages}.

  3) The label $b_1$ corresponding to the coordinate $\beta_1$
  takes two values from $\{ -1, 0, 1 \}$ depending on the value of $a_i$,
  see Fig. \ref{fig_loctets}. Indeed, if $a_i = 1$ then $b_1 = (\gamma, \beta_1) \neq 1$, otherwise the triangle
  $\{\alpha_i, \beta_1, \gamma \}$ contains exactly two dotted edges, i.e., contains $\widetilde{A}_3$,
  contradicting Lemma A.1 from \cite{St10}. Thus, $b_1 = -1$ or $b_1 = 0$.
  Respectively, we have linkage diagrams $\gamma^{\vee}_{ij}(2)$ or $\gamma^{\vee}_{ij}(8)$.
  If $a_i = -1$ then $b_1 = (\gamma, \beta_1) \neq -1$, otherwise the triangle
  $\{\alpha_i, \beta_1, \gamma \}$ does not contain any dotted edges. Thus, $b_1 = 1$ or $b_1 = 0$.
  Respectively, we have linkage diagrams $\gamma^{\vee}_{ij}(7)$ or $\gamma^{\vee}_{ij}(1)$.
\qed
~\\

 \begin{corollary}
   \label{corol_loctet}
    1) Any linkage diagram containing non-zero $\alpha$-label belongs to one of the loctets of the
       linkage system.

    2) Any linkage diagram of the loctet uniquely determines whole loctet.

    3) If two loctets have one common linkage diagram they coincide.

    4) Every linkage diagram from the linkage system either belongs to one of the loctets
    or is $\beta$-unicolored.

 \end{corollary}
    \PerfProof  Statements 1) and 4) follows from headings 2), 3) of Proposition \ref{prop_generating_diagr};
     statements 2) and 3) follow from heading 1) of Proposition \ref{prop_generating_diagr}.
    \qed
~\\

  The loctets of types
  $L_{12}$, $L_{13}$, $L_{23}$ are the {\bf main construction blocks} used for every
  linkage system, see all figures in  Fig. \ref{D4a1_linkages} - Fig. \ref{Dk_al_wind_rose}.

According to Proposition \ref{prop_generating_diagr}, any linkage diagram of a loctet
gives rise to the whole loctet.
By Theorem \ref{th_B_less_2}, to obtain all loctets associated with the given Carter diagram
it suffices to find linkage diagrams $\gamma^{\vee}_{ij}(n)$ for a certain fixed $n \in \{1, 2, \dots, 8 \}$  satisfying the inequality:
\begin{equation*}
       \mathscr{B}^{\vee}_L(\gamma^{\vee}_{ij}(n)) < 2.
\end{equation*}
The number of different loctets is defined by a number of different linkage diagrams $\gamma^{\vee}(n)$
for given fixed $n$, where $1 \leq n \leq 8$.
In that follows, we enumerate loctets by linkage diagrams $\gamma^{\vee}(8)$.
In Section \ref{sec_calc_gamma_8}, as an example, we show how to calculate
linkage diagrams $\gamma^{\vee}(8)$  of all loctets of $E_6(a_1)$. By Tables \ref{sol_inequal_1}-\ref{sol_inequal_6} one can recover the complete calculation of all linkage diagrams $\gamma^{\vee}(8)$ of all loctets 
for Carter diagrams $\Gamma \in \mathsf{C4} \coprod \mathsf{DE4}$.
The linkage diagrams $\gamma^{\vee}_{ij}(6)$ per every component and every loctet are listed
in Table \ref{tab_seed_linkages_6} for all linkage systems.

 The partial Cartan matrices $B_L$ and the inverse matrices $B^{-1}_L$
 for all Carter diagrams are presented in
 Tables \ref{tab_partial Cartan_1} - \ref{tab_partial Cartan_2}.

%% file: 6homogen.tex
\section{\sc\bf Enumeration of linkage diagrams, loctets and linkage systems}
  \label{sec_enum_loctets}
 In this section we demonstrate some calculation examples of linkage diagrams, loctets and
 linkage systems. These calculations are based on findings of Section \ref{sec_loctets_diagr}.
\subsection{Calculation of linkage diagrams $\gamma^{\vee}(8)$}
  \label{sec_calc_gamma_8}
 It seems a little easier to calculate the $8$th linkage diagram (we calculate it for every Carter diagram loctet) rather than to calculate any other linkage diagram of
 a loctet since $8$th linkage diagram contains $3$ zeroes among coordinates $\{ \alpha_1, \alpha_2, \alpha_3, \beta_1 \}$. We have
\begin{equation}
 \label{dual_linkage_1}
   \gamma^{\vee}(8) =
  \begin{cases}
     \{ a_1, a_2, a_3, 0,  b_2 \} \text{ for } D_5(a_1),  \\
     \{ a_1, a_2, a_3, 0,  b_2, b_3\} \text{ for } D_6(a_1), E_6(a_1), E_6(a_2),  \\
     \{ a_1, a_2, a_3, a_4, 0,  b_2 \} \text{ for } D_6(a_2), \\
     \{ a_1, a_2, a_3, a_4, 0,  b_2, b_3\} \text{ for } D_7(a_1), D_7(a_2), \\
     \{ a_1, a_2, a_3, 0,  b_2, b_3, b_4\}  \text { for }  E_7(a_1)$,  $E_7(a_2)$, $E_7(a_3)$, $E_7(a_4),
  \end{cases} \\
\end{equation}
   where $a_i = 0,  a_j = 0, a_k = 1, b_1 = 0$ and
   $\{i, j, k \} = \{1, 2, 3 \}$  for type $L_{ij}$.
For types $L_{12}$, $L_{13}$, $L_{23}$, there is the same quadratic form $q(\gamma^{\vee}(8))$  as the part of the inequality $\mathscr{B}^{\vee}_L(\gamma^{\vee}(8)) < 2$. The quadratic form $q(\gamma^{\vee}(8))$ is determined by principal submatrix associated with coordinates $\beta_2$, $\beta_3$, namely,
\begin{equation}
  \label{q_part}
  q(\gamma^{\vee}(8)) = d_{{\beta_2}{\beta_2}}b_2^2 + 2d_{{\beta_2}{\beta_3}}b_2{b_3} + d_{{\beta_3}{\beta_3}}b_3^2,
\end{equation}
where $d_{ij}$ is the $\{i, j\}$ slot of the inverse matrix $B^{-1}_L$.
The quadratic terms related to coordinates $\alpha_4$ or $\beta_4$ should be supplemented in the respective cases, see \eqref{dual_linkage_1}. The linear part $l(\gamma^{\vee}(8))$ and the free term $f(\gamma^{\vee}(8))$ of the inequality $\mathscr{B}^{\vee}_L(\gamma^{\vee}(8)) < 2$ are as follows:
\begin{equation*}
  l(\gamma^{\vee}(8)) =
    \begin{cases}
      2(d_{{\alpha_1}{\beta_2}}b_2 + d_{{\alpha_1}{\beta_3}}b_3) \text{ for } L_{23}, \\
      2(d_{{\alpha_2}{\beta_2}}b_2 + d_{{\alpha_2}{\beta_3}}b_3) \text{ for } L_{13}, \\
      2(d_{{\alpha_3}{\beta_2}}b_2 + d_{{\alpha_3}{\beta_3}}b_3) \text{ for } L_{12}, \\
    \end{cases}, \quad
  f(\gamma^{\vee}(8)) =
    \begin{cases}
      d_{{\alpha_1}{\alpha_1}}  \text{ for } L_{23}, \\
      d_{{\alpha_2}{\alpha_2}}  \text{ for } L_{13}, \\
      d_{{\alpha_2}{\alpha_3}}  \text{ for } L_{12}. \\
    \end{cases}
\end{equation*}
  We calculate the case $E_6(a_1)$. By Tables \ref{sol_inequal_1}-\ref{sol_inequal_4} 
  one can recover calculation for the remaining diagrams $\Gamma \in \mathsf{C4}$.

 \subsubsection{Calculation example for diagram $E_6(a_1)$}
   \label{calc_examp_E6a1}
 Here, $q = 4(b_2^2 + b_2{b_3} + b_3^2)$, see \eqref{q_part}.

  a) {\it Loctets $L_{12}$, $\gamma^{\vee}_{12}(8) = \{0, 0, 1, 0, b_2, b_3 \}$.}
\begin{equation*}
  \label{loctet_calc_E6a1}
  \begin{split}
    \mathscr{B}^{\vee}_L(\gamma^{\vee}_{12}(8)) =
    &  \frac{1}{3}(10 + 2(4b_2 + 5b_3) + 4(b_2^2 + b_2{b_3} + b_3^2)) < 2, \quad \text {i.e.,} \\
    &   \frac{1}{3}(2(b_2 + b_3)^2 + 2(b_2 + 2)^2 + 2(b_3 + \frac{5}{2})^2 - \frac{21}{2}) < 2, \\
    &   \frac{2}{3}((b_2 + b_3)^2 + (b_2 + 2)^2 + (b_3 + \frac{5}{2})^2) < 2 + \frac{7}{2} = \frac{11}{2}, \\
    &  (b_2 + b_3)^2 + (b_2 + 2)^2  + (b_3 + \frac{5}{2})^2 < \frac{33}{4}. \\
   \end{split}
\end{equation*}
  We get  $\gamma^{\vee}_{12}(8) =  \{ 0, 0, 1, 0, 0, -1 \}$  or $\gamma^{\vee}_{12}(8) = \{ 0, 0, 1, 0, -1, -1 \}$.
~\\

  b) {\it Loctets $L_{13}$, $\gamma^{\vee}_{13}(8) = \{0, 1, 0, 0, b_2, b_3 \}$.}
\begin{equation*}
  \begin{split}
    \mathscr{B}^{\vee}_L(\gamma^{\vee}_{13}(8)) =
    &  \frac{1}{3}(4 + 2(-b_2 + b_3) + 4(b_2^2 + b_2{b_3} + b_3^2)) < 2, \quad \text {i.e.,} \\
    &  \frac{1}{3}(2(b_2 + b_3) + 2(b_2 - \frac{1}{2})^2 + 2(b_3 + \frac{1}{2})^2 + 4 - \frac{1}{2} - \frac{1}{2})) < 2,  \\
    &  (b_2 + b_3)^2 + (b_2 - \frac{1}{2})^2  + (b_3 + \frac{1}{2})^2 < \frac{3}{2}. \\
   \end{split}
\end{equation*}
   We get  $\gamma^{\vee}_{13}(8) =  \{ 0, 1, 0, 0, 0, 0 \}$  or  $\gamma^{\vee}_{13}(8) = \{ 0, 1, 0, 0,  1, -1 \}$.
~\\

c) {\it Loctets $L_{23}$, $\gamma^{\vee}_{23}(8) = \{1, 0, 0, 0, b_2, b_3 \}$.}
\begin{equation*}
  \begin{split}
    \mathscr{B}^{\vee}_L(\gamma^{\vee}_{23}(8)) =
    &  \frac{1}{3}(4 + 2(b_2 + 2b_3) + 4(b_2^2 + b_2{b_3} + b_3^2)) < 2, \quad \text {i.e.,} \\
    &  (b_2 + b_3)^2 + (b_2 + \frac{1}{2})^2  + (b_3 + 1)^2 < \frac{9}{4}. \\
   \end{split}
\end{equation*}
 Here,  $\gamma^{\vee}_{23}(8) =  \{ 0, 1, 0, 0, 0, 0 \}$  or $\gamma^{\vee}_{23}(8) = \{ 0, 1, 0, 0,  0, -1 \}$.
 All loctets of the linkage system $E_6(a_1)$ are depicted in Fig. \ref{E6a1_linkages}.

\subsection{Calculation of the $\beta$-unicolored linkage diagrams}
  \label{sec_calc_homog}
 Now we consider $\beta$-unicolored linkage diagrams.
 Let  $a_i = (\alpha_i, \gamma)$,  $b_i = (\beta_i, \gamma)$, where $i = 1,2,3,\dots$, be
 linkage labels of the linkage $\gamma$. Since $\gamma$ is $\beta$-unicolored linkage,
 we have $a_i = 0$ for every $i = 1,2,3, \dots$. In addition, we note that
 \begin{equation}
   \label{b1_is_0}
       b_1 = (\gamma, \beta_1) = 0.
 \end{equation}
 Eq. \eqref{b1_is_0} holds for Carter diagrams containing a $D_5(a_1)$ or $D_4$ with predefined numbering of vertices $\{ \alpha_1, \alpha_2, \alpha_3, \beta_1 \}$ as in Fig. \ref{D5a1_D4_pattern},
 since otherwise the linkage diagram contains $5$-vertex
 subdiagram  $\{ \alpha_1, \alpha_2, \alpha_3, \beta_1, \gamma \}$ that is the extended Dynkin diagram $\widetilde{D}_4$.
 In other words, \eqref{b1_is_0} holds for all Carter diagrams from Tables \ref{tab_partial Cartan_1}, \ref{tab_partial Cartan_2} except for $D_4(a_1)$, and \eqref{b1_is_0} also holds for Dynkin diagrams $E_n$, $D_n$.
 Note that
 \begin{equation*}
   \mathscr{B}^{\vee}(\gamma^{\vee}) = \mathscr{B}^{\vee}(-\gamma^{\vee}),
 \end{equation*}
 and solving the inequality $\mathscr{B}^{\vee}(\gamma^{\vee}) < 2$
 we can assume that $b_2 > 0$, or $b_2 =0, b_3 > 0$, etc.
 We present here calculations only for cases $E_6(a_1)$, $E_6(a_2)$ and $E_7(a_1)$.
 By Tables \ref{homog_inequal_1}-\ref{homog_inequal_3} one can recover the remaining cases.

\subsubsection{$\beta$-unicolored linkage diagrams for $E_6(a_1)$ and $E_6(a_2)$}
  \label{sec_homog_E6a1a2}
For these cases, the $\beta$-unicolored linkage diagrams
coincide since the principal submatrices associated with coordinates $\beta_2$, $\beta_3$
for cases $E_6(a_1)$ and $E_6(a_2)$ coincide. This submatrix is the following $2\times2$ submatrix of $B^{-1}_L$:
\begin{equation*}
  \label{homog_E6a1_E6a2}
   \frac{1}{3}
   \left [
  \begin{array}{cc}
     4 & 2 \\
     2 & 4 \\
  \end{array}
   \right ],
\end{equation*}
see Table \ref{tab_partial Cartan_1}. By \eqref{eq_p_less_2} from Theorem \ref{th_B_less_2}, we have
\begin{equation}
  \label{homog_E6a1_E6a2_2}
    \mathscr{B}^{\vee}(\gamma^{\vee}) =
  \frac{2}{3}(2b_2^2 + 2b_3^2 + 2b_2{b_3}) < 2, \quad \text {i.e.,} \quad
     b_2^2 + b_3^2 + (b_2 + b_3)^2 < 3.
\end{equation}
 There are exactly $6$ solutions of the inequality \eqref{homog_E6a1_E6a2_2}, the corresponding
linkage diagrams are:
\begin{equation}
  \label{homog_E6a1_E6a2_3}
   \begin{array}{lll}
      & \{ 0, 0, 0, 0, 0, 1 \},  & \{ 0, 0, 0, 0, 0, -1 \}, \\
      & \{ 0, 0, 0, 0, 1, 0 \},  & \{ 0, 0, 0, 0, -1, 0 \}, \\
      & \{ 0, 0, 0, 0, 1, -1 \},  & \{ 0, 0, 0, 0, -1, 1 \}, \\
   \end{array}
\end{equation}
 see linkage diagrams located outside of the loctets in linkage systems $E_6(a_1)$, $E_6(a_2)$,
 Fig. \ref{E6a1_linkages},  Fig. \ref{E6a2_linkages}.

 \subsubsection{$\beta$-unicolored linkage diagrams for $E_7(a_1)$}
 The principal $3\times3$ submatrices  of $B^{-1}_L$ associated with coordinates $\beta_2$, $\beta_3$, $\beta_4$
 is as follows:
\begin{equation*}
  \label{homog_E7a1}
   \frac{1}{2}
   \left [
  \begin{array}{ccc}
     3 & 2 & 1 \\
     2 & 4 & 2 \\
     1 & 2 & 3 \\
  \end{array}
   \right ],
\end{equation*}
see Table \ref{tab_partial Cartan_2}. Then we have
 \begin{equation}
  \label{homog_E7a1_2}
   \begin{split}
  \mathscr{B}^{\vee}(\gamma^{\vee}) = &
    \frac{1}{2}(3b_2^2 + 4b_3^2 + 3b_4^2 + 4b_2{b_3} + 2b_2{b_4} + 4b_3{b_4}) < 2, \quad
    \text {i.e.,}\\
  & 2(b_2 + b_3)^2 + 2(b_3 + b_4)^2 + (b_2 + b_4)^2 < 4.
  \end{split}
\end{equation}
 There are exactly $8$ solutions of the inequality \eqref{homog_E7a1_2}, the corresponding
  linkage diagrams are:
\begin{equation*}
  \label{homog_E7a1_3}
   \begin{array}{lll}
      & \{ 0, 0, 0, 0, 1, -1, 0 \},  & \{ 0, 0, 0, 0, -1, 1, 0 \}, \\
      & \{ 0, 0, 0, 0, 0, 1, -1 \},  & \{ 0, 0, 0, 0, 0, -1, 1 \}, \\
      & \{ 0, 0, 0, 0, 0, 0, 1 \},  & \{ 0, 0, 0, 0, 0, 0, -1 \}, \\
      & \{ 0, 0, 0, 0, 1, 0, 0 \},  & \{ 0, 0, 0, 0, -1, 0, 0 \}, \\
   \end{array}
\end{equation*}
 see linkage diagrams located outside of the loctets in linkage systems $E_7(a_1)$,
 Fig. \ref{E6a1_linkages}.

%% file: 7Dynkin.tex
\subsection{Linkage systems for simply-laced Dynkin diagrams}

In this section we extend the previous results to simply-laced Dynkin diagrams.
First, we find such simply-laced Dynkin diagrams that each of them determines only one conjugacy class.

\begin{figure}[H]
\centering
\includegraphics[scale=0.8]{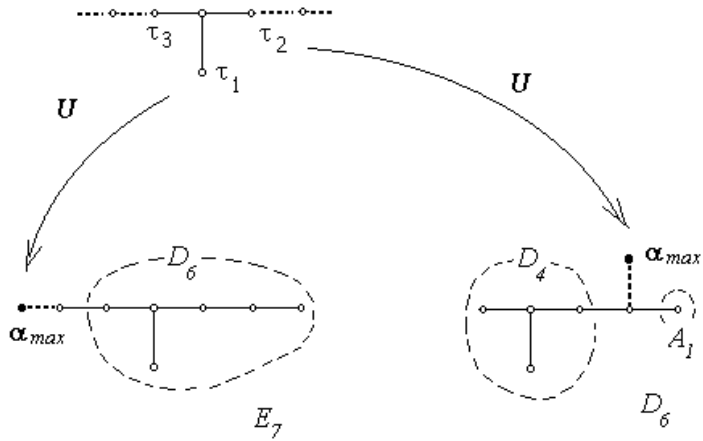}
\caption{\hspace{3mm}The mapping $U: \varPhi(\Gamma) \longrightarrow \varPhi(E_7)$}
\label{E7_D6_D4__3roots}
\end{figure}

\begin{remark}[on isolated roots]
 \label{rem_isolated}
{\rm
  1) Let $\alpha, \beta$ be two roots from $\varPhi(\Gamma)$, $\alpha \perp \beta$.
     The root $\beta$ is said to be {\it isolated root in sense of orthogonality to $\alpha$}
     or, for brevity, {\it isolated root for $\alpha$}
     if any root $\gamma$ connected to $\beta$ (i.e., $(\gamma, \beta) \neq 0$) is non-orthogonal to $\alpha$.
     For example, $\alpha_{max}(D_6)$, the maximal root from $\varPhi(D_6)$  is orthogonal to
     $\varPhi(D_4) \oplus \{ \tau_1 \}$, where $\tau_1$ is the single root in  $\varPhi(A_1)$, see Fig. \ref{E7_D6_D4__3roots}. The root $\tau_1$ is isolated for $\alpha_{max}(D_6)$.

  2)  Let $\varPhi(\Gamma),  \varPhi(\Gamma') \subset W$, where $W$ is the Weyl group, let $U \in W$ be the map
  $U : \varPhi(\Gamma) \longrightarrow \varPhi(\Gamma')$. The root $\beta$ is isolated
  for $\alpha$ if and only if $U\beta$ is isolated for $U\alpha$. Indeed, since $U$ preserves orthogonality, we have:
  \begin{equation*}
    \begin{split}
    & \{ \beta \perp \alpha,  \quad \beta \text{ is connected to } \gamma, \quad \gamma  \not\perp \alpha \}
    \Longleftrightarrow \\
    & \{ U\beta \perp U\alpha,  \quad U\beta \text{ is connected to } U\gamma, \quad U\gamma  \not\perp U\alpha \}.
    \end{split}
  \end{equation*}
 \qed
}
\end{remark}

\begin{proposition} 
  \label{diagr_El_Dl}
   Let $\{ \tau_1, \dots, \tau_l \}$, $\{ \tau'_1, \dots, \tau'_l \}$
   be two root subsets (with not necessarily simple roots) corresponding to one of the
   Dynkin diagram $\Gamma = E_6, E_7, E_8, D_l$ (considered as Carter diagrams). There is $U \in W$ such that
     \begin{equation}
       \label{map_U_Dynkin}
         U\tau_i = \tau'_i, \text{ where } i = 1,\dots, l.
     \end{equation}
   (That means that every Dynkin diagram $E_6, E_7, E_8, D_l$  determines only one conjugacy class.)
\end{proposition}

\PerfProof

     1) \underline{Case $\Gamma = E_6, E_7, E_8$.} First, we will show that two triples of orthogonal roots
     \begin{equation}
       \label{two_sets_3roots}
         S = \{ \tau_1, \tau_2, \tau_3 \}, \qquad S'= \{ \tau'_1, \tau'_2, \tau'_3 \}
     \end{equation}
     forming the root system $D_4$ (i.e.,  such that any $\tau_i$ (resp. $\tau'_i$) is adjacent
     to the branch point of the diagram $\Gamma$) are equivalent under the Weyl group $W(\Gamma)$.
     By \cite[Lemma 11]{Ca72}, \cite[Corollary A.6]{St10} any two sets of $3$ orthogonal
     roots in $\varPhi(E_6)$ (resp. $\varPhi(E_8)$) are equivalent under $W(E_6)$ (resp. $W(E_8)$).
     For $\varPhi(E_7)$, this property is not correct,
     there are two  sets of $3$ orthogonal roots in $\varPhi(E_7)$, which are not equivalent
     under $W(E_7)$, see  \cite[Corollary A.6]{St10}. The reason of this obstacle
     in the case $E_7$ is the \lq\lq unlucky \rq\rq location of
     the maximal root in the root system of $D_6$, see Fig. \ref{E7_D6_D4__3roots}.
     To solve this we use Corollary A.5 from \cite{St10}. According to this corollary any two sets of $2$ orthogonal roots in $\varPhi(E_7)$  are equivalent. Let $U \in W(E_7)$ map $\tau_1$ to $\alpha_{max}(E_7)$, the maximal  root in $E_7$, and $\tau_2$ to $\alpha_{max}(D_6)$, see Fig. \ref{E7_D6_D4__3roots}.  Then $\tau_3$ is mapped into a root of $\varPhi(D_4) \oplus \{ \tau_1\}$, $\alpha_1 \in \varPhi(A_1)$. By Remark \ref{rem_isolated},
     since $\alpha_1 \in \varPhi(A_1)$ is isolated for $\alpha_{max}(D_6)$,
     and $\tau_3$ is not isolated for $\tau_2$ then  $\tau_3$ is mapped into a root of $\varPhi(D_4)$. We can choose $U$ in such a way that $U\tau_3$ is mapped to any root in $\varPhi(D_4)$, for example, to $\alpha_{max}(D_4)$.
  ~\\
     The same is true for the set $S'$ in \eqref{two_sets_3roots}. Therefore, sets $S$ and $S'$ are equivalent.
     If sets $S$ and $S'$ from \eqref{two_sets_3roots} are equivalent the branch points
     which bind the roots of each of these sets are mapped by the mapping $U$ from \eqref{map_U_Dynkin},
     see Section 4.6
     from \cite{St10}. By Lemma 4.7 from \cite{St10}, the map $U$ can be extended to
     the map \eqref{map_U_Dynkin} for every pair $\{\tau_i, \tau'_i \}$.
\begin{figure}[H]
\centering
\includegraphics[scale=0.8]{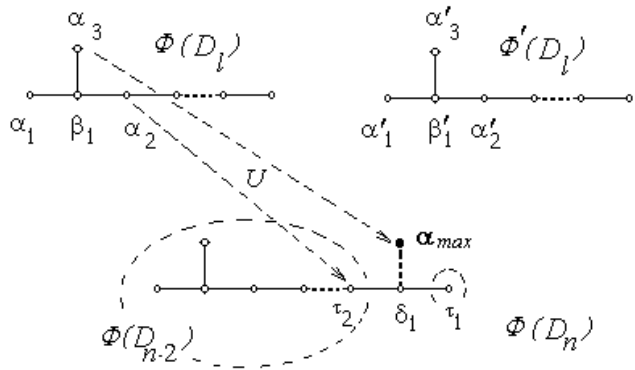}
\caption{\hspace{3mm}The mapping $U: \varPhi(D_l) \longrightarrow \varPhi'(D_l)$}
\label{Dm_Dl__3roots}
\end{figure}

     2) \underline{Case $\Gamma = D_l$.} Let us show that two roots subsets $\varPhi(D_l)$ and $\varPhi'(D_l)$
     lying in the Weyl group $W$ are equivalent. In the case $\varPhi(D_l) , \varPhi'(D_l) \subset \varPhi(E_n)$,
     where $n \leq 8$, we repeat the same arguments as in heading 1). Let now $\varPhi(D_l) , \varPhi'(D_l) \subset \varPhi(D_n)$, where $l \leq n$.
     Let us take such $U \in W$ that $U : \alpha_3 \longrightarrow \alpha_{max}$, see Fig. \ref{Dm_Dl__3roots}.
     Then $\alpha_1, \alpha_2$ are mapped by $U$ into $\varPhi(D_{n-2}) \oplus \{ \tau_1 \}$, see \cite[\S A.4]{St10}.
     We have $U\alpha_1 \perp \alpha_{max}$, $U\alpha_2 \perp \alpha_{max}$.
     By Remark \ref{rem_isolated}, since $\alpha_1$ is isolated for $\alpha_3$ then $\alpha_1$ is mapped into $\tau_1$,
     and $U\alpha_2 \in \varPhi(D_{n-2})$:
  $$
     U\alpha_1 = \tau_1, \quad U\alpha_3 = \alpha_{max}, \quad  U\alpha_2 \in \varPhi(D_{n-2}).
  $$
  Applying a certain mapping $U' \in W({D_{n-2}})$
  to the mapping $U$, we change $U$ in such a way that $U\alpha_2 = \tau_2$ and $U\alpha_1$, $U\alpha_3$ are not changed.
  Therefore, $S = \{ \alpha_1, \alpha_2, \alpha_3 \} \subset \varPhi(D_l)$
  is mapped  into $\{ \tau_1, \tau_2, \alpha_{max} \}$. Similarly, there exists
  $U'$ mapping $S' = \{ \alpha'_1, \alpha'_2, \alpha'_3 \} \subset \varPhi'(D_l)$ into $\{ \tau_1, \tau_2, \alpha_{max} \}$. Then $U^{-1}U'$ maps $S$ into $S'$.
  Further, applying  Lemmas $4.6$ and $4.7$ from \cite{St10}, we get the extension $\tilde{U}$ of $U^{-1}U'$,
  so that $\tilde{U}: \varPhi(D_l) \longrightarrow \varPhi'(D_l)$.
\qed

\begin{remark}
  \label{two_class_Al}
{\rm
    Note that Proposition \ref{diagr_El_Dl} does not hold for $A_l$.
    There are two conjugacy classes of type $A_5$ in $W(E_7)$ and in $W(D_6)$,
    two conjugacy classes of type $A_7$ in $W(E_8)$, see \cite[p. 31]{Ca72} and \cite[Lemma 27]{Ca72}.
    Moreover, for the Carter diagram $A_3$, there exist three conjugacy classes in $W(D_4)$.
   } \qed
\end{remark}

In order to build the linkage systems for Dynkin diagrams $E_n$, $D_n$,
by Proposition \ref{diagr_El_Dl},  we can use the technique of the partial Cartan matrix,
linkage diagrams and loctets from Section \ref{sec_Cartan_matr_ccl}. Note that the partial Cartan matrix for Dynkin diagrams coincides with the usual Cartan matrix ${\bf B}$ associated with the given Dynkin diagram. Since $E_8$ does not have linkage diagrams, see Remark \ref{rem_l_less_8}, heading 1), we are interested only
in $E_6$, $E_7$, $D_n$. In cases $E_6$, $E_7$, $D_5$, $D_6$, $D_7$, for the Cartan matrix ${\bf B}$ and its inverse  ${\bf B}^{-1}$, see Table \ref{tab_partial Cartan_3}.
One can recover the complete calculation  of linkage diagrams $\gamma^{\vee}(8)$ of all loctets
of $E_6$, $E_7$, $D_5$, $D_6$, $D_7$  by means of Tables \ref{sol_inequal_5}, \ref{sol_inequal_6}.
The $\beta$-unicolored linkage diagrams look as follows:
\begin{equation*}
  \gamma^{\vee} =
  \begin{cases}
       \{0, 0, 0, 0, b_2, b_3 \} \text{ for } E_6, \\
       \{0, 0, 0, a_4, 0, b_2, b_3 \} \text{ for } E_7, \\
       \{0, 0, 0, 0, b_2 \}  \text{ for } D_5, \\
       \{0, 0, 0, a_4, 0, b_2 \} \text{ for } D_6, \\
       \{0, 0, 0, a_4, 0, b_2, b_3 \} \text{ for } D_7,
  \end{cases}
\end{equation*}
see Tables \ref{homog_inequal_1}-\ref{homog_inequal_3}. The $\beta$-unicolored linkage diagrams
are located outside of all loctets, see Fig. \ref{E6pure_loctets_norm_gam8},
Fig. \ref{E7pure_linkage_system}, Fig. \ref{D5pure_loctets}, Fig. \ref{D6pure_loctets}.

Note that for $E_6$, the principal matrix associated with coordinates $\beta_2$, $\beta_3$ coincide with
the principal matrix for the Carter diagrams $E_6(a_1)$, $E_6(a_2)$, see Section \ref{sec_homog_E6a1a2},
and, consequently, $\beta$-unicolored linkage diagrams coincide with these diagrams for
$E_6(a_1)$, $E_6(a_2)$, see $6$ solutions \eqref{homog_E6a1_E6a2_3} of the inequality \eqref{homog_E6a1_E6a2_2}.

\subsubsection{The Dynkin diagrams of $E$-type and $D$-type}

\begin{figure}[H]
\centering
\includegraphics[scale=0.8]{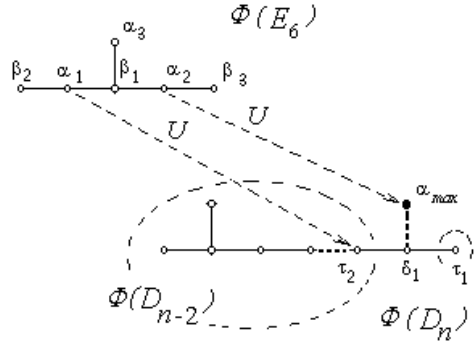}
\caption{\hspace{3mm}The conjugacy class $E_6$ is not contained in any $W(D_n)$}
\label{E6_not_in_Dl}
\end{figure}

The relation between Dynkin diagrams of $E$-type and $D$-type is asymmetrical in the following sense:
conjugacy classes $D_i$ are contained in $W(E_{i+1})$ for $i = 5,6,7$; the reverse
inclusions, however, are not true, as we see from the following:
\begin{proposition}
   \label{prop_E6_Dk}
    The conjugacy classes $E_i$, where $i = 6,7,8$, are not contained in any $W(D_n)$.
\end{proposition}
\PerfProof
   It suffices to prove the statement for $E_6$.
   Suppose, $E_6$ is contained in $W(D_n)$. Then subset $S = \{ \alpha_1, \alpha_2, \alpha_3 \}$
   from $E_6$ also belongs to $W(D_n)$, see Fig. \ref{E6_not_in_Dl}.  Let us choose the map $U \in W(D_n)$ in such a way that
   $U$ maps $\alpha_2 \in \varPhi(E_6)$ into $\alpha_{max} \in \varPhi(D_n)$.
   Then $\alpha_1$, $\alpha_3$ are mapped into  $\varPhi(D_{n-2}) \oplus \{ \tau_1 \}$, where
   $\tau_1$ is a single root from $A_1$,  see Fig. \ref{E6_not_in_Dl}. By Remark \ref{rem_isolated},
   the isolated root $\alpha_3$ is mapped into the isolated root $\tau_1$. Further, $U$ can be selected
   in such a way that $U$ maps $\alpha_1$ into $\tau_2$, and $U\alpha_2$, $U\alpha_3$ are not changed.
   According to \cite[\S4.6]{St10}, $U$ maps $\beta_1$ into $\delta_1$.
  ~\\
  Then the edge
  $\{ \beta_1, \alpha_2 \} \in E_6$ is mapped into the edge $\{ \delta_1, \alpha_{\max} \} \in D_n$.
  One can see that the image $U\beta_3 \neq U\beta_1 = \delta_1$ is connected to $U\alpha_2 = \alpha_{max}$, and, consequently,
  the root $\delta_1$ is contained in the support of $U\beta_3$, i.e., $U\beta_3 \not\perp \tau_1$.
  This contradicts the orthogonality of their preimages: $\beta_3 \perp \alpha_3$.
  \qed

%% file: 5a_extendable.tex
\section{\sc\bf Projection of linkage systems}
  \label{sec_proj}
\subsection{Simply extendable Carter diagrams}
  \label{sect_extendable}
\begin{figure}[H]
\centering
\includegraphics[scale=0.5]{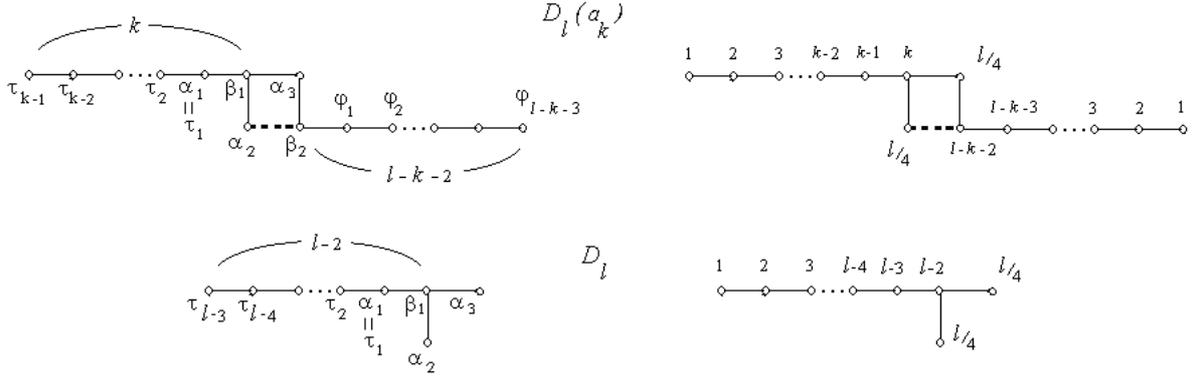}
\caption{
\hspace{3mm}The numerical labels (in the right-hand side) are diagonal
elements of $B^{-1}_L$}
\label{Dk_al_Carter_diagr_in_prop}
\end{figure}

 We say that the Carter diagram $\Gamma$ is {\it simply extendable in the vertex $\tau_p$} if
 the new diagram obtained by the extra vertex $\tau_{l+1}$ together with the additional connection edge
  $\{ \tau_p, \tau_{l+1} \}$  is also the Carter diagram. The extra vertex with given property
  is called {\it simply extendable}.
  We will show that extensibility in the vertex is closely
 associated with the value the diagonal element $b^{\vee}_{\tau_p,\tau_p}$  of the matrix $B_L^{-1}$, the inverse of the partial Cartan matrix $B_L$.
 \begin{proposition}
   \label{prop_determinants}
    1) The determinant of partial Cartan matrix $B_L$ is as follows:
  \begin{equation*}
      \det{B_L} = 
        \begin{cases}
           l + 1 \text{ for } A_l, \text{ where } l \geq 2; \text{ here}, B_L = {\bf B}, \\
           4 \text{ for } D_l, \text{ where } l \geq 4, \\
           4 \text{ for } D_l(a_k), \text{ where } l \geq 4. 
        \end{cases}
  \end{equation*}  
    2) Let $b^{\vee}_{\eta,\eta}$ be diagonal elements of the inverse matrix $B^{-1}_L$,
     where
      \begin{equation}
        \label{diag_elem_1}
        \begin{split}
          & \eta \in \{ \tau_1, \dots, \tau_{k-1}, \varphi_1, \dots, \varphi_{l-k-3},
          \alpha_2, \alpha_3, \beta_1, \beta_2\},  \text{ for } \Gamma = D_l(a_k),  \\
          & \eta \in \{ \tau_1, \dots, \tau_{l-3},
          \alpha_2, \alpha_3, \beta_1, \beta_2\},  \text{ for } \Gamma = D_l,
        \end{split}
      \end{equation}
     see Fig. \ref{Dk_al_Carter_diagr_in_prop}.
     For $\Gamma = D_l(a_k)$ or $\Gamma = D_l$,
      \begin{equation}
        \label{bi_eq_1}
     b^{\vee}_{\eta,\eta} =
     \begin{cases}
          & \displaystyle\frac{l}{4}, \text{ for } \eta = \alpha_2, \text{ or } \eta = \alpha_3, \\
          & d, \text{ for } \eta \text{ given in }  \eqref{diag_elem_1}.
     \end{cases}
     \end{equation}
Here, $d-1$ is the number of vertices remaining in the chain $A_{d-1}$ after removing the vertex $\eta$.
\end{proposition}
\PerfProof  1) This statement is easily verified for $A_1$, $A_2$, $D_4$ and $D_5$.
   By induction,
 \begin{equation*}
   \begin{split}
    & \det{{\bf B}(A_{l+1})} = 2\det{{\bf B}(A_{l})} -  \det{{\bf B}(A_{l-1})} =
        2(l+1) - l = l + 2, \\
    & \det{{\bf B}(D_{l+1})} = 2\det{{\bf B}(D_{l})} -  \det{{\bf B}(D_{l-1})} = 4,
   \end{split}
 \end{equation*}
 where  ${\bf B}(A_{l})$ (resp. ${\bf B}(D_{l})$)  is the Cartan matrix for $A_{l}$ (resp. $D_{l}$).

 The decomposition of $\det{B_L(D_{l+1}(a_k))}$ leads
 to the matrix $l \times l$ and another matrix $(l - 1) \times (l - 1)$ obtained by
 removing the $i$-th line and $i$-th column for some $i$. The matrix $(l - 1) \times (l - 1)$ corresponds
 to the Carter diagram either $D_{l-1}(a_k)$ or $D_{l-1}$. In both cases
 the determinant is $4$. By induction, we have
 \begin{equation*}
   \begin{split}
    & \det{B_L(D_{l+1}(a_k))} = 2\det{B_L(D_{l}(a_k))} -  \det{B_L(D_{l-1}(a_k))} = 4,
       \text{ for } k > 2 \text{ or } (l+1) - k - 2 > 2, \\
    & \det{B_L(D_{l+1}(a_k))} = 2\det{B_L(D_{l}(a_k))} -  \det{{\bf B}(D_{l-1})} = 4,
       \text{ for } k \leq 2, \text{ and } (l+1) - k - 2 \leq 2. \\
    \end{split}
 \end{equation*}

  2) For $\Gamma = D_l(a_k)$, we have
   \begin{equation}
    b^{\vee}_{\eta,\eta} =
    \begin{cases}
    & \displaystyle\frac{\det {B_L(D_{l}(a_k))}}{\det {B_L(D_{l+1}(a_k))}} = \frac{4}{4} = 1,
        \text{ for } $d = 1$, \\
    & \\
    & \displaystyle\frac{\det{\bf B}(A_{d-1}) \det B_L(D_{l-d}(a_k))}
                {\det B_L(D_{l+1}(a_k))} = \frac{d \cdot 4}{4} = d,  \text{ for } d > 1, \\
    & \\
    &  \displaystyle\frac{\det {{\bf B}(A_{l-1})}}{\det {B_L(D_{l+1}(a_k))}} = \frac{l}{4},  \text{ for }
      \eta = \alpha_2, \alpha_3.
     \end{cases}
   \end{equation}
   For $\Gamma = D_l$, we have
   \begin{equation}
    b^{\vee}_{\eta,\eta} =
    \begin{cases}
    & \displaystyle\frac{\det {{\bf B}(D_{l})}}{\det {{\bf B}(D_{l+1})}} = \frac{4}{4} = 1,
        \text{ for } $d = 1$, \\
    & \\
    & \displaystyle\frac{\det{\bf B}(A_{d-1}) \det {\bf B}(D_{l-d})}
                {\det {\bf B}(D_{l+1})} = \frac{d \cdot 4}{4} = d,  \text{ for } d > 1, \\
    & \\
    &  \displaystyle\frac{\det {{\bf B}(A_{l-1})}}{\det {{\bf B}(D_{l+1})}} = \frac{l}{4},  \text{ for }
      \eta = \alpha_2, \alpha_3.
     \end{cases}
   \end{equation}
   \qed

\begin{remark}
  \label{rem_val_le_2}
  {\rm
    For $D_l(a_k)$, where $l \geq 8$, we have $b^{\vee}_{\eta,\eta} < 2$ only for endpoints $\tau_{k-1}$ and
    $\varphi_{l - k -3}$, see Fig. \ref{Dk_al_Carter_diagr_in_prop}. For $D_l$, where $l \geq 8$, we have $b^{\vee}_{\eta,\eta} < 2$ only for endpoint $\tau_{l-3}$, see Fig. \ref{Dk_al_Carter_diagr_in_prop}.
    For $D_4(a_k)$, $D_5(a_k)$, $D_6(a_k)$, $D_7(a_k)$ and $D_4$, $D_5$, $D_6$, $D_7$,
    we have $b^{\vee}_{\eta,\eta} < 2$ also for $\eta = \alpha_2, \alpha_3$.
    }
\end{remark}

 \begin{proposition}
  \label{prop_crit_sim_ext}
   The Carter diagram $\Gamma$ is simply extendable in the vertex $\tau_p$ if and only if
  \begin{equation}
    b^{\vee}_{\tau_p,\tau_p} < 2,
  \end{equation}
   where $b^{\vee}_{\tau_p,\tau_p}$ is the diagonal element (corresponding to the vertex $\tau_p$) of the matrix $B_L^{-1}$.
 \end{proposition}

 \PerfProof This is a direct consequence of Theorem \ref{th_B_less_2}. Indeed, the linkage labels vector
 $\gamma^{\vee}$ corresponding to the simply extendable vertex $\tau_p$ is the vector with the unit in the place $\tau_p$
 and zeros in remaining places. Then $\mathscr{B}^{\vee}_L(\gamma^{\vee}) = b^{\vee}_{\tau_p,\tau_p}$. \qed

\begin{figure}[H]
\centering
\includegraphics[scale=0.7]{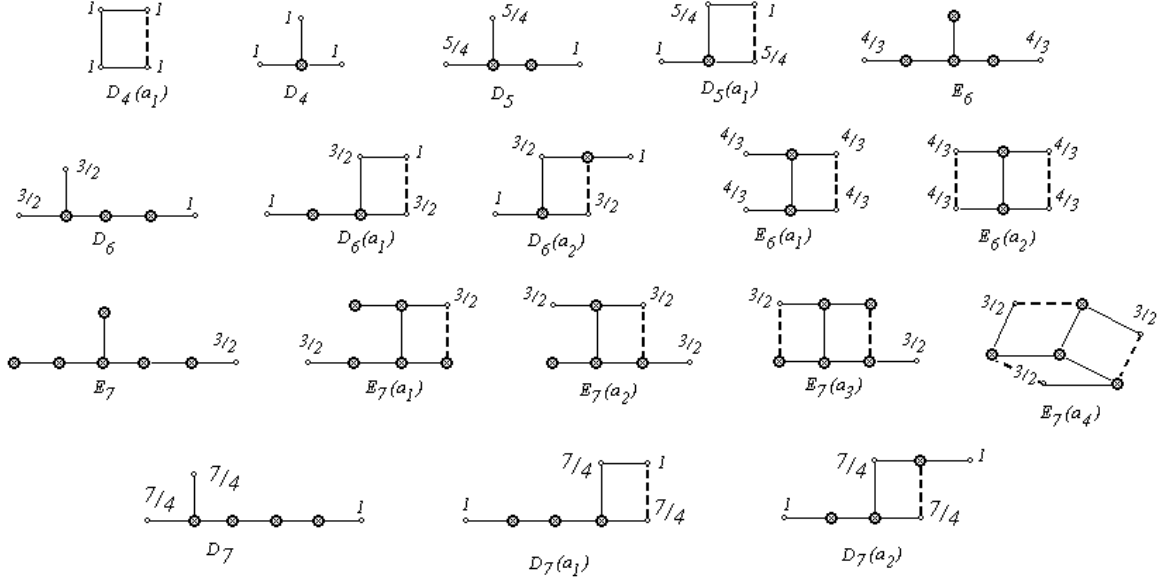}
\caption{
\hspace{3mm}Simply extendable vertices are marked by numerical values, which are the diagonal elements
of $B^{-1}_L$. The vertices marked in bold are not simply extendable.}
\label{examples_ext}
\end{figure}

\begin{remark}{\rm
  1) Not every Carter diagram is an simple extension of any other Carter diagram.
  For example, $E_6(a_2)$ and $E_7(a_4)$ are not extension of any Carter diagram.

 2) For $l < 8$, Proposition \ref{prop_crit_sim_ext} can be also checked by means of Tables \ref{tab_partial Cartan_1}-\ref{tab_partial Cartan_3}.
 For $l \geq 8$, Proposition \ref{prop_crit_sim_ext} can be also derived from Proposition \ref{prop_determinants} and Remark \ref{rem_val_le_2}.
  }
\end{remark}

\subsection{Theorem on projection of linkage systems}
  Let $\mathscr{L}(\Gamma)$  the linkage system associated with the Carter diagram $\Gamma$.
  In the case, where the Carter diagram $\widetilde\Gamma$ is an
  simple extension of other Carter diagram $\Gamma$, we construct the projection of
  the linkage diagrams $\widetilde\gamma^{\vee} \in \mathscr{L}(\widetilde\Gamma)$
  to the linkage diagrams $\gamma^{\vee} \in \Gamma$.

\begin{theorem}
  \label{th_map_linkage_systems}
  Let the Carter diagram $\widetilde\Gamma$ be the simple extension of the Carter diagram $\Gamma$ in the vertex $\tau_p$,
  such that the vertex $\tau_p$ is connected to the vertex $\tau_{l+1}$.
  Let $\widetilde\gamma^{\vee}$ be a certain linkage diagram for $\widetilde\Gamma$, and $\gamma^{\vee}$
  be the vector obtained from $\widetilde\gamma^{\vee}$ by removing the coordinate $\tau_{l+1}$.
  If $\gamma^{\vee} \neq 0$ then $\gamma^{\vee}$ is the linkage diagram for $\Gamma$.
\end{theorem}

   \PerfProof
   {\it Step 1}. According to Theorem \ref{th_B_less_2}, it suffices to prove that
  \begin{equation}
     \label{two_ineq}
     \widetilde{\mathscr{B}}^{\vee}_L(\widetilde\gamma^{\vee}) < 2 \quad
     \Longrightarrow \quad
     \mathscr{B}^{\vee}_L(\gamma^{\vee}) < 2,
  \end{equation}
  where $\widetilde{\mathscr{B}}^{\vee}_L$  the inverse quadratic form associated with $\widetilde\Gamma$.
  The connection between partial Cartan matrices for $\Gamma$ and $\widetilde\Gamma$ is as follows:
  \begin{equation}
    \label{B_relation_1}
     \widetilde{B}_L =
       \left (
         \begin{array}{cc}
           B_L  &   \vec{v} \\
           {}^t\vec{v} & 2 \\
         \end{array}
       \right ),  \text{ where }
       \vec{v} = \left (
         \begin{array}{c}
           0 \\
           \dots \\
           -1 \\
           \dots \\
           0 \\
         \end{array}
         \right )
         \begin{array}{c}
           ~ \\
           ~ \\
           \tau_p \\
           ~ \\
           ~ \\
         \end{array} .
  \end{equation}
   According to \eqref{conn_L_1}, let
   \begin{equation}
    \label{B_relation_2}
       \widetilde{B}_L^{-1}\widetilde\gamma^{\vee} = \widetilde\gamma_L = \left (
         \begin{array}{c}
           \vec{\varphi} \\
           u \\
         \end{array}
         \right )
         \begin{array}{c}
           ~ \\
           \tau_{l+1} \\
         \end{array}         ,
  \end{equation}
  where $\vec\varphi$ a certain vector of the dimension $l$, and $u$ is a rational number. Then
   \begin{equation}
    \label{B_relation_3}
       \widetilde\gamma^{\vee} = \widetilde{B}_L \widetilde\gamma_L =
         \widetilde{B}_L \left (
           \begin{array}{c}
              \vec{\varphi} \\
             u \\
            \end{array}
           \right ) =
       \left (
           \begin{array}{c}
             B_L  \vec{\varphi} + u \vec{v}  \\
             -\varphi_{\tau_p} + 2u \\
             \end{array}
       \right ),
  \end{equation}
   where $\varphi_{\tau_p}$ is the coordinate $\tau_p$  of $\vec{\varphi}$.
  From \eqref{B_relation_2} and \eqref{B_relation_3} we have
  \begin{equation}
    \label{B_relation_4}
    \widetilde{\mathscr{B}}^{\vee}_L(\widetilde\gamma^{\vee}) =
    \langle \widetilde{B}^{-1}_L \widetilde\gamma^{\vee}, \widetilde\gamma^{\vee}   \rangle =
    \langle \widetilde\gamma_L, \widetilde{B}_L \widetilde\gamma_L \rangle =
    \langle {B}_L \vec{\varphi}, \vec{\varphi}  \rangle - 2\varphi_{\tau_p}u + 2u^2.
  \end{equation}
  Thus the property $\widetilde{\mathscr{B}}^{\vee}_L(\widetilde\gamma^{\vee}) < 2$ from \eqref{two_ineq} is
  equivalent to
  \begin{equation}
    \label{B_boxed_1}
    \boxed{\langle {B}_L \vec{\varphi}, \vec{\varphi} \rangle
       - 2\varphi_{\tau_p}u + 2u^2 < 2.}
  \end{equation}

  {\it Step 2}. Now, let $\gamma^{\vee}$ be the vector obtained from $\widetilde\gamma^{\vee}$
  by removing $\tau_{l+1}$, i.e. $\tau_{l+1} = 0$ in $\gamma^{\vee}$.
  By \eqref{B_relation_3}  we have
    \begin{equation}
      \label{B_relation_6}
         \left (
           \begin{array}{c}
              \gamma^{\vee} \\
              0 \\
            \end{array}
           \right )
           \begin{array}{c}
              ~ \\
              \tau_{l+1} \\
            \end{array} =
         \left (
           \begin{array}{c}
              B_L \vec{\delta} \\
              0 \\
            \end{array}
           \right ) =
         \widetilde{B}_L \left (
           \begin{array}{c}
              \vec{\varphi} \\
             u \\
            \end{array}
           \right ) =
       \left (
           \begin{array}{c}
             B_L  \vec{\varphi} + u \vec{v}  \\
             -\varphi_{\tau_p} + 2u \\
             \end{array}
       \right ),
  \end{equation}
  where $\vec{\delta}$ is a ceratin vector of the dimension $l$.

  \begin{equation}
      \label{B_relation_7}
        \vec{\delta} =
           B_L^{-1}(B_L\vec{\varphi} + u\vec{v}) = \vec{\varphi} + u{B}_L^{-1}\vec{v}.
  \end{equation}
  By \eqref{B_relation_1}   ${B}_L^{-1}\vec{v} = -\vec{b}_{\tau_p}$, where $\vec{b}_{\tau_p}$ is the
  $\tau_p$-th column of ${B}_L^{-1}$. Since $\langle u \vec{v}, \vec{\varphi} \rangle =  -u\varphi_{\tau_p}$,
  we have
    \begin{equation}
      \label{B_relation_8}
        \begin{split}
        \mathscr{B}^{\vee}_L(\gamma^{\vee}) =
        \langle \gamma^{\vee}, B_L^{-1}\gamma^{\vee} \rangle = &
        \langle B_L\vec{\delta}, \vec{\delta} \rangle =
        \langle B_L  \vec{\varphi} + u \vec{v}, \vec{\varphi} - u\vec{b}_{\tau_p} \rangle = \\
        & \langle B_L  \vec{\varphi}, \vec{\varphi} \rangle - u\varphi_{\tau_p} -
          \langle B_L \vec{\varphi}, u\vec{b}_{\tau_p} \rangle  + u^2(\vec{b}_{\tau_p})_{\tau_p} = \\
        & \langle B_L  \vec{\varphi}, \vec{\varphi} \rangle - u\varphi_{\tau_p} -
          u\langle \vec{\varphi}, B_L\vec{b}_{\tau_p} \rangle  + u^2(\vec{b}_{\tau_p})_{\tau_p},
        \end{split}
  \end{equation}
 where $(\vec{b}_{\tau_p})_{\tau_p}$ is the $(\tau_p, \tau_p)$ slot of the matrix $B_L^{-1}$, i.e.,
    \begin{equation}
    \begin{split}
      \label{B_relation_9}
        \langle B_L\vec{\delta}, \vec{\delta} \rangle = &
         \langle B_L  \vec{\varphi}, \vec{\varphi} \rangle - u\varphi_{\tau_p} -
          u\langle \vec{\varphi}, B_L\vec{b}_{\tau_p} \rangle  + u^2{({B}_L^{-1})}_{\tau_p, \tau_p} = \\
         & \langle B_L  \vec{\varphi}, \vec{\varphi} \rangle - u\varphi_{\tau_p} -
          u\langle \vec{\varphi}, -\vec{v} \rangle  + u^2{({B}_L^{-1})}_{\tau_p, \tau_p} = \\
         & \langle B_L  \vec{\varphi}, \vec{\varphi} \rangle - 2u\varphi_{\tau_p}
           + u^2{({B}_L^{-1})}_{\tau_p, \tau_p}. \\
    \end{split}
  \end{equation}
  Thus the property $\mathscr{B}^{\vee}_L(\gamma^{\vee}) < 2$ from \eqref{two_ineq} is
  equivalent to
  \begin{equation}
    \label{B_boxed_2}
    \boxed{\langle {B}_L \vec{\varphi}, \vec{\varphi} \rangle
       - 2\varphi_{\tau_p}u + u^2{({B}_L^{-1})}_{\tau_p, \tau_p} < 2.}
  \end{equation}
  Since $\widetilde\Gamma$ is simply extendable in the vertex $\tau_p$, then by Proposition \ref{prop_crit_sim_ext}  we get ${({B}_L^{-1})}_{\tau_p, \tau_p} < 2$. Then inequality \eqref{B_boxed_2} follows from inequality \eqref{B_boxed_1}  and \eqref{two_ineq} is proven. \qed

  \begin{remark}
     \label{rem_mapping}{\rm
  1) If $\gamma^{\vee}$ obtained from $\widetilde\gamma^{\vee}$ by removing the coordinate $\tau_{l+1}$ is the
  linkage diagram for $\Gamma$, we construct the linkage diagram projection from
  the linkage system of $\widetilde\Gamma$ to the linkage system of $\Gamma$:
  \begin{equation}
    \label{f_map_1}
     f : \widetilde\gamma^{\vee} \longrightarrow \gamma^{\vee}.
  \end{equation}

  2) If $\gamma^{\vee}$ (resp. $\widetilde\gamma^{\vee}$) belongs to any loctet $L$ (resp. $\widetilde{L}$) of the linkage system $\mathscr{L}(\Gamma)$ (resp. $\mathscr{L}(\widetilde\Gamma)$) associated with the Carter diagram $\Gamma$ (resp. $\widetilde\Gamma$)
  then the projection \eqref{f_map_1} is extended to the projection of the loctet $\widetilde{L}$ onto the loctet $L$,  see Fig. \ref{inclusions_E6_E7} and Fig. \ref{inclusions_D6a_D5a1}. 
  It follows from \eqref{dual_refl} and the fact that $s_{\alpha_i}$, where $i = 1,2,3$,
  and $s_{\beta_1}$ act by the same way on all coordinates of $\gamma^{\vee}$ and 
  $\widetilde\gamma^{\vee}$ except the coordinate $\tau_{l+1}$.

  3) Two loctets of $\mathscr{L}(\widetilde\Gamma)$ can be mapped by $f$ onto the same loctet of $\mathscr{L}(\Gamma)$.
     For example, the following pairs of loctets of $D_6(a_1)$ are mapped onto the same loctet of $D_5(a_1)$.
  \begin{equation}
    \label{mapping_to_same}
      \begin{split}
       f : L^d_{13},  L^c_{13}  \longrightarrow L^b_{13}
                  \quad (L^d_{13},  L^c_{13} \subset D_6(a_1); L^b_{13}  \subset D_5(a_1)), \\
       f : L^d_{12},  L^c_{12}  \longrightarrow L^b_{12}
                  \quad (L^d_{12},  L^c_{12} \subset D_6(a_1); L^b_{12}  \subset D_5(a_1)), \\
       f : L^a_{12},  L^b_{12}  \longrightarrow L^a_{12}
                   \quad (L^a_{12},  L^b_{12} \subset D_6(a_1); L^a_{12}  \subset D_5(a_1)), \\
       f : L^a_{13},  L^b_{13}  \longrightarrow L^a_{13}
                    \quad (L^a_{13},  L^b_{13} \subset D_6(a_1); L^a_{13}  \subset D_5(a_1)), \\
      \end{split}
  \end{equation}
  see Fig. \ref{inclusions_D6a_D5a1}, Fig. \ref{D5a1_linkages} and Fig. \ref{D6a1_linkages}.
 \begin{figure}[H]
\centering
\includegraphics[scale=0.5]{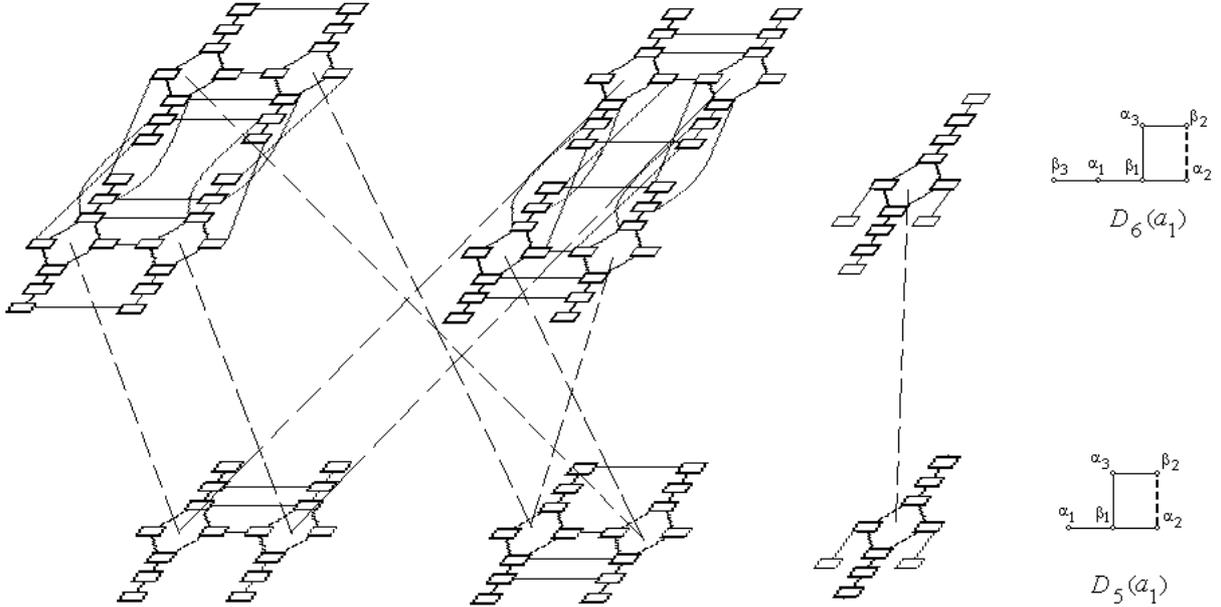}
\caption{
\hspace{3mm}There are pairs of $D_6(a_1)$ loctets, which are mapped (dotted lines) onto the same $D_5(a_1)$ loctet}
\label{inclusions_D6a_D5a1}
\end{figure}

  4) Let $\widetilde\gamma^{\vee}_1$ (resp. $\widetilde\gamma^{\vee}_2$)
  be a linkage diagram of $\widetilde\Gamma$ such that the coordinate $\tau_{l+1} = 1$ (resp. $\tau_{l+1} = -1$),
  and all remaining coordinates are zero. Vectors $\widetilde\gamma^{\vee}_{1,2}$ constitute $2$-element kernel of
  the projection \eqref{f_map_1}. That is why the third component of $D_l(a_k)$ contains by $2$ linkage diagrams more than $D_{l-1}(a_k)$, see Fig. \ref{inclusions_D6a_D5a1} (they are the top and bottom linkage diagrams in the component).
  }
  \end{remark}

\subsection{The linkage systems $D_l(a_k), D_l$ for $l > 7$}
 We observe that $D$-type components of linkage systems of
 $D_5(a_1)$,  $D_6(a_1)$, $D_6(a_2)$, $D_7(a_1)$, $D_7(a_2)$,
 see Fig. \ref{D5a1_linkages}, Fig. \ref{D6a1_linkages},  Fig. \ref{D6a2_linkages},
 Fig. \ref{D7a1_D7a2_D7pu_loctets_comp3} have the same shape.
 We will show that the linkage systems $D_l(a_k)$, where $l > 7$, are of the same shape which we call
 the {\it wind rose of linkages}, see Fig. \ref{Dk_al_wind_rose}.
 Similarly, $D$-type components of linkage systems of $D_5$, $D_6$, $D_7$,
 see Fig. \ref{D5pure_loctets},  Fig. \ref{D6pure_loctets},
 Fig. \ref{D7a1_D7a2_D7pu_loctets_comp3}  have the same shape.
 We will show that the linkage system $D_l$, where $l > 7$, are of the same shape,
 see Fig. \ref{Dlpu_linkages}.

 \begin{lemma}[on moving triangles]
   \label{lem_triangles}
   1) Let the diagram $\Gamma$ contain the chain $\{ \alpha_1, \beta_1, \alpha_2, \beta_2, ... \}$,
   ($\Gamma$ may be $D_l(a_k), D_l, \dots$) and let $\gamma^{\vee}$ be such a linkage diagram
   that only $3$ coordinates are non-zero, namely $(\gamma^{\vee})_{\alpha_1} \neq 0$ and two other non-zero
   coordinates together with $\gamma$ constitute a triangle as in Fig. \ref{elim_triangle},(a).
   Let us consider linkage diagrams obtained by means of sequential transformations as follows:
   \begin{equation}
     \label{eq_transf_tr}
       \begin{split}
        & \gamma_b = s_{\alpha_k}(\gamma) = \gamma - \alpha_k, \\
        & \gamma_c = s_{\beta_{k-1}}(\gamma_b) = \gamma - \alpha_k - \beta_{k-1}, \\
        & \dots, \\
        & \gamma_d = s_{\alpha_2}(\gamma_c - \dots - \beta_2) = \gamma_c - \dots - \beta_2 - \alpha_2, \\
        & \gamma_e = s_{\beta_1}(\gamma_c - \dots - \beta_2 - \alpha_2) = \gamma_c - \dots - \beta_2 - \alpha_2 - \beta_1.
       \end{split}
   \end{equation}
   Transformations \eqref{eq_transf_tr} preserves the inverse quadratic form $\mathscr{B}^{\vee}$
   on the corresponding linkage vectors:
 \begin{equation}
    \gamma_a^{\vee}, \gamma_b^{\vee}, \gamma_c^{\vee}, \dots, \gamma_d^{\vee}, \gamma_e^{\vee}.
 \end{equation}

   2) By transformations \eqref{eq_transf_tr},
   the triangle (describing the linkage diagram $\gamma^{\vee}$) is shifted to the left:
   $(a) \Rightarrow (b) \Rightarrow (c) \Rightarrow \dots \Rightarrow (d)$, see Fig. \ref{elim_triangle},(b),(c),(d).
   In the last step $(d) \Rightarrow (e)$, the triangle is eliminated, see Fig. \ref{elim_triangle},(e).

   3) Linkage diagrams depicted in Fig. \ref{elim_triangle} can not occur for $\Gamma = D_l(a_k), D_l$.

\begin{figure}[H]
\centering
\includegraphics[scale=0.55]{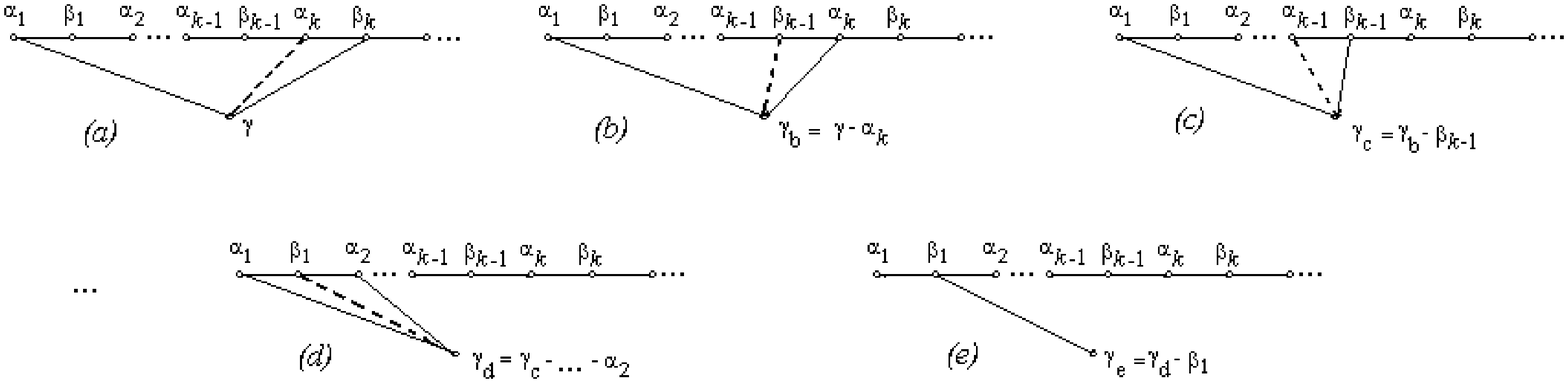}
\caption{
\hspace{3mm}The moving triangles}
\label{elim_triangle}
\end{figure}

    4) Let $\gamma^{\vee}$ be such a linkage diagram for $\Gamma = D_l(a_k)$ (resp. $\Gamma = D_l$)
   that only $4$ coordinates are non-zero, namely, $(\gamma^{\vee})_{\alpha_1} \neq 0$ and three other non-zero
   coordinates together with $\gamma_f$ (resp. $\gamma_h$) constitute two triangles as in Fig. \ref{elim_triangle_2},(f) (resp. \ref{elim_triangle_2},(h)).
   Let us consider linkage diagram obtained by means of transformation $s_{\beta_1}$ as follows:
   \begin{equation}
     \label{eq_transf_tr_2}
      \gamma_g = s_{\beta_1}(\gamma_f),  \quad (resp. \quad \gamma_i = s_{\beta_1}(\gamma_h)) \\
   \end{equation}
   Transformations \eqref{eq_transf_tr_2} preserves the inverse quadratic form $\mathscr{B}^{\vee}$
   on the corresponding linkage vectors:
 \begin{equation}
    \gamma_f^{\vee}, \gamma_g^{\vee},  \quad (resp. \quad \gamma_h^{\vee}, \gamma_i^{\vee}).
 \end{equation}

\begin{figure}[H]
\centering
\includegraphics[scale=0.6]{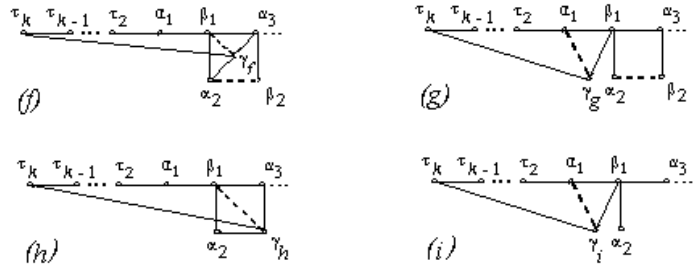}
\caption{
\hspace{3mm}The moving triangles, II}
\label{elim_triangle_2}
\end{figure}
   Linkage diagrams depicted in Fig. \ref{elim_triangle_2},(f),(g) (resp. Fig. \ref{elim_triangle_2},(h),(i)) can not occur for $\Gamma = D_l(a_k)$ (resp. $D_l$).
 \end{lemma}

   \PerfProof
   1) It follows from Proposition \ref{prop_link_diagr_conn}.

   2) Let us check $(b) \Rightarrow (c)$, where
   $\gamma_c = s_{\beta_{k-1}}(\gamma_b) = \gamma_b - \beta_{k-1}$:
  \begin{equation*}
        (\gamma_c, \alpha_k) = -1 + 1 = 0, \quad
        (\gamma_c, \beta_{k-1}) = 1 - 2 = -1, \quad
        (\gamma_c, \alpha_{k-1})  = 0 + 1 = 1,
   \end{equation*}
  remaining inner products $(\gamma_c, \tau)$ are not changed. 
  In the case (e), we have $s_{\gamma_e}(\gamma_d) = \gamma_d - \beta_1$, and
  \begin{equation*}
        (\gamma_e, \alpha_2) =  1 - 1 = 0, \quad
        (\gamma_e, \beta_1)  =  1 - 1 = 0, \quad
        (\gamma_e, \alpha_1) =  1 - 2 = -1.
  \end{equation*}
  remaining inner $(\gamma_e, \tau)$ are not changed.

  3) Suppose one of linkage diagrams in Fig. \ref{elim_triangle} occurs for $\Gamma = D_l(a_k), D_l$.
     Then, in particular, $\gamma_e$ from  Fig. \ref{elim_triangle} is the valid linkage diagram for $\Gamma$.
     Hence, $\Gamma$ has the simple extension in $\beta_1$. According to Proposition \ref{prop_crit_sim_ext}, 
     it can be only if $b^{\vee}_{\beta_1, \beta_1} < 2$. However, by Proposition \ref{prop_determinants},
     we have $b^{\vee}_{\beta_1, \beta_1} = 2$.

  4) Let us check $(f) \Rightarrow (g)$, where
   $\gamma_g = s_{\beta_1}(\gamma_f) = \gamma_f - \beta_1$:
  \begin{equation*}
        (\gamma_g, \alpha_2) = 1 - 1 = 0, \quad
        (\gamma_g, \alpha_3) = 1 - 1 = 0, \quad
        (\gamma_g, \beta_1)  = 1 - 2 = -1,
        (\gamma_g, \tau_k)  = (\gamma_f, \tau_k) = -1,
   \end{equation*}
  remaining inner products $(\gamma_g, \tau)$ are not changed. Thus, we obtain the triangle
   $\{\alpha_1, \beta_1, \gamma_g \}$, i.e, we come to one of cases in
  Fig. \ref{elim_triangle}. As in heading 3) we get simple extension in $\beta_1$, that can not be.
\qed

 \begin{proposition}
  \label{prop_Dl_ak}
    1) Four vectors $\gamma^{\vee}_i$  \eqref{endp_linkages} (each of them contains
       the single non-zero integer number in the slot $i$)  are the linkage labels vectors for the $D_l(a_k)$.
      \begin{equation}
        \label{endp_linkages}
        \gamma^{\vee}_{\tau^{+}_{k-1}} =
        \left (
         \begin{array}{c}
            0 \\
            \dots \\
            1 \\
            \dots \\
            0 \\
         \end{array}
        \right ),
        \quad
      \gamma^{\vee}_{\tau^{-}_{k-1}} =
        \left (
         \begin{array}{c}
            0 \\
            \dots \\
            -1 \\
            \dots \\
            0 \\
         \end{array}
        \right ),
         \quad
      \gamma^{\vee}_{\varphi^{+}_{l - k - 2}} =
        \left (
         \begin{array}{c}
            0 \\
            \dots \\
            1 \\
            \dots \\
            0 \\
         \end{array}
        \right ),
      \quad
      \gamma^{\vee}_{\varphi^{-}_{l - k - 2}} =
        \left (
         \begin{array}{c}
            0 \\
            \dots \\
            -1 \\
            \dots \\
            0 \\
         \end{array}
        \right ).
      \end{equation}

 2) Linkage diagrams for the general case $D_l(a_k)$, where $l > 7$, are presented in Fig. \ref{Dk_al_linkages}.
    The linkage system $D_l(a_k)$, where $l > 7$, is depicted in Fig. \ref{Dk_al_wind_rose}.

 3) Linkage diagrams and the linkage system $D_l$ for $l > 7$, are presented in Fig. \ref{Dlpu_linkages}.
 \end{proposition}

 \PerfProof

  1) According to Proposition \ref{prop_determinants}, heading 2),
  for $\gamma^{\vee}_{\tau} = \gamma^{\vee}_{\tau^{\pm}_{k-1}},  \gamma^{\vee}_{\varphi^{\pm}_{l - k - 2}}$,
  we have
 \begin{equation}
   \label{eq_1_Dlak}
   \mathscr{B}^{\vee}_L(\gamma^{\vee}_{\tau}) = 1.
 \end{equation}
  Eq. \eqref{eq_1_Dlak} means that $\mathscr{B}^{\vee}_L(\gamma^{\vee}_{\tau}) < 2$, and
  by Theorem \ref{th_B_less_2} we get that vectors $\gamma^{\vee}_{\tau}$ are linkage labels vectors.

  2) Further, we apply $s^{*}_{i}$ as depicted in Fig. \ref{Dk_al_wind_rose}.
     This way we obtain all linkage diagrams in Fig. \ref{Dk_al_linkages}.
     The linkage diagrams in Fig. \ref{Dk_al_linkages} constitute the linkage
     system depicted in Fig. \ref{Dk_al_wind_rose}. 
     We show that there are no other linkage diagrams, except those which
     are shown in Fig. \ref{Dk_al_linkages}. Assume, that the statement already shown
     for $\Gamma = D_l(a_k)$, and we consider $\widetilde\Gamma = D_{l+1}(a_k)$ (or, similarly, for $D_l(a_{k+1})$).
     Let $\widetilde\gamma^{\vee}$ be the linkage diagram from the linkage system of $\widetilde\Gamma$,     
     and $\widetilde\gamma^{\vee}_{\tau}$ be the coordinate corresponding to the simply extendable vertex $\tau$ of 
     $\widetilde\Gamma$, see Remark \ref{rem_val_le_2}. 
     Let us consider two cases $a) \widetilde\gamma^{\vee}_{\tau} = 0 $ and $b)
     \widetilde\gamma^{\vee}_{\tau} \neq 0$.

  a) Case $\widetilde\gamma^{\vee}_{\tau} = 0$. Suppose, that the linkage diagram $\widetilde\gamma^{\vee}$ is outside  of the list of Fig. \ref{Dk_al_linkages}. Removing the connection
  $\{ \gamma, \tau \}$ does not change the geometric form of  $\gamma^{\vee}$ (since $(\gamma, \tau)$ = 0).
  Therefore, $\gamma^{\vee}$ is also outside of the list of Fig. \ref{Dk_al_linkages}.
  According to Theorem \ref{th_map_linkage_systems}, we have that $\gamma^{\vee}$ is the linkage diagram for $\Gamma$,
  that contradicts to the induction hypothesis, that $\widetilde\gamma^{\vee}$ is necessarily from
  the list of Fig. \ref{Dk_al_linkages}.

  b) Case $\widetilde\gamma^{\vee}_{\tau} \neq 0$.
    Let us remove the connection $\{ \gamma, \tau \}$, we get one of the following cases:

  \quad b.1) $\gamma = 0$, then $\widetilde\gamma$ is simple extension in the vertex $\tau$,
      $\widetilde\gamma = \widetilde\gamma^{\vee}_{\tau^{\pm}_k}$ for the case $\widetilde\Gamma = D_{l}(a_{k+1})$,
      or $\widetilde\gamma = \widetilde\gamma^{\vee}_{\varphi^{\pm}_{l-k-1}}$ for the case $\widetilde\Gamma = D_{l+1}(a_k)$, see Fig. \ref{Dk_al_linkages}.

   \quad b.2) $\gamma$ is one of triangles in Fig. \ref{Dk_al_linkages}. Since $(\gamma, \tau ) \neq 0$, then
     $\widetilde\gamma$ is one of the cases of Lemma \ref{lem_triangles} and in Fig. \ref{elim_triangle}
     that can not occur.

   \quad b.3) $\gamma$ is one linkage diagrams $3,4,5,6$ in in Fig. \ref{Dk_al_linkages}.
      Since $(\gamma, \tau ) \neq 0$, then
      $\widetilde\gamma$ is the linkage diagram as in Lemma \ref{lem_triangles}, heading 4) and
      Fig. \ref{elim_triangle_2},(f),(g) that can not occur.

  3) We repeat here arguments of the heading 2). In the case b.3), $\widetilde\gamma$ is the linkage diagram as in Lemma \ref{lem_triangles}, heading 4) and Fig. \ref{elim_triangle_2},(h),(i) that can not occur.
   \qed

%% file: A4matricesB.tex
\newpage
\section{\sc\bf The inverse matrix $B^{-1}_L$, linkage diagrams $\gamma^{\vee}(8)$ and inequality $\mathscr{B}^{\vee}_L(\gamma^{\vee}) < 2$}
   \label{sec_inverse}

\subsection{The partial Cartan matrix $B_L$ and the inverse matrix $B^{-1}_L$}
  \label{sec_inv_matr}
~\\
 \begin{table}[H]
\tiny
  \centering
  \renewcommand{\arraystretch}{1.5}
\begin{tabular} {|c|c|c|}
  \hline
      The Carter   & The partial Cartan         &   The inverse          \cr
        diagram    &  matrix $B_L$              &   matrix $B^{-1}_L$    \\
  \hline  
     & & \\
     $\begin{array}{c} \includegraphics[scale=0.6]{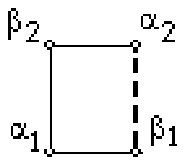} \\
        {\bf D_4(a_1)} \end{array}$  &
     $\begin{array}{c}
       \left [
   \begin{array}{cccc}
     2 & 0 &   -1 & -1  \\
     0 & 2 &   1 & -1 \\
     -1 & 1 &  2 & 0  \\
     -1 & -1 &   0 & 2  \\
  \end{array}
  \right ]
     \end{array}$   &
    $\begin{array}{c}
    \frac{1}{2} \left [ \begin{array}{cccc}
     2 & 0 &   1 & 1  \\
     0 & 2 &   -1 & 1 \\
     1 & -1 &  2 & 0  \\
     1 &  1 &   0 & 2  \\
  \end{array} \right ]
     \end{array}$    \\
    & & \\
    \hline
    & & \\
     $\begin{array}{c} \includegraphics[scale=0.6]{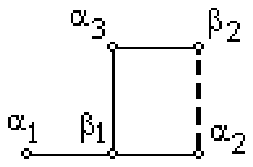} \\
        {\bf D_5(a_1) = D_5(a_2)} \end{array}$ &
     $\begin{array}{c}
       \left [
   \begin{array}{ccccc}
     2 & 0 &   0 & -1 & 0  \\
     0 & 2 &   0 & -1 & 1\\
     0 & 0 &   2 & -1 & -1  \\
     -1 & -1 &  -1 & 2 & 0  \\
      0 &  1 &   -1  & 0 & 2 \\
  \end{array}
  \right ]
     \end{array}$   &
      $\begin{array}{c}
    \frac{1}{4} \left [ \begin{array}{ccccc}
     4 & 2 &   2 & 4 & 0  \\
     2 & 5 &   1 & 4 & -2\\
     2 & 1 &   5  & 4 & 2  \\
     4 & 4 &   4 & 8 & 0  \\
     0 & -2 &  2  & 0 & 4 \\
  \end{array} \right ]
     \end{array}$    \\
    & & \\
\hline
\end{tabular}
  \vspace{2mm}
  \caption{\small\hspace{3mm} The partial Cartan matrix $B_L$ and the inverse matrix $B^{-1}_L$ for Carter
   diagrams with the number of vertices $l < 7$}
  \label{tab_partial Cartan_1}
\end{table}

 \begin{table}[H]
\tiny
  \centering
  \renewcommand{\arraystretch}{1.4}
  \begin{tabular} {|c|c|c|}
  \hline
      The Carter   & The partial Cartan         &   The inverse          \cr
        diagram    &  matrix $B_L$              &   matrix $B^{-1}_L$    \\
     \hline
     & & \\
          $\begin{array}{c} \includegraphics[scale=0.6]{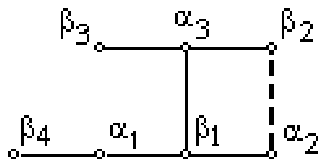} \\
        {\bf E_7(a_1)} \end{array}$  &
 $\begin{array}{c}
      \left [
   \begin{array}{ccccccc}
      2 & 0 &   0 & -1 & 0  & 0 & -1 \\
     0 & 2 &   0 & -1 & 1   & 0 & 0 \\
     0 & 0 &   2 & -1 & -1  & -1 & 0 \\
     -1 & -1 &  -1  & 2 & 0 & 0 & 0 \\
      0 &  1 &  -1  & 0 & 2 & 0 & 0 \\
      0 &  0 &  -1  & 0 & 0 & 2 & 0 \\
     -1 &  0 &  0   & 0 & 0 & 0 & 2 \\
  \end{array}
  \right ]
     \end{array}$   &
   $\begin{array}{c}
      \frac{1}{2}\left [
   \begin{array}{ccccccc}
      8 & 4 &  8 &  10 & 2 & 4 & 4 \\
      4 & 4 &  4 &  6 &  0 & 2 & 2 \\
      8 & 4 & 12 &  12 & 4 & 6 &  4 \\
      10 & 6 & 12 & 15 & 3 & 6  & 5 \\
      2 &  0 & 4 &  3 & 3 & 2  & 1 \\
      4 & 2 &  6 &  6 & 2 & 4  & 2 \\
      4 & 2 &  4 &  5 & 1 & 2  & 3 \\
  \end{array}
  \right ]
     \end{array}$   \\
     & & \\
     \hline
    & & \\
     $\begin{array}{c} \includegraphics[scale=0.6]{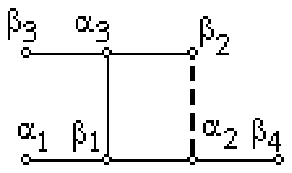} \\
        {\bf E_7(a_2)}
     \end{array}$ &
     $\begin{array}{c}
      \left [
   \begin{array}{ccccccc}
     2 & 0 &   0 & -1 &  0  & 0 & 0 \\
     0 & 2 &   0 & -1 &  1  & 0 & -1 \\
     0 & 0 &   2 & -1 & -1  & -1 & 0 \\
     -1 & -1 &  -1  & 2 & 0 & 0 & 0 \\
      0 &  1 &  -1  & 0 & 2 & 0 & 0 \\
      0 &  0 &  -1  & 0 & 0 & 2 & 0 \\
      0 & -1 &   0  & 0 & 0 & 0 & 2 \\
  \end{array}
  \right ]
     \end{array}$   &
     $\begin{array}{c}
      \frac{1}{2}\left [
   \begin{array}{ccccccc}
      4 & 4 &  4 &  6 & 0 & 2 & 2 \\
      4 & 8 &  4 &  8 & -2 & 2 & 4 \\
      4 & 4 &  8 &  8 & 2 & 4 &  2 \\
      6 & 8 & 8 &  12 & 0 & 4  & 4 \\
      0 & -2 & 2 &  0 & 3 & 1  & -1 \\
      2 &  2 & 4 &  4 & 1 & 3  & 1 \\
      2 &  4 & 2 &  4 & -1 & 1  & 3 \\
  \end{array}
  \right ]
     \end{array}$   \\
\hline
     \hline
\end{tabular}
  \vspace{2mm}
  \caption{\small\hspace{3mm} (cont.) The partial Cartan matrix $B_L$ and the inverse matrix $B^{-1}_L$ for Carter
   diagrams with the number of vertices $l = 7$}
  \label{tab_partial Cartan_2}
\end{table}

 \begin{table}[H]
\tiny
  \centering
  \renewcommand{\arraystretch}{1.5}
  \begin{tabular} {|c|c|c|}
  \hline
      The Carter   & The Cartan         &   The inverse          \cr
        diagram    &  matrix ${\bf B}$              &   matrix ${\bf B}^{-1}$    \\
  \hline
      & & \\
     $\begin{array}{c} \includegraphics[scale=0.6]{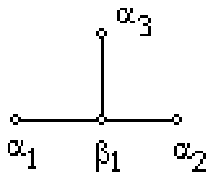} \\
        {\bf D_4} \end{array}$ &
     $\begin{array}{c}
       \left [
   \begin{array}{cccc}
     2 & 0 & 0  &   -1  \\
     0 & 2 & 0  &   -1  \\
     0 & 0 & 2  &   -1  \\
    -1 & -1 & -1 &   2  \\
  \end{array}
  \right ]
     \end{array}$   &
      $\begin{array}{c}
    \frac{1}{2} \left [ \begin{array}{cccc}
     2 & 1 & 1 & 2  \\
     1 & 2 & 1 & 2  \\
     1 & 1 & 2 & 2  \\
     2 & 2 & 2 & 4  \\
  \end{array} \right ]
     \end{array}$    \\
    & & \\
  \hline  
    & & \\
     $\begin{array}{c} \includegraphics[scale=0.6]{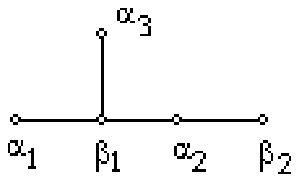} \\
        {\bf D_5} \end{array}$ &
     $\begin{array}{c}
       \left [
   \begin{array}{ccccc}
     2 & 0 & 0  &   -1 &  0 \\
     0 & 2 & 0  &   -1 & -1 \\
     0 & 0 & 2  &   -1 &  0 \\
    -1 & -1 & -1 &   2 &  0 \\
     0 & -1 & 0  &   0 & 2  \\
  \end{array}
  \right ]
     \end{array}$   &
      $\begin{array}{c}
    \frac{1}{4} \left [ \begin{array}{ccccc}
     5 & 4 & 3 & 6 & 2  \\
     4 & 8 & 4 & 8 & 4  \\
     3 & 4 & 5 & 6 & 2  \\
     6 & 8 & 6 & 12 & 4 \\
     2 & 4 & 2 & 4 & 4 \\
  \end{array} \right ]
     \end{array}$    \\
    & & \\
\hline
\end{tabular}
  \vspace{2mm}
  \caption{\small\hspace{3mm} (cont.) The Cartan matrix ${\bf B}$ and the inverse matrix ${\bf B}^{-1}$ for the conjugacy classes $E_6$, $E_7$, $D_4$, $D_5$, $D_6$, $D_7$}
  \label{tab_partial Cartan_3}
\end{table}

%% file: A5solutions.tex
\newpage
\subsection{The linkage diagrams $\gamma^{\vee}_{ij}(8)$ and solutions of inequality $\mathscr{B}^{\vee}_L(\gamma^{\vee}_{ij}(8)) < 2$}
~\\

\begin{table}[H]
  \centering
  \renewcommand{\arraystretch}{1.4}
  \small
  \begin{tabular} {||c|c|c|l|l||}
  \hline \hline
     The Carter &  $L_{ij}$  & $\mathscr{B}^{\vee}_L(\gamma^{\vee}_{ij}(8)) = p < 2$
                                & The linkage diagram  & $p$ \cr
      diagram &    &            &  \qquad \qquad $\gamma^{\vee}_{ij}(8)$  &  \\
    \hline  
                  &  &   $q = 4(b_2^2 + b_2{b_3} + b_3^2)$ & & \cr
                  \cline{3-3}
                  &  $L_{12}$
                  &  $\frac{1}{3}(10 + 2(4b_2 + 5b_3) + q) < 2$, or
                  &  $\{ 0, 0, 1, 0, 0, -1 \}$  &  \cr
                  &
                  &  $(b_2 + b_3)^2 + (b_2 + 2)^2  + (b_3 + \frac{5}{2})^2 < \frac{33}{4}$
                  &  $\{ 0, 0, 1, 0, -1, -1 \}$ & \cr
                  \cline{2-4}
    $E_6(a_1)$    &  $L_{13}$
                  &  $\frac{1}{3}(4 + 2(-b_2 + b_3) + q) < 2$, or
                  &  $\{ 0, 1, 0, 0, 0, 0 \}$  & $\frac{4}{3}$ \cr
                  &
                  &  $(b_2 + b_3)^2 + (b_2 - \frac{1}{2})^2  + (b_3 + \frac{1}{2})^2 < \frac{3}{2}$
                  &  $\{ 0, 1, 0, 0,  1, -1 \}$ &   \cr
                  \cline{2-4}
                  &  $L_{23}$
                  &  $\frac{1}{3}(4 + 2(b_2 + 2b_3) + q) < 2$,
                  &  $\{ 1, 0, 0, 0, 0, 0 \}$  &  \cr
\end{tabular}
  \vspace{2mm}
  \caption{\hspace{3mm}The linkage diagrams $\gamma^{\vee}_{ij}(8)$ obtained as solutions of inequality
    $\mathscr{B}^{\vee}_L(\gamma^{\vee}_{ij}(8)) < 2$}
  \label{sol_inequal_1}
  \end{table}

\begin{table}[H]
  \centering
  \renewcommand{\arraystretch}{1.6}
  \small
  \begin{tabular} {||c|c|c|l|l||}
  \hline \hline
      The Carter &  $L_{ij}$  & $\mathscr{B}^{\vee}_L(\gamma^{\vee}_{ij}(8)) = p < 2$
                                & The linkage diagram  & $p$ \cr
      diagram &    &            &  \qquad \qquad $\gamma^{\vee}_{ij}(8)$  &  \\
    \hline  
                  &  &   $q = 4b_2^2 + 3b_3^2 + 3b_4^2 + 4b_2{b_3} - 4b_2{b_4} - 2b_3{b_4}$ & & \cr
                  \cline{3-3}
                  &  $L_{12}$
                  &  $\frac{1}{2}(4 + 2(2b_2 + 2b_3) + q) < 2$, or
                  &  $\{ 0, 0, 1, 0, 0, -1, 0 \}$  &  \cr
                  &
                  &  $2(b_2 + b_3 + 1)^2 + (b_3 - b_4)^2 + 2(b_4 - b_2)^2 < 2$
                  &  $\{ 0, 0, 1, 0, -1, 0, -1 \}$ & \cr
                  \cline{2-4}
    $E_7(a_3)$    &  $L_{13}$
                  &  $\frac{1}{2}(8 + 2(-4b_2 -2b_3 + 4b_4) + q) < 2$, or
                  &  $\{ 0, 1, 0, 0, 0, 0, -1 \}$  & $\frac{3}{2}$ \cr
                  &
                  &  $2(b_2 + b_3 - \frac{1}{2})^2 +
                       (b_3 - b_4 - 1)^2  + 2(b_4 -b_2 + \frac{3}{2})^2 < 2$
                  &  $\{ 0, 1, 0, 0, 1,  0, -1 \}$ &   \cr
                  \cline{2-4}
                  &  $L_{23}$
                  &  $\frac{1}{2}(4 + 2(-2b_2 - 2b_3 + 2b_4) + q) < 2$, or
                  &  $\{ 1, 0, 0, 0, 0, 1, 0 \}$  &  \cr
\hline\hline
    \end{tabular}
  \vspace{2mm}
  \caption{\hspace{3mm}(cont.) The linkage diagrams $\gamma^{\vee}_{ij}(8)$ obtained as solutions of inequality  $\mathscr{B}^{\vee}_L(\gamma^{\vee}_{ij}(8)) < 2$}
  \label{sol_inequal_2}
  \end{table}

  \begin{table}[H]
  \centering
  \renewcommand{\arraystretch}{1.7}
  \small
  \begin{tabular} {||c|c|c|l|l||}
  \hline \hline
      The Carter &  $L_{ij}$  & $\mathscr{B}^{\vee}_L(\gamma^{\vee}_{ij}(8)) = p < 2$
                                & The linkage diagram  & $p$ \cr
      diagram &    &            &  \qquad \qquad $\gamma^{\vee}_{ij}(8)$  &  \\
     \hline\hline  
                  &  &   $q = 2a_4^2 + 4b_2^2 + 4a_4{b_2}$ & & \cr
                  \cline{3-3}
                  & $L_{12}$
                  & $\frac{1}{2}(3 + 2(a_4 + 2b_2) +  q) < 2$, or
                  & $\{ 0, 0, 1, 0, 0, 0 \}$ & \cr
                  &
                  & $(a_4 + 2b_2 + 1)^2 + a_4^2 < 2$
                  & $\{ 0, 0, 1, 0, 0, -1 \}$ &  \cr
                  & & & $\{ 0, 0, 1, 1, 0, -1 \}$ & \cr
                  & & & $\{ 0, 0, 1, -1, 0, 0 \}$ & \cr
                  \cline{2-4}
       $D_6(a_2)$ & $L_{13}$
                  & $\frac{1}{2}(3 + 2(-a_4 - 2b_2) +  q) < 2$, or
                  & $\{ 0, 1, 0, 0, 0, 0 \}$ & $\frac{3}{2}$ \cr
\hline\hline
    \end{tabular}
  \vspace{2mm}
  \caption{\hspace{3mm}(cont.) The linkage diagrams $\gamma^{\vee}_{ij}(8)$ obtained as solutions of inequality $\mathscr{B}^{\vee}_L(\gamma^{\vee}_{ij}(8)) < 2$}
  \label{sol_inequal_3}
  \end{table}

  \begin{table}[H]
  \centering
  \renewcommand{\arraystretch}{1.7}
  \small
  \begin{tabular} {||c|c|c|l|l||}
  \hline \hline
      The Carter &  $L_{ij}$  & $\mathscr{B}^{\vee}_L(\gamma^{\vee}_{ij}(8)) = p < 2$
                                & The linkage diagram  & $p$ \cr
      diagram &    &            &  \qquad \qquad $\gamma^{\vee}_{ij}(8)$  &  \\
     \hline\hline  
                  &  &   $q = 4a_4^2 + 8b_2^2 + 4b_3^2 + 8a_4{b_2}$ & & \cr
                  \cline{3-3}
                  & $L_{12}$
                  & $\frac{1}{4}(7 + 2(2a_4 + 4b_2 + 2b_3) +  q) < 2$, or
                  & $\{ 0, 0, 1, 1, 0, -1, 0 \}$ & \cr
                  &
                  & $(a_4 + b_3 + \frac{1}{2})^2 + (b_2 + \frac{1}{2})^2 + (b_3 + \frac{1}{2})^2 < 1$
                  &     $\{ 0, 0, 1, 1, 0, -1, -1 \}$ &  \cr
                  & & & $\{ 0, 0, 1, 0, 0, -1, 0 \}$ & \cr
                  & & & $\{ 0, 0, 1, 0, 0, 0, -1 \}$ & \cr
                  & & & $\{ 0, 0, 1, 0, 0,  0, 0 \}$ & \cr
\hline\hline
    \end{tabular}
  \vspace{2mm}
  \caption{\hspace{3mm}(cont.) The linkage diagrams $\gamma^{\vee}_{ij}(8)$ obtained as solutions of inequality $\mathscr{B}^{\vee}_L(\gamma^{\vee}_{ij}(8)) < 2$}
  \label{sol_inequal_4}
  \end{table}


  \begin{table}[H]
  \centering
  \renewcommand{\arraystretch}{1.4}
  \small
  \begin{tabular} {||c|c|c|l|l||}
  \hline \hline
     The Carter &  $L_{ij}$  & $\mathscr{B}^{\vee}_L(\gamma^{\vee}_{ij}(8)) = p < 2$
                                & The linkage diagram  & $p$ \cr
      diagram &    &            &  \qquad \qquad $\gamma^{\vee}_{ij}(8)$  &  \\
    \hline  
                  &  &   $q = 4b_2^2 + 4b_2{b_3} + 4b_3^2$ & & \cr
                  \cline{3-3}
                  &  $L_{12}$
                  &  $\frac{1}{3}(6 + 2(3b_2 + 3b_3) + q) < 2$, or
                  &  $\{ 0, 0, 1, 0, 0, -1 \}$  &  \cr
                  &
                  &  $(b_2 + b_3)^2 + (b_2 + \frac{3}{2})^2  + (b_3 + \frac{3}{2})^2 < \frac{9}{2}$
                  &  $\{ 0, 0, 1, 0, -1, 0 \}$ & \cr
                  \cline{2-4}
    $E_6$         &  $L_{13}$
                  &  $\frac{1}{3}(10 + 2(4b_2 + 5b_3) + q) < 2$, or
                  &  $\{ 0, 1, 0, 0, 0, -1 \}$  & $\frac{4}{3}$ \cr
                  &
                  &  $(b_2 + b_3)^2 + (b_2 + 2)^2  + (b_3 + \frac{5}{2})^2 < \frac{33}{4}$
                  &  $\{ 0, 1, 0, 0,  -1, -1 \}$ &   \cr
                  \cline{2-4}
                  &  $L_{23}$
                  &  $\frac{1}{3}(10 + 2(5b_2 + 4b_3) + q) < 2$,
                  &  $\{ 0, 1, 0, 0, -1, 0 \}$  &  \cr
\end{tabular}
  \vspace{2mm}
  \caption{\hspace{3mm}(cont.) The linkage diagrams $\gamma^{\vee}_{ij}(8)$ obtained as solutions of inequality $\mathscr{B}^{\vee}_L(\gamma^{\vee}_{ij}(8)) < 2$}
  \label{sol_inequal_5}
  \end{table}

  \begin{table}[H]
  \centering
  \renewcommand{\arraystretch}{1.7}
  \small
  \begin{tabular} {||c|c|c|l|l||}
  \hline \hline
      The Carter &  $L_{ij}$  & $\mathscr{B}^{\vee}_L(\gamma^{\vee}_{ij}(8)) = p < 2$
                                & The linkage diagram  & $p$ \cr
      diagram &    &            &  \qquad \qquad $\gamma^{\vee}_{ij}(8)$  &  \\
     \hline\hline  
                  &  &   $q = 8a_4^2 + 12b_2^2 + 4b_3^2 + 16b_2{a}_4 + 8b_2{b}_3 + 8b_3{a}_4$ & & \cr
                  \cline{3-3}
                  & $L_{12}$
                  & $\frac{1}{4}(7 + 4(2a_4 + 3b_2 + b_3) +  q) < 2$, or
                  & $\{ 0, 0, 1, 0, 0, 0, 0 \}$ & \cr
                  &
                  & $(a_4 + b_2 + b_3 + \frac{1}{2})^2 + (a_4 + b_2 + \frac{1}{2})^2 + (b_2 + \frac{1}{2})^2 < 1$
                  &     $\{ 0, 0, 1, 0, 0, -1, 1 \}$ &  \cr
                  & & & $\{ 0, 0, 1, 1, 0, -1, -1 \}$ & \cr
                  & & & $\{ 0, 0, 1, -1, 0, 0, 0 \}$ & \cr
                  & & & $\{ 0, 0, 1, 1, 0,  -1, 0 \}$ & $\frac{7}{4}$ \cr
\hline\hline
    \end{tabular}
  \vspace{2mm}
  \caption{\hspace{3mm}(cont.) The linkage diagrams $\gamma^{\vee}_{ij}(8)$ obtained as solutions of inequality $\mathscr{B}^{\vee}_L(\gamma^{\vee}_{ij}(8)) < 2$}
  \label{sol_inequal_6}
  \end{table}

\newpage
\subsection{$\beta$-unicolored linkage diagrams. Solutions of inequality $\mathscr{B}^{\vee}_L(\gamma^{\vee}) < 2$}
~\\

\begin{table}[H]
  \centering
  \renewcommand{\arraystretch}{1.6}
  \small
  \begin{tabular} {||c|l|l|l||}
  \hline \hline
      Diagram & \multicolumn{2}{c|}{$\mathscr{B}^{\vee}_L(\gamma^{\vee}) = p < 2$,}  & $p$ \cr
              & \multicolumn{2}{c|}{$\gamma^{\vee}$ -  $\beta$-unicolored linkage diagram} &  \\
    \hline\hline  
                  &  \multicolumn{2}{c|}{$\frac{1}{3}(4b_2^2 + 4b_3^2 + 4b_2{b_3}) < 2$, or} &  \cr
    $E_6(a_1)$, $E_6(a_2)$
                  &  \multicolumn{2}{c|}{$(b_2 + b_3)^2 + b_2^2  + b_3^2 < 3$} & \cr
    \cline{2-3}
                  &  $\{ 0, 0, 0, 0, 0, 1 \}$ \quad\quad\quad\quad\quad\quad & $\{ 0, 0, 0, 0, 0, -1 \}$ & $\frac{4}{3}$ \cr
                  &  $\{ 0, 0, 0, 0, 1, 0 \}$  & $\{ 0, 0, 0, 0, -1, 0 \}$ &  \cr
                  &  $\{ 0, 0, 0, 0, 1, -1 \}$  & $\{ 0, 0, 0, 0, -1, 1 \}$ & \cr
      \hline
                  &  \multicolumn{2}{c|}
      {$\frac{1}{2}(3b_2^2 + 4b_3^2 + 3b_4^2 + 4b_2{b_3} + 2b_2{b_4} + 4b_3{b_4}) < 2$, or} &  \cr
\hline\hline
    \end{tabular}
  \vspace{2mm}
  \caption{\hspace{3mm}(cont.) $\beta$-unicolored linkage diagrams obtained as solutions of inequality
    $\mathscr{B}^{\vee}_L(\gamma^{\vee}) < 2$}
  \label{homog_inequal_1}
  \end{table}

\begin{table}[H]
  \centering
  \renewcommand{\arraystretch}{1.7}
  \small
  \begin{tabular} {||c|l|l|l||}
  \hline \hline
      Diagram & \multicolumn{2}{c|}{$\mathscr{B}^{\vee}_L(\gamma^{\vee}) = p < 2$,}  & $p$ \cr
              & \multicolumn{2}{c|}{$\gamma^{\vee}$ -  $\beta$-unicolored linkage diagram} &  \\
      \hline\hline
       $D_5(a_1)$ &  \multicolumn{2}{c|}
      {$\frac{1}{4}(4b_2^2) < 2 ,\quad $ or $\quad b_2^2  < 2$} &  \cr
    \cline{2-3}
      & $\{ 0, 0, 0, 0, 1\}$  \quad\quad\quad\quad & $\{ 0, 0, 0, 0, -1\}$ & $1$ \cr
      \hline
       $D_6(a_1)$ &  \multicolumn{2}{c|}
      {$\frac{1}{2}(2b_2^2 + 2b_3^2) < 2 ,\quad $ or $\quad b_2^2 + b_3^2  < 2$} &  \cr
    \cline{2-3}
      & $\{ 0, 0, 0, 0, 1, 0 \}$   & $\{ 0, 0, 0, 0, -1, 0 \}$ & $1$ \cr
      & $\{ 0, 0, 0, 0, 0, 1 \}$  & $\{ 0, 0, 0, 0, 0, -1 \}$ & \cr
      \hline
        &  \multicolumn{2}{c|}
      {$\frac{1}{2}(2a_4^2 + 4b_1^2 + 4b_2^2 + 4a_4{b_2}) < 2,$ or} &  \cr
\hline\hline
         \end{tabular}
  \vspace{2mm}
  \caption{\hspace{3mm}(cont.) $\beta$-unicolored linkage diagrams obtained as solutions of inequality
    $\mathscr{B}^{\vee}_L(\gamma^{\vee}) < 2$}
  \label{homog_inequal_2}
  \end{table}

  \begin{table}[H]
  \centering
  \renewcommand{\arraystretch}{1.6}
  \small
  \begin{tabular} {||c|l|l|l||}
  \hline \hline
      Diagram & \multicolumn{2}{c|}{$\mathscr{B}^{\vee}_L(\gamma^{\vee}) = p < 2$,}  & $p$ \cr
              & \multicolumn{2}{c|}{$\gamma^{\vee}$ -  $\beta$-unicolored linkage diagram} &  \\
    \hline\hline  
                  &  \multicolumn{2}{c|}{$\frac{1}{3}(4b_2^2 + 4b_3^2 + 4b_2{b_3}) < 2$, or} &  \cr
    $E_6$
                  &  \multicolumn{2}{c|}{$(b_2 + b_3)^2 + b_2^2  + b_3^2 < 3$} & \cr
    \cline{2-3}
    & \multicolumn{2}{c|}{(coincide with $E_6(a_1)$, $E_6(a_2)$)} & \cr
                  &  $\{ 0, 0, 0, 0, 0, 1 \}$ \quad\quad\quad\quad\quad\quad & $\{ 0, 0, 0, 0, 0, -1 \}$ & $\frac{4}{3}$ \cr
                  &  $\{ 0, 0, 0, 0, 1, 0 \}$  & $\{ 0, 0, 0, 0, -1, 0 \}$ &  \cr
                  &  $\{ 0, 0, 0, 0, 1, -1 \}$  & $\{ 0, 0, 0, 0, -1, 1 \}$ & \cr
\hline\hline
    \end{tabular}
  \vspace{2mm}
  \caption{\hspace{3mm}(cont.) $\beta$-unicolored linkage diagrams obtained as solutions of inequality
    $\mathscr{B}^{\vee}_L(\gamma^{\vee}) < 2$}
  \label{homog_inequal_3}
  \end{table}

\newpage
\subsection{Linkage diagrams $\gamma^{\vee}_{ij}(6)$ per loctets and components}
  \label{sec_diagr_per_comp}
~\\

\begin{table}[H]
  \small
  \centering
  \renewcommand{\arraystretch}{1}
  \begin{tabular} {|| c | c || l | l | l ||}
  \hline \hline
    Diagram & Comp. & \multicolumn{3}{c||}{Linkage diagrams $\gamma^{\vee}_{ij}(6)$ for the loctet of type $L_{ij}$}    \\
            \cline{3-5}
            &           & \multicolumn{1}{c|}{Type $L_{12}: (1, 1, 0, -1, ...)$}
                        & \multicolumn{1}{c|}{Type $L_{13}: (1, 0, 1, -1, ...)$}
                        & \multicolumn{1}{c||}{Type $L_{23}: (0, 1, 1, -1, ...)$} \\
    \hline \hline
    $D_5(a_1)$ & $1$  & $(1, 1, 0, -1, 0) \quad L^a_{12}$
                      & $(1, 0, 1, -1, 0) \quad  L^a_{13}$  & -  \cr
               & $2$  & $(1, 1, 0, -1, 1) \quad  L^b_{12}$
                      & $(1, 0, 1, -1, -1) \quad  L^b_{13} $  & -  \cr
               & $3$  & - & - &  $(0, 1, 1, -1, 0) \quad  L^a_{23}$  \\
   \hline
    $D_6(a_1)$ & $1$  & $(1, 1, 0, -1, 1, -1) \quad L^d_{12}$
                      & $(1, 0, 1, -1, -1, -1) \quad L^d_{13}$  & -  \cr
               &      & $(1, 1, 0, -1, 0, 0) \quad L^a_{12}$
                      & $(1, 0, 1, -1, 0, 0) \quad  L^a_{13}$    & -  \cr
               & $2$  & $(1, 1, 0, -1, 0, -1) \quad L^b_{12}$
                      & $(1, 0, 1, -1, 0, -1) \quad L^b_{13}$   & -  \cr
\end{tabular}
  \vspace{2mm}
  \caption{\hspace{3mm}Linkage diagrams $\gamma^{\vee}_{ij}(6)$ for the Carter diagrams
  from the class $\mathsf{C4}$, $n < 8$.
   For $D_6(a_2)$, $D_7(a_1)$, $D_7(a_2)$, the length of the $\alpha$-set is $4$}
  \label{tab_seed_linkages_6}
  \end{table}

%% file: A6linkages.tex
\section{\sc\bf The linkage systems for the Carter diagrams}
 \label{sect_linkage_diagr}
~\\
  The linkage systems are similar to the weight systems (= weight diagrams)
  of the irreducible representations of the semisimple Lie algebras.
~\\
~\\
\subsection{The linkage systems $D_4(a_1)$, $D_5(a_1)$, $D_6(a_1)$, $D_6(a_2)$}
~\\

\begin{figure}[H]
\centering
\includegraphics[scale=0.8]{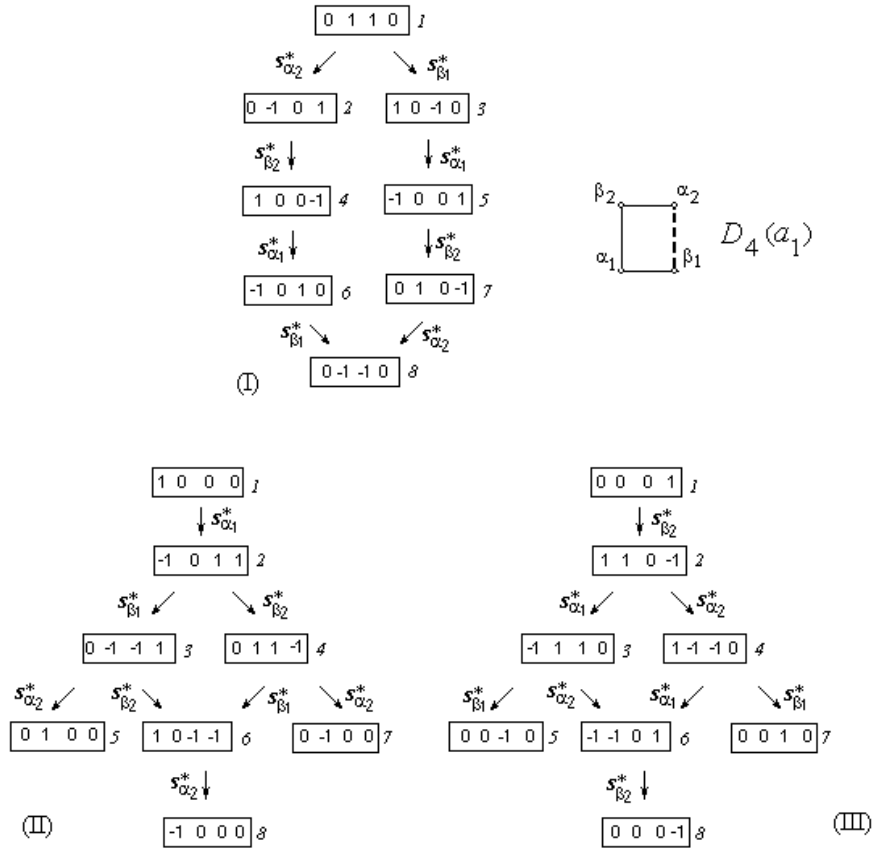}
\caption[\hspace{3mm}Three components of the linkage system $D_4(a_1)$, 3 components]{\hspace{3mm}Three components of the linkage system $D_4(a_1)$. There are $24$ linkage diagrams in the case $D_4(a_1)$}
\label{D4a1_linkages}
\end{figure}

\begin{figure}[H]
\centering
\includegraphics[scale=1.5]{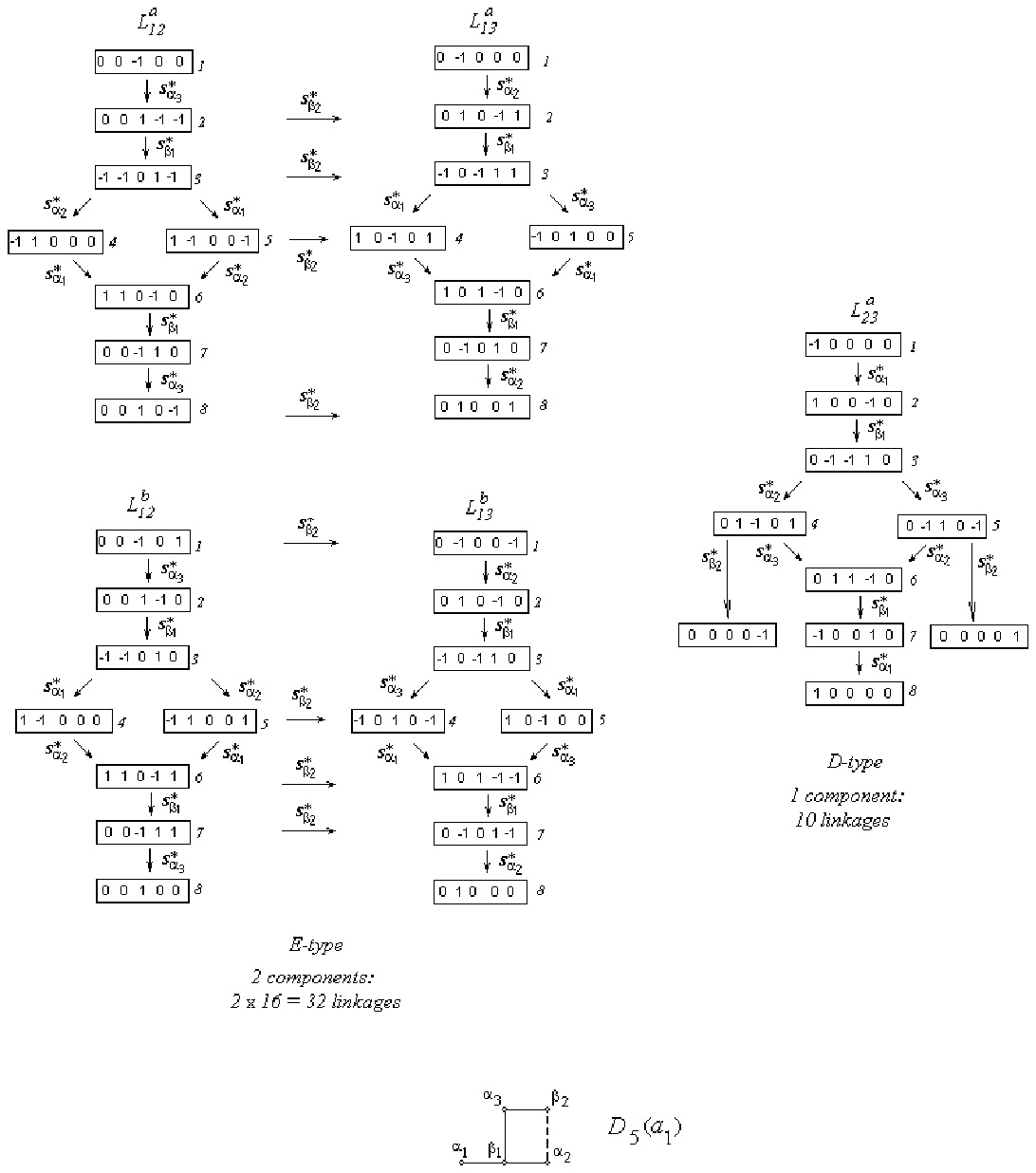}
\vspace{3mm}
\caption[\hspace{3mm}The linkage system $D_5(a_1)$, $3$ components, $5$ loctets]{\hspace{3mm}The linkage system  $D_5(a_1)$. There are one component
of the $D$-type containing $10$ linkage diagrams, and two components of the $E$-type containing
$2\times16 = 32$ elements}
\label{D5a1_linkages}
\end{figure}

\begin{figure}[H]
\centering
\includegraphics[scale=1.3]{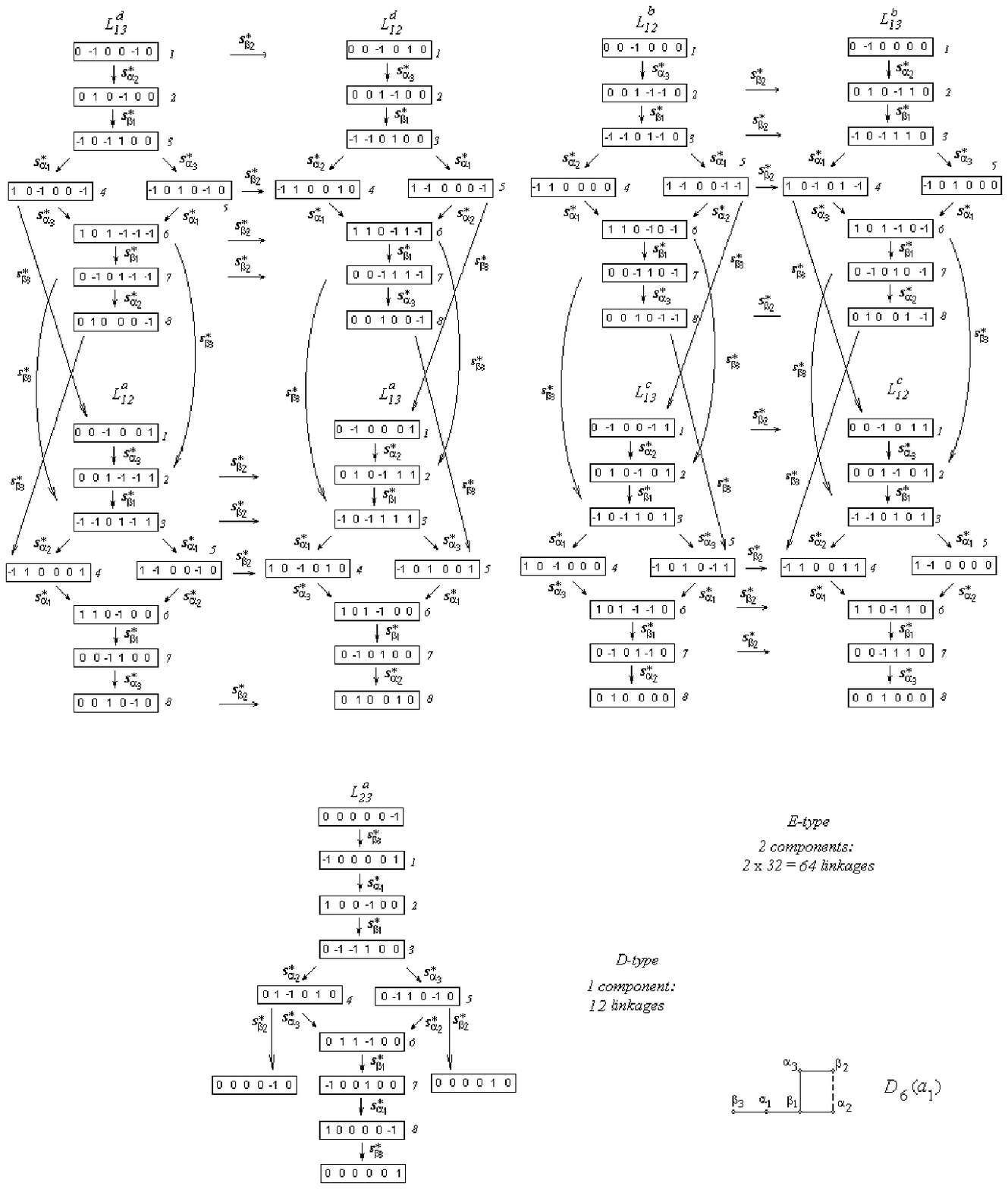}
\vspace{3mm}
\caption[\hspace{3mm}The linkage system $D_6(a_1)$, $3$ components, $9$ loctets]{\hspace{3mm}The linkage system  $D_6(a_1)$.
There are $12$ linkage diagrams, $1$ loctet in the single $D$-type component,
and $2\times32 = 64$ linkage diagrams, $8$ loctets in two $E$-type components}
\label{D6a1_linkages}
\end{figure}

\begin{figure}[H]
\centering
\includegraphics[scale=1.3]{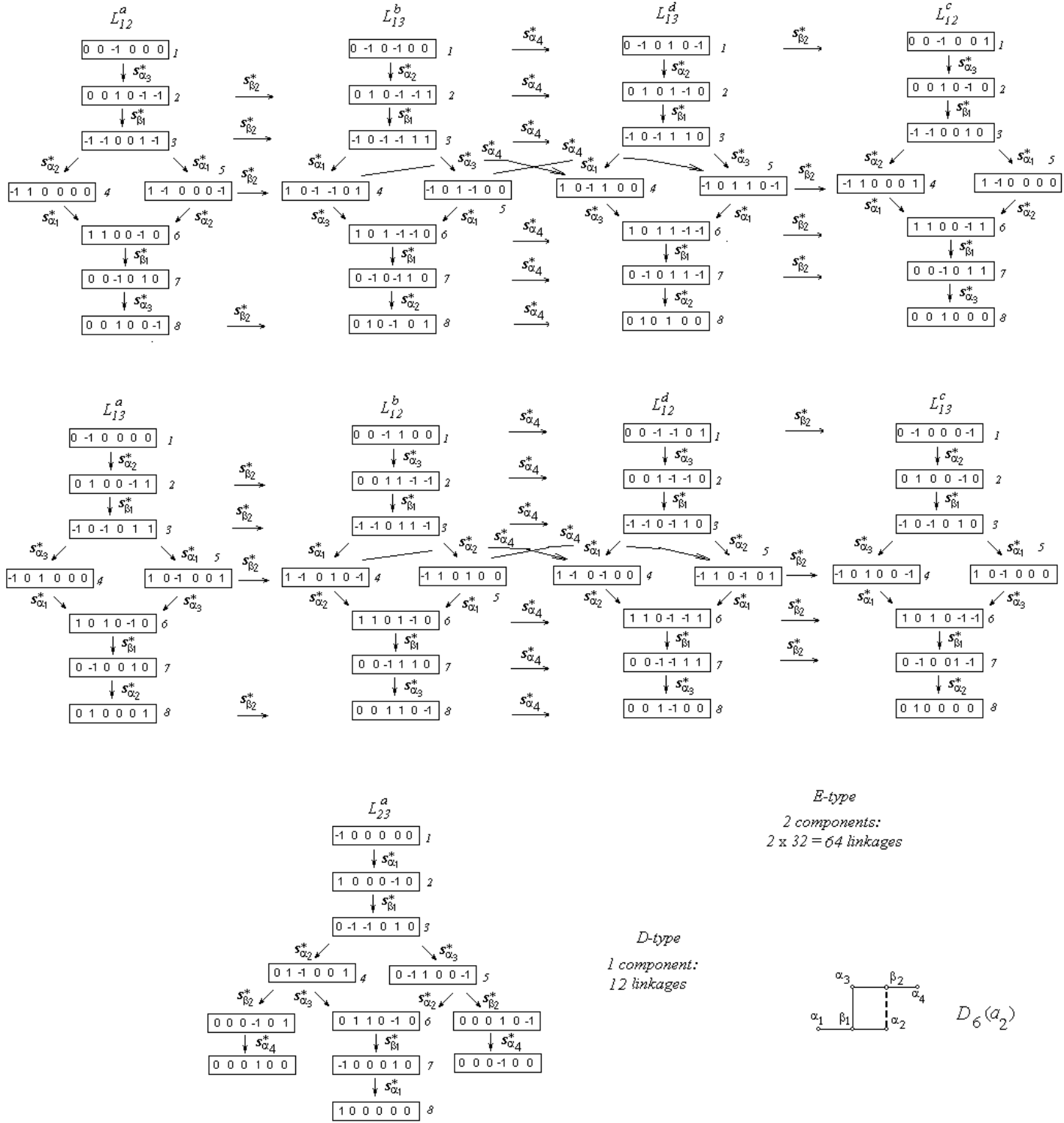}
\vspace{3mm}
\caption[\hspace{3mm}The linkage system $D_6(a_2)$, $3$ components, $9$ loctets]{\hspace{3mm}The linkage system  $D_6(a_2)$.
There are $12$ linkages, $1$ loctet in the single $D$-type component,
and $2\times32 =64$ linkages, $8$ loctets in two $E$-type components}
\label{D6a2_linkages}
\end{figure}

\newpage
\subsection{The linkage systems $E_6(a_1)$, $E_6(a_2)$, $E_7(a_1)$, $E_7(a_2)$, $E_7(a_3)$, $E_7(a_4)$}
~\\

\begin{figure}[H]
\centering
\includegraphics[scale=1.4]{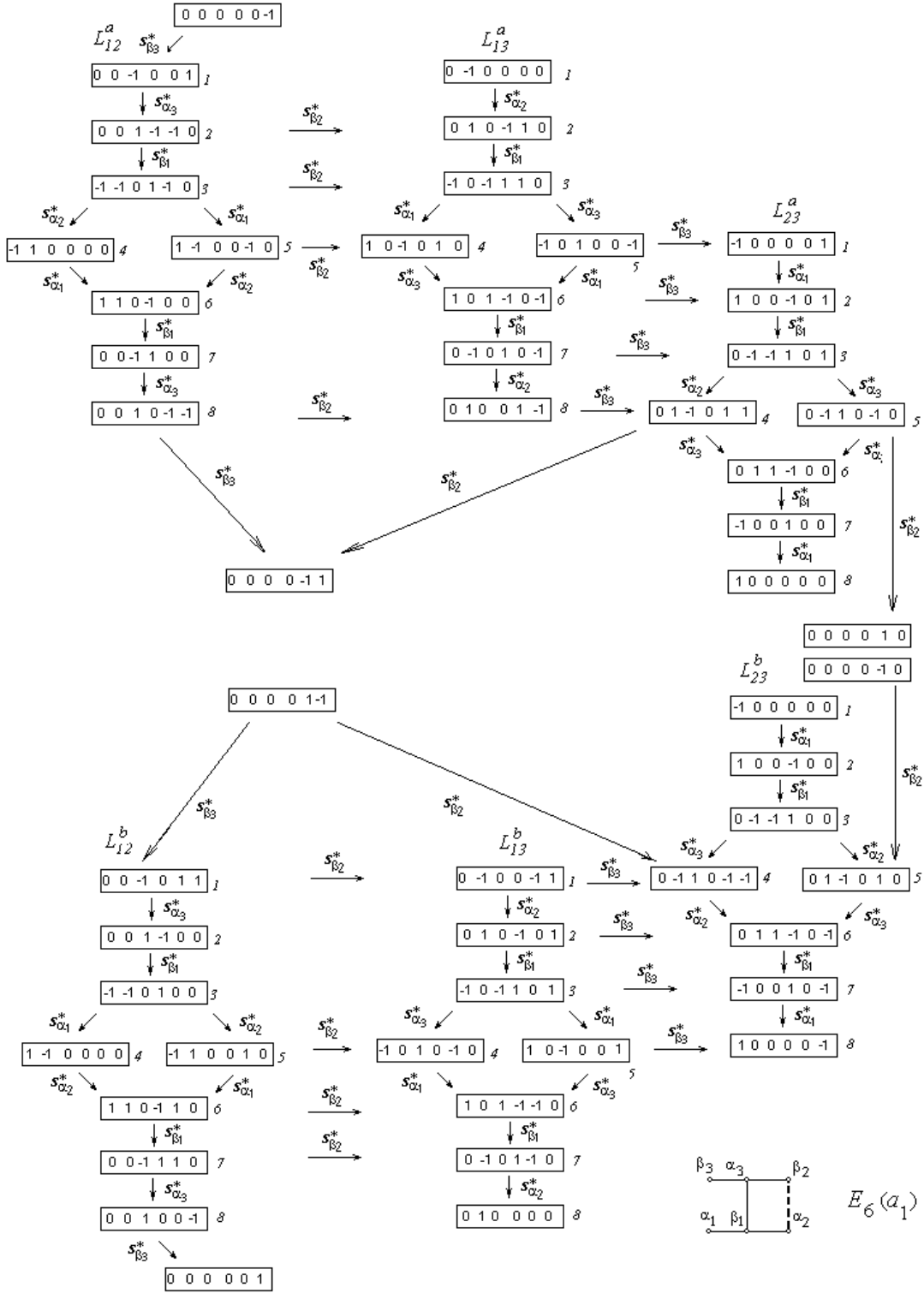}
\vspace{3mm}
\caption{\hspace{3mm}The linkage system $E_6(a_1)$, two components, $54$ linkage diagrams, $6$ loctets}
\label{E6a1_linkages}
\end{figure}

\begin{figure}[H]
\centering
\includegraphics[scale=1.5]{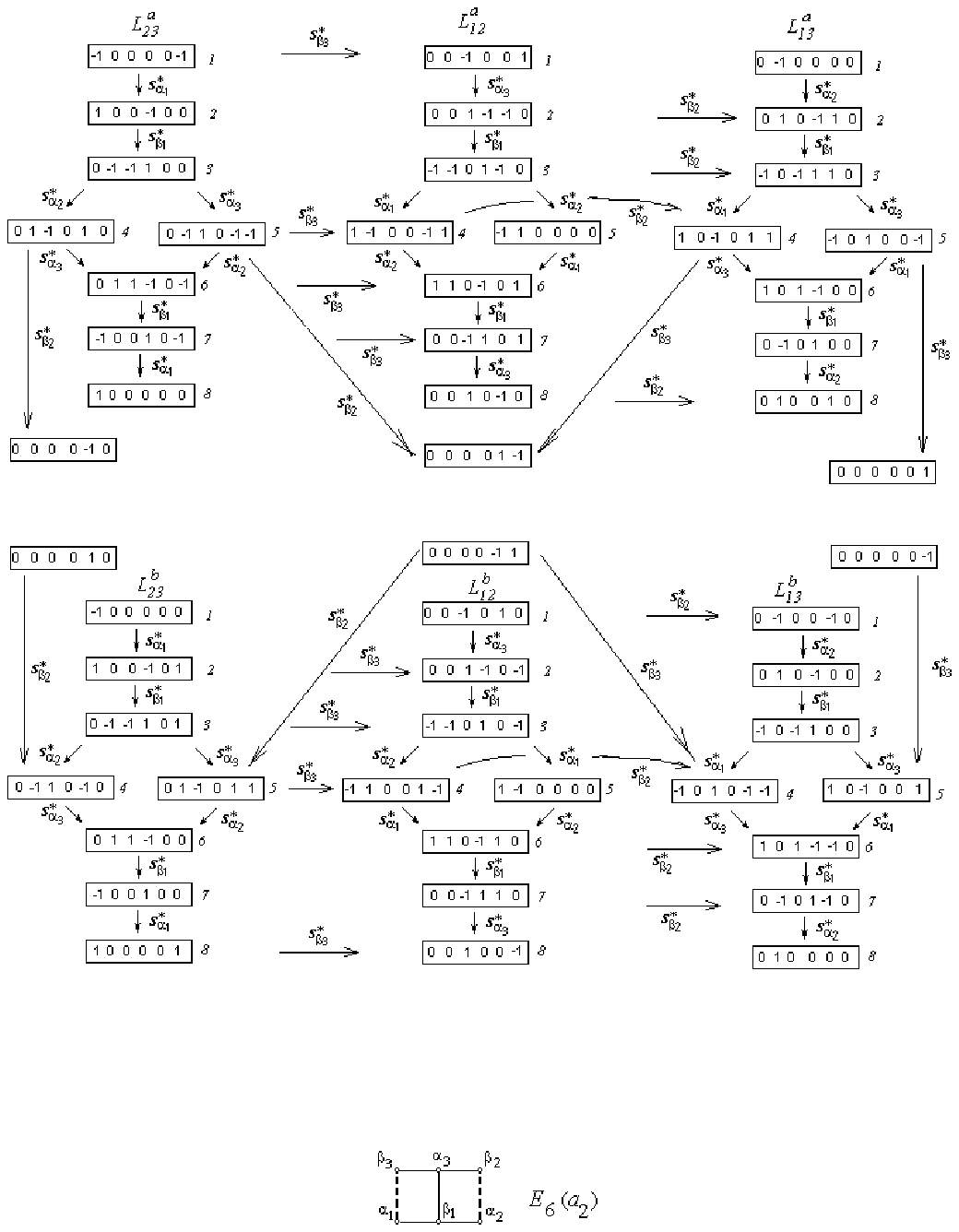}
\vspace{3mm}
\caption{\hspace{3mm}The linkage system $E_6(a_2)$, two components, $54$ linkage diagrams, $6$ loctets}
\label{E6a2_linkages}
\end{figure}

\begin{figure}[H]
\centering
\includegraphics[scale=1.4]{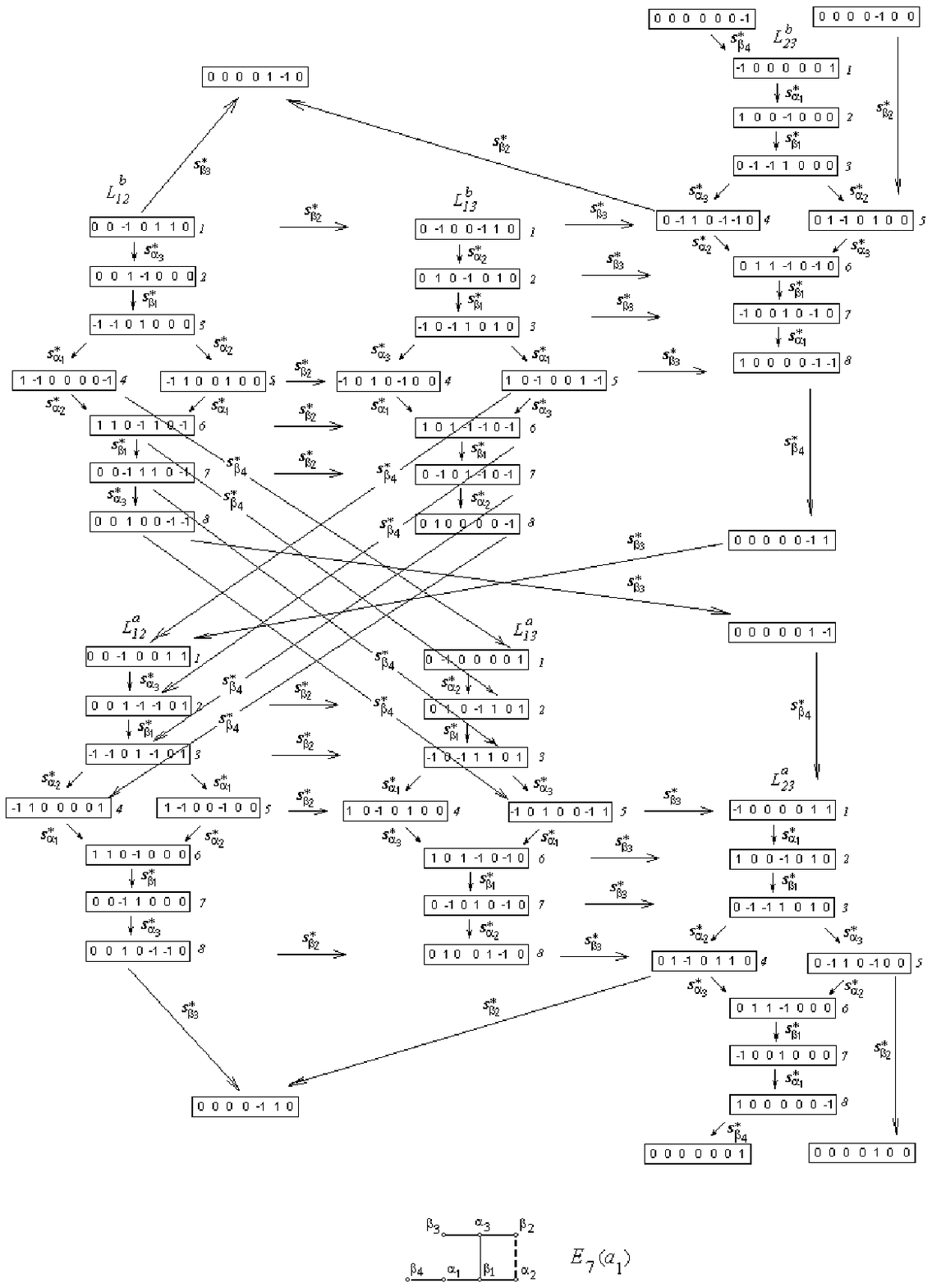}
\caption{\hspace{3mm}The linkage system $E_7(a_1)$, one component, $56$ linkage diagrams, $6$ loctets}
\label{E7a1_linkages}
\end{figure}

\begin{figure}[H]
\centering
\includegraphics[scale=1.4]{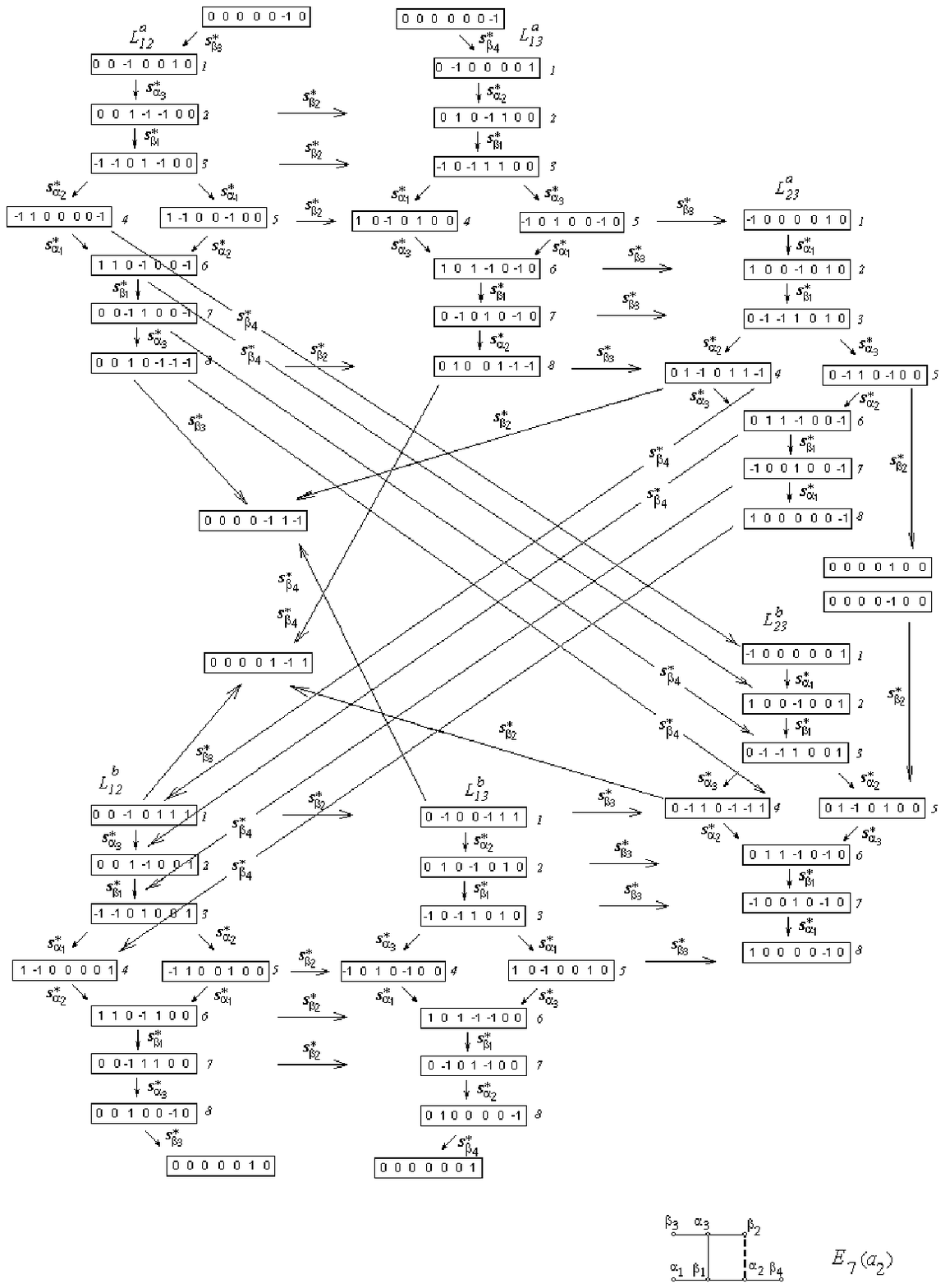}
\caption{\hspace{3mm}The linkage system $E_7(a_2)$, one component, $56$ linkage diagrams, $6$ loctets}
\label{E7a2_linkages}
\end{figure}

\begin{figure}[H]
\centering
\includegraphics[scale=1.4]{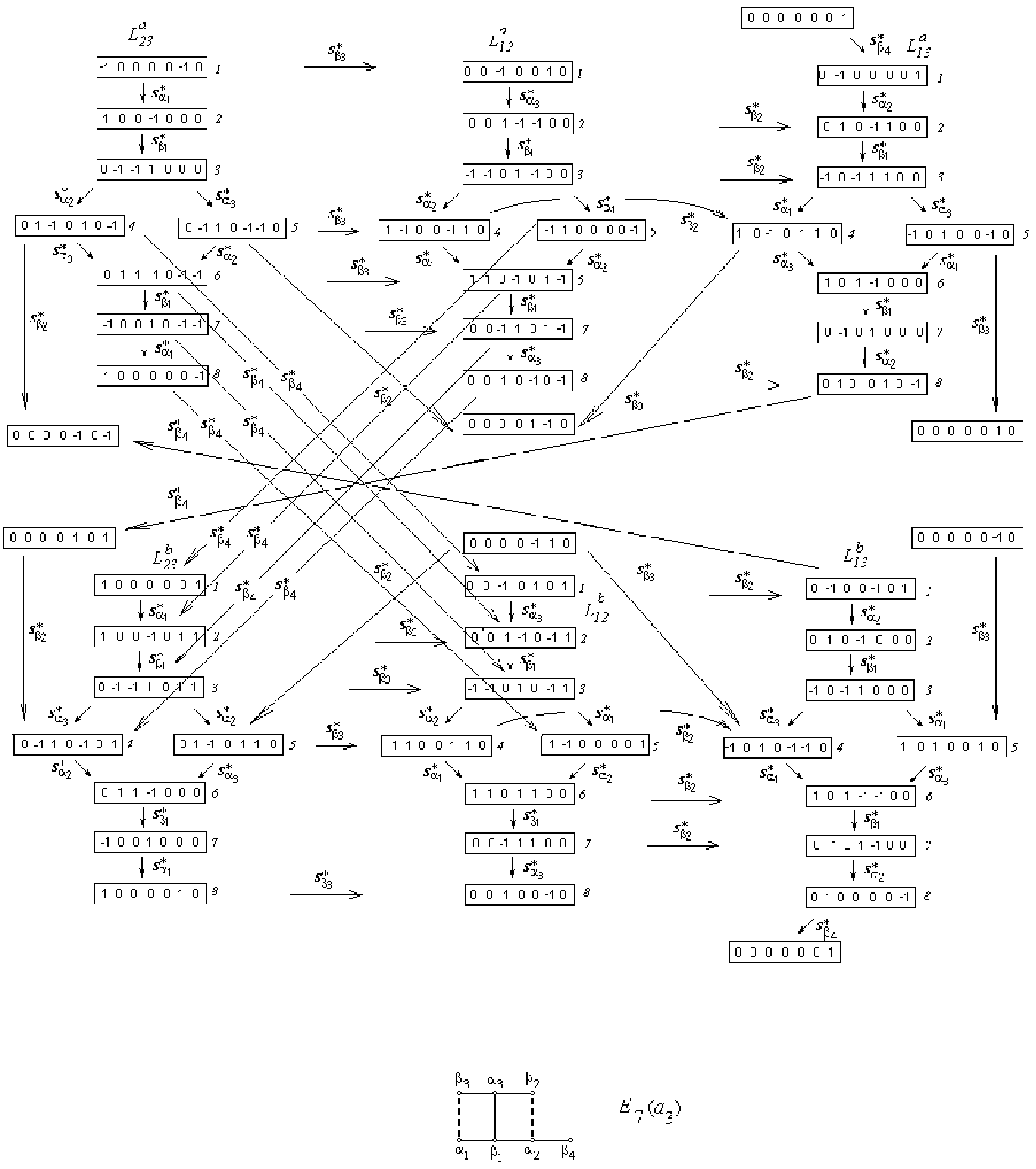}
\caption{\hspace{3mm}The linkage system $E_7(a_3)$, one component, $56$ linkage diagrams, $6$ loctets}
\label{E7a3_linkages}
\end{figure}

\begin{figure}[H]
\centering
\includegraphics[scale=1.4]{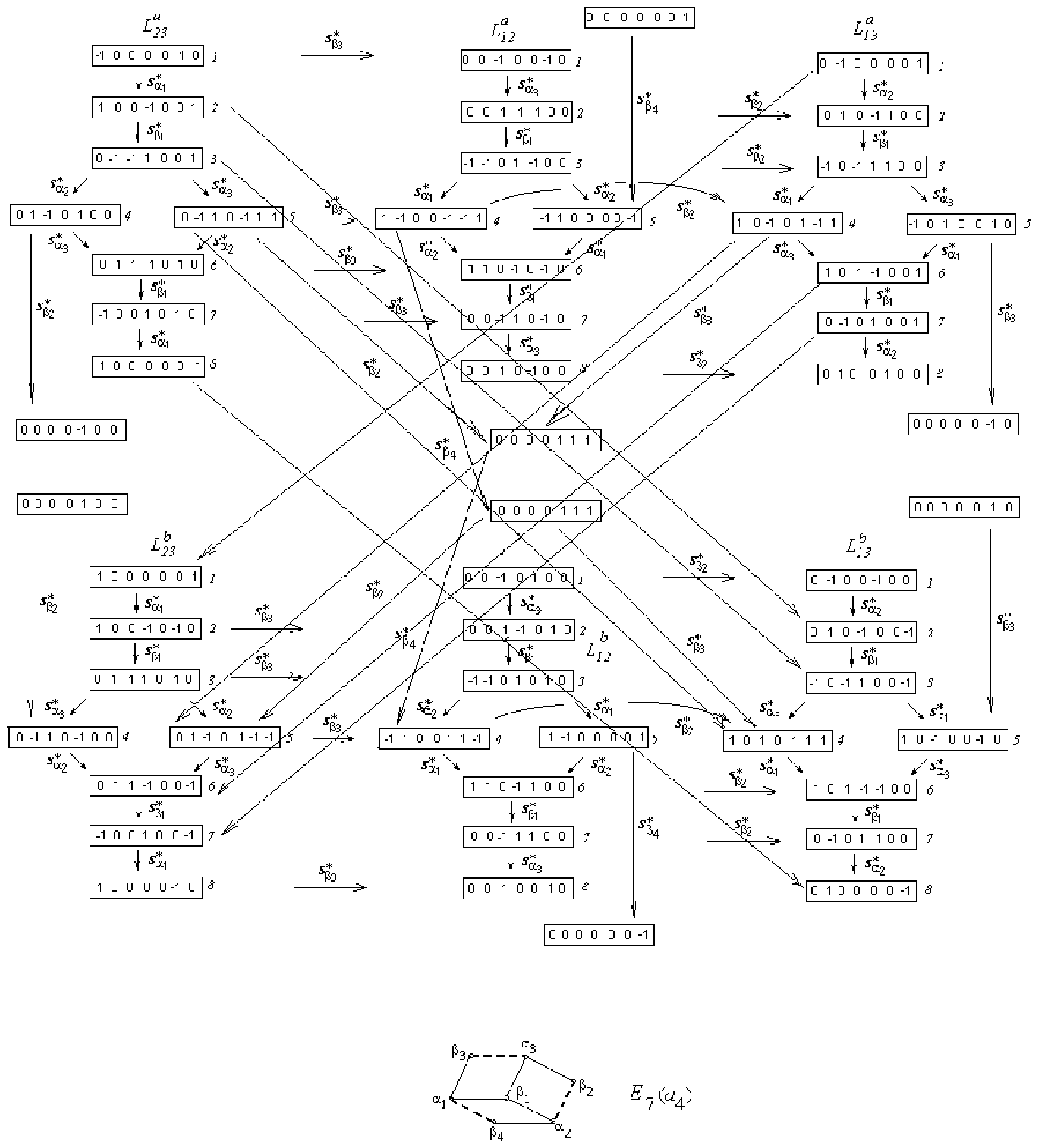}
\caption{\hspace{3mm}The linkage system $E_7(a_4)$,  one component, $56$ linkage diagrams,
 $6$ loctets}
\label{E7a4_linkages}
\end{figure}

\newpage
\subsection{The linkage systems and weight systems for $E_6$, $E_7$, $D_5$, $D_6$}
~\\

\begin{figure}[H]
\centering
\includegraphics[scale=0.9]{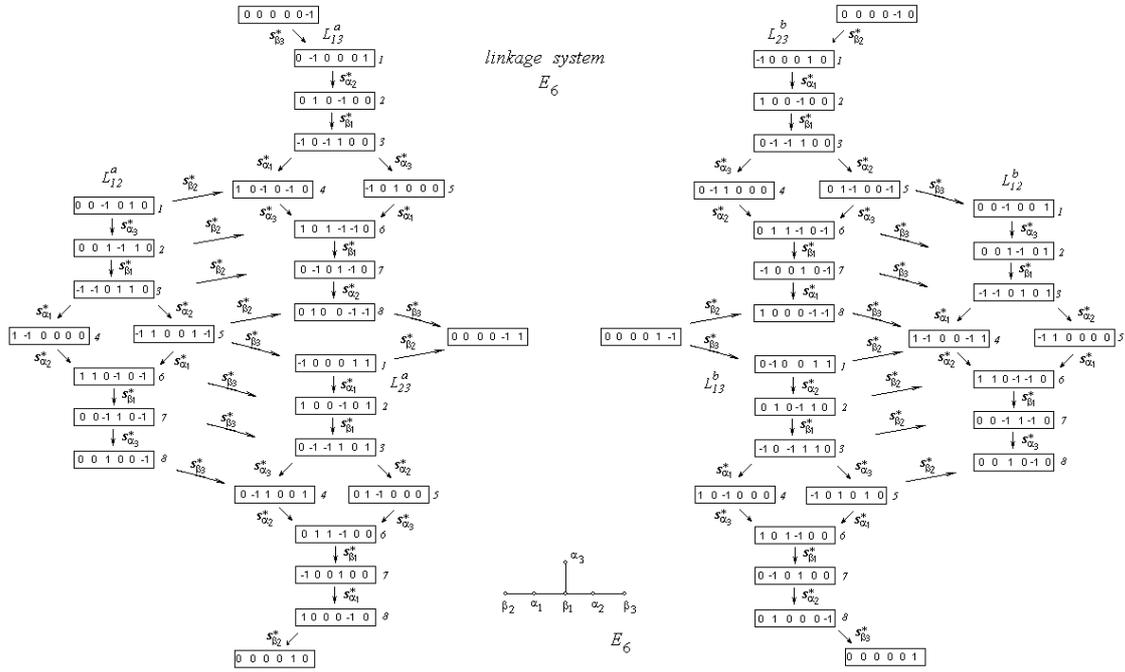}
\caption{\hspace{3mm}The linkage system $E_6$, two components, each of which contains $27$ elements  }
\label{E6pure_loctets_norm_gam8}
\end{figure}

\begin{figure}[H]
\centering
\includegraphics[scale=0.9]{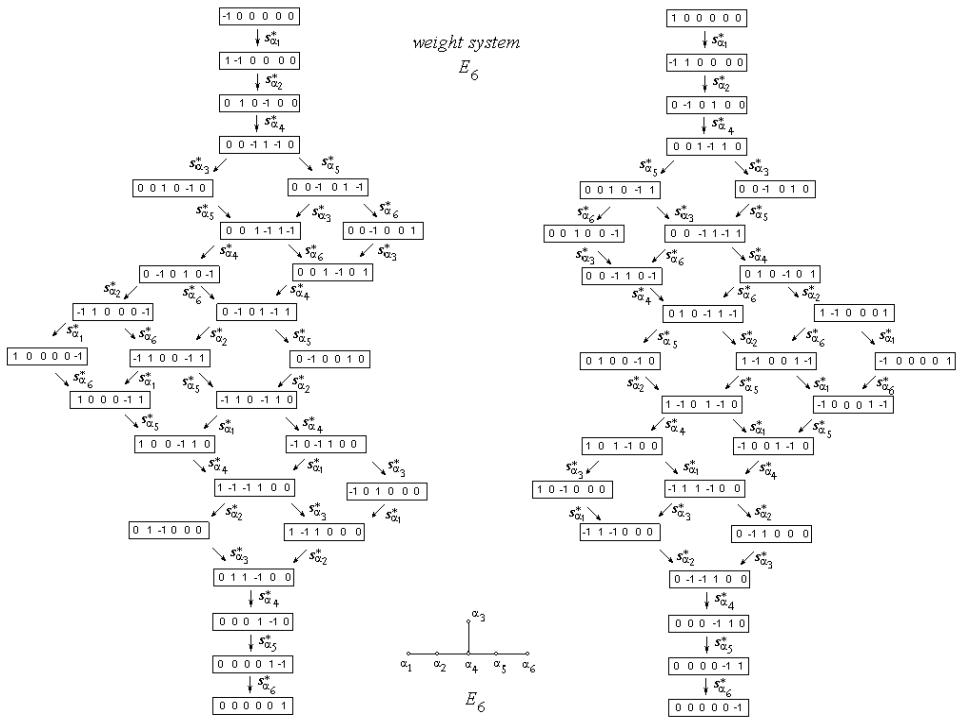}
\caption[\hspace{3mm}The weight systems of ${\bf 27}$ and $\overline{\bf 27}$  of the semisimple Lie algebra $E_6$]{\hspace{3mm}The weight systems (= weight diagrams) of the fundamental representations ${\bf 27}$ and $\overline{\bf 27}$  of the semisimple Lie algebra $E_6$}
\label{27_weight_diagr_E6__2comp}
\end{figure}


\begin{figure}[H]
\centering
\includegraphics[scale=1.2]{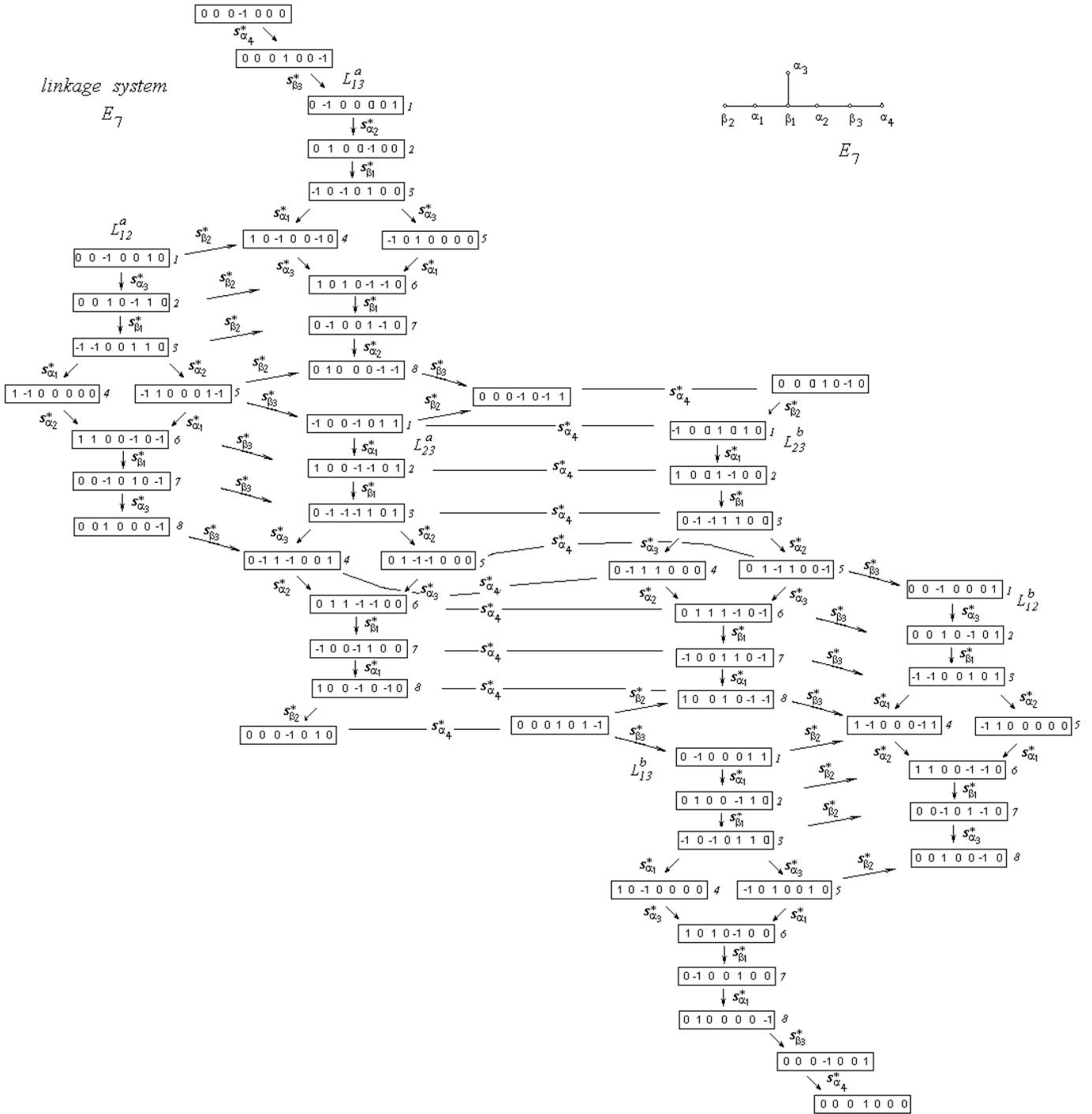}
\caption{\hspace{3mm}The linkage system $E_7$, one component, $56$ elements  }
\label{E7pure_linkage_system}
\end{figure}

\begin{figure}[H]
\centering
\includegraphics[scale=1.2]{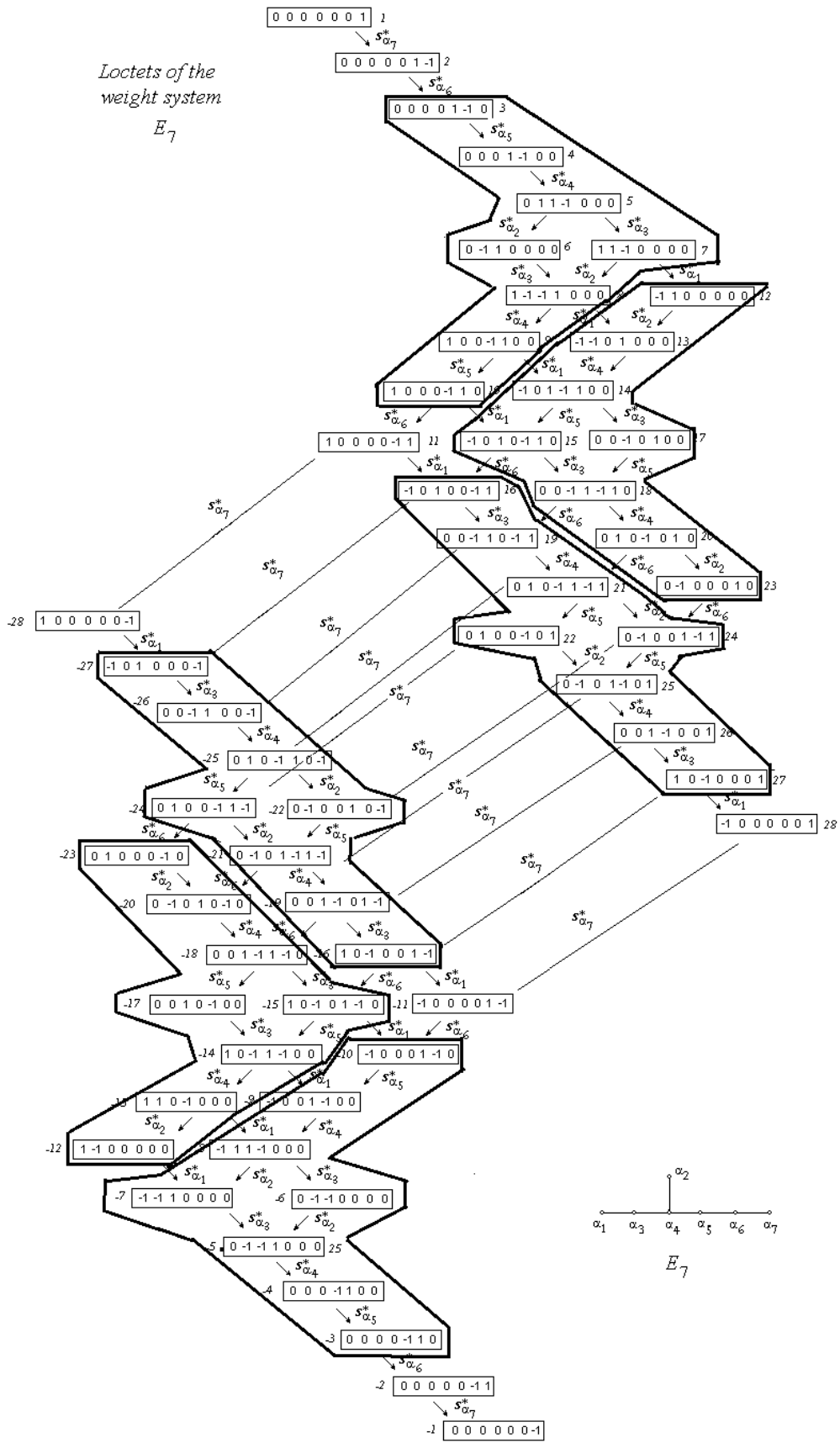}
\vspace{3mm}
\caption{\hspace{3mm}Loctets in the weight system of the fundamental representations ${\bf 56}$  of  $E_7$}
\label{56_weight_diagr_E7}
\end{figure}


\begin{figure}[H]
\centering
\includegraphics[scale=0.9]{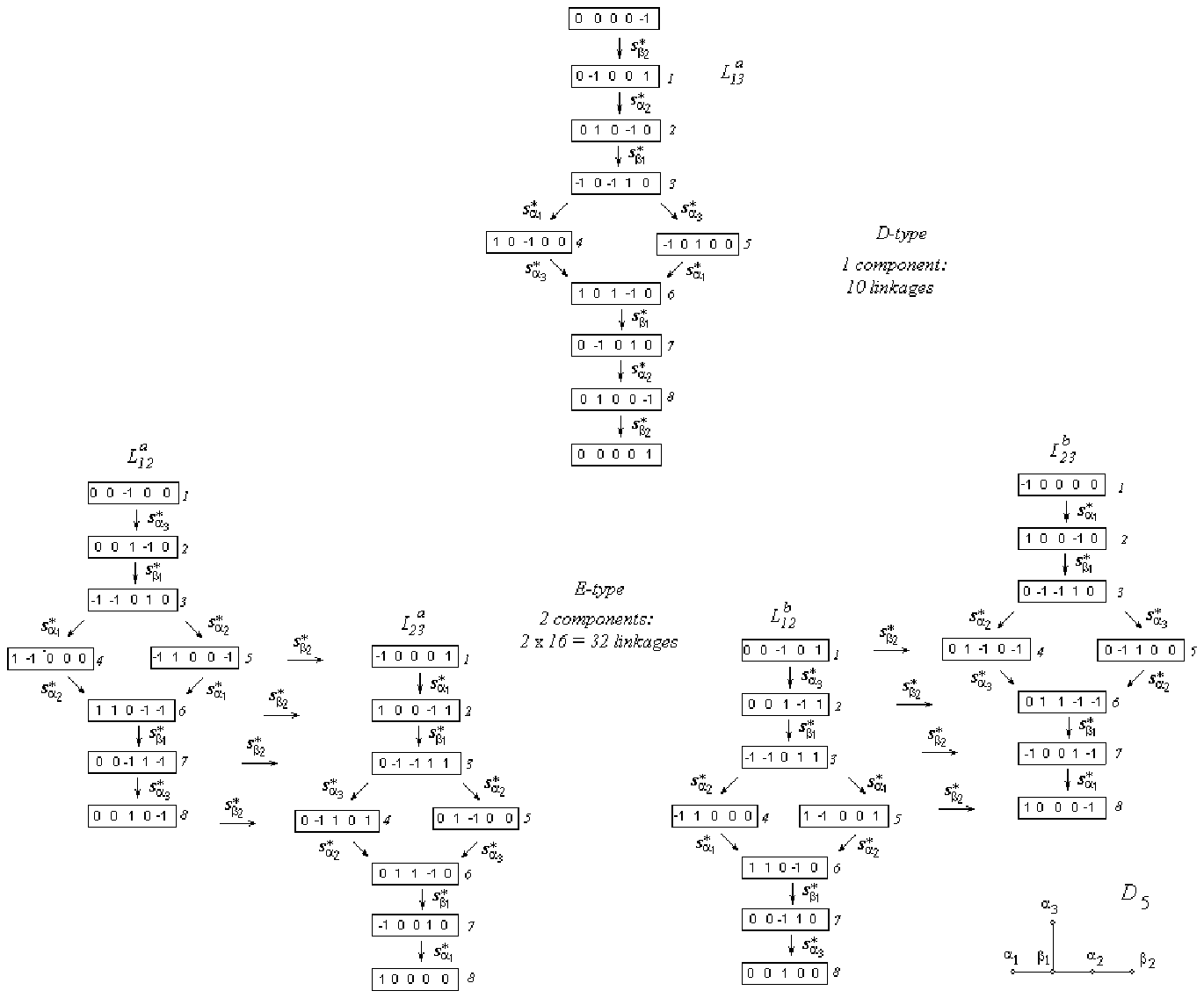}
\vspace{3mm}
\caption[\hspace{3mm}Linkage system for $D_5$, 3 components, 10 + 2x16 linkages]{\hspace{3mm}The linkage system  $D_5$. There are one component
of the $D$-type containing $10$ linkage diagrams and two components of the $E$-type containing
$2\times16 = 32$ elements}
\label{D5pure_loctets}
\end{figure}

\begin{figure}[H]
\centering
\includegraphics[scale=0.9]{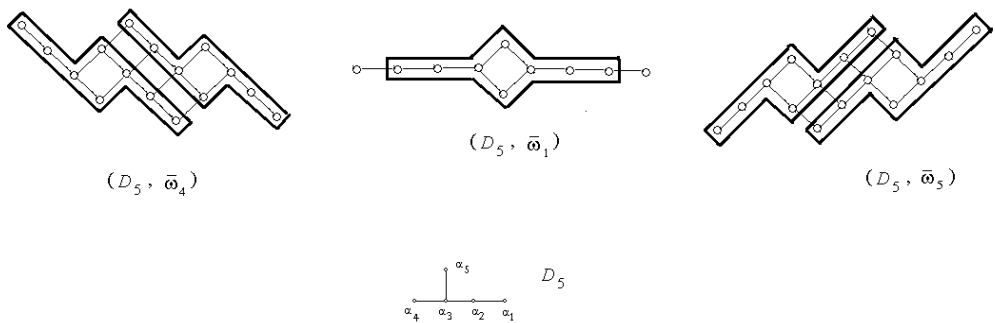}
\vspace{3mm}
\caption[\hspace{3mm}Loctets in the weight system for representations $(D_5, \overline{\omega}_1)$, 
 $(D_5, \overline{\omega}_4)$ and $(D_5, \overline{\omega}_5)$]
 {\hspace{3mm}Loctets in the weight system for $3$ fundamental representations of $D_5$:
  $(D_5, \overline{\omega}_1$), $(D_5, \overline{\omega}_4$) and $(D_5, \overline{\omega}_5$)  }
\label{16x2_10_weight_diagr_D5}
\end{figure}

\begin{figure}[H]
\centering
\includegraphics[scale=0.85]{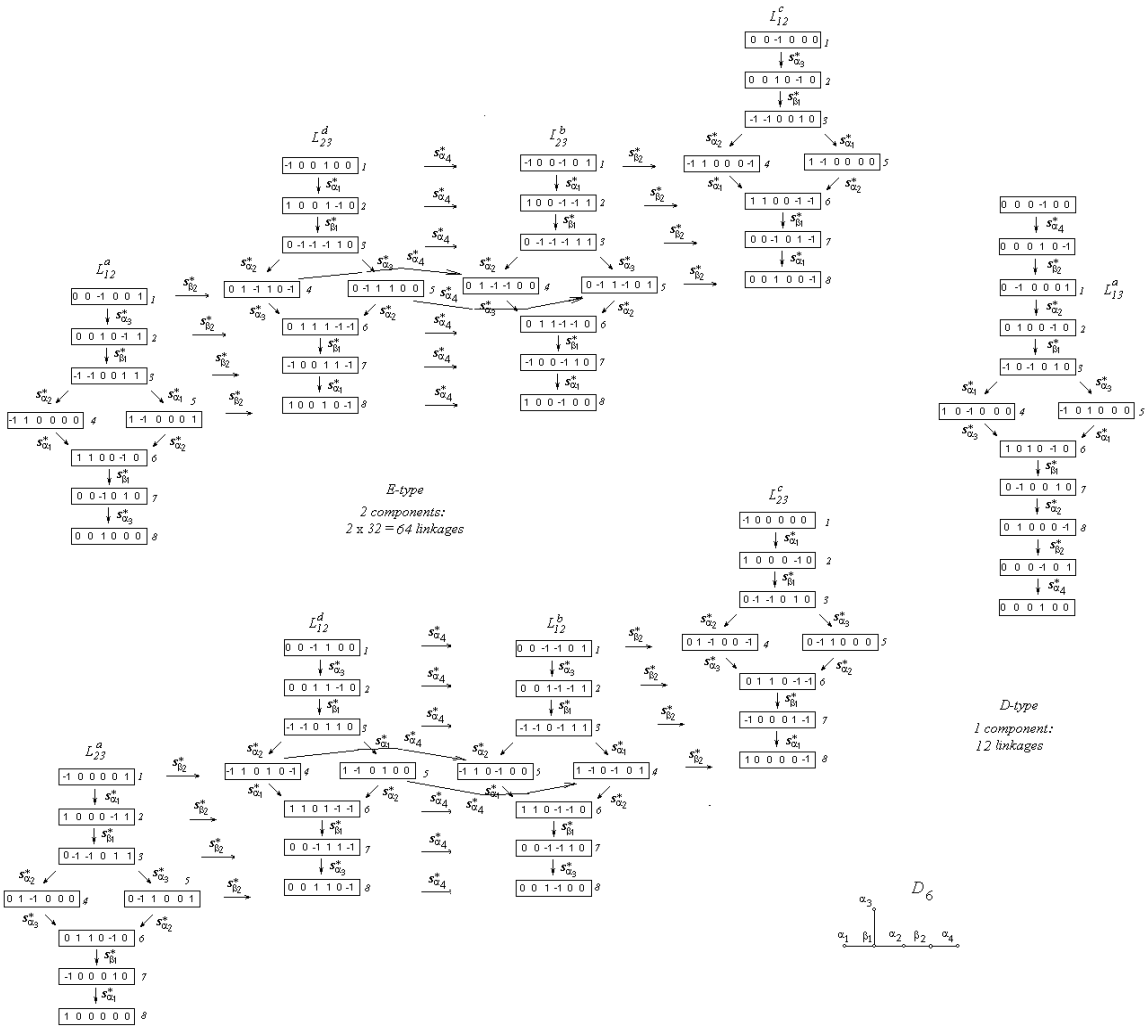}
\vspace{3mm}
\caption[\hspace{3mm}Linkage system for $D_6$, 3 components, 12 + 2x32 linkages]{\hspace{3mm}The linkage system  $D_6$, one $D$-type component containing $12$ elements and two components of the $E$-type containing $2\times32 = 64$ elements}
\label{D6pure_loctets}
\end{figure}

\begin{figure}[H]
\centering
\includegraphics[scale=0.85]{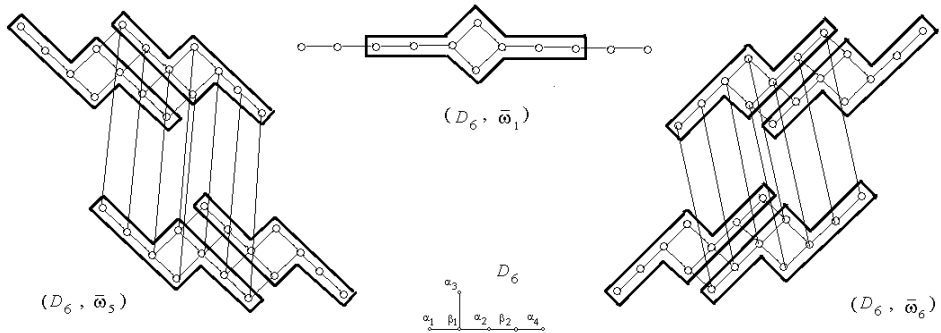}
\vspace{3mm}
\caption[\hspace{3mm}Loctets in the weight system for representations $(D_6, \overline{\omega}_i)$, 
$(D_6, \overline{\omega}_5$) and $(D_6, \overline{\omega}_6$)]{\hspace{3mm}Loctets in the weight system for $3$ fundamental representations of $D_6$:
  $(D_6, \overline{\omega}_1$), $(D_6, \overline{\omega}_5$) and $(D_6, \overline{\omega}_6$)  }
\label{32x2_12_weight_diagr_D6}
\end{figure}

\newpage
\subsection{The linkage systems $D_7(a_1)$, $D_7(a_2)$, $D_7$}
~\\

\begin{figure}[H]
\centering
\includegraphics[scale=0.42,  angle=270]{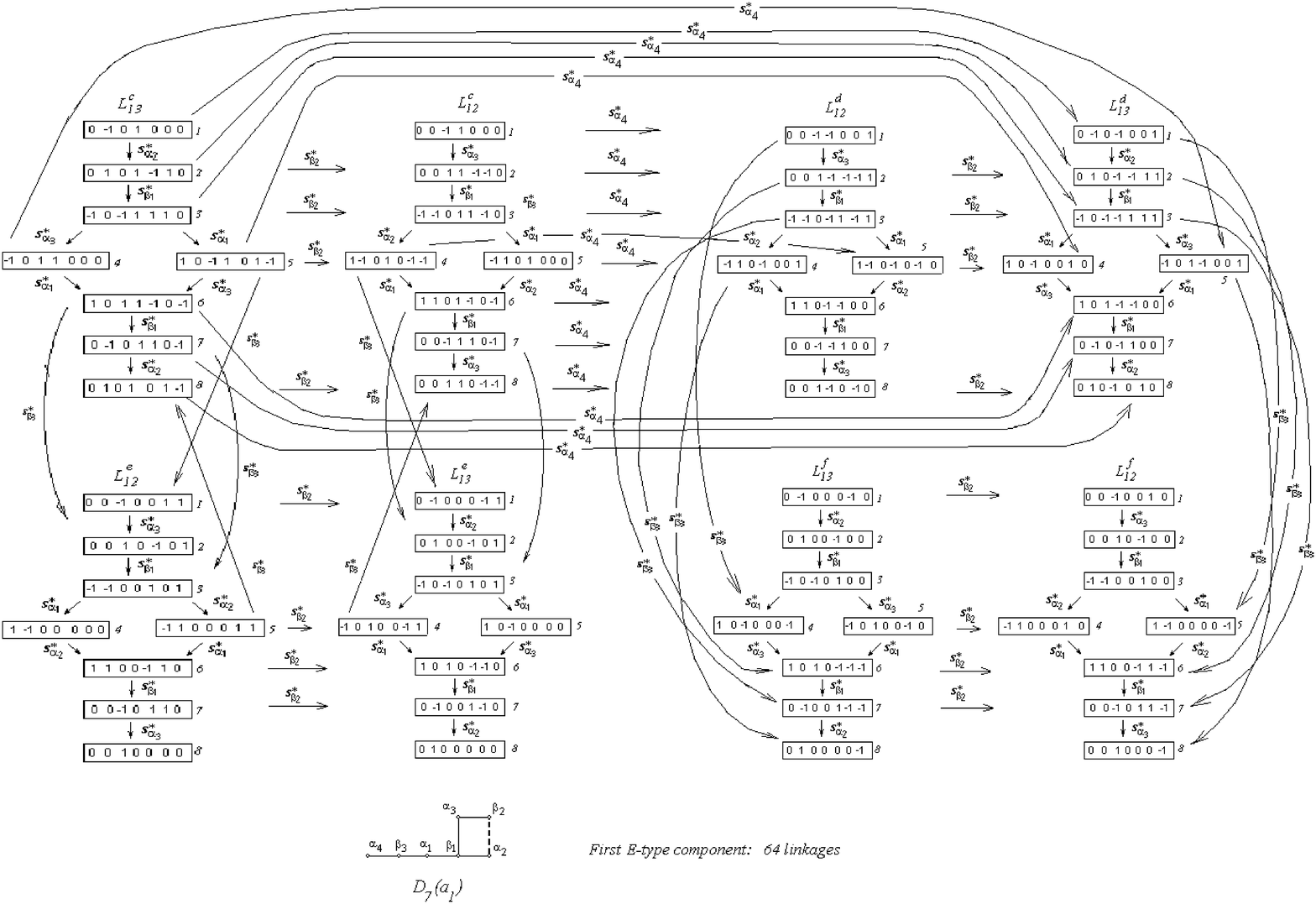}
\caption[\hspace{3mm}The linkage system $D_7(a_1)$, first component of $E$-type]{\hspace{3mm}The linkage system  $D_7(a_1)$, $1$st $E$-type component , $64$ linkage diagrams, $8$ loctets}
\label{D7a1_linkages_cdef}
\end{figure}

\begin{figure}[H]
\centering
\includegraphics[scale=0.47,  angle=270]{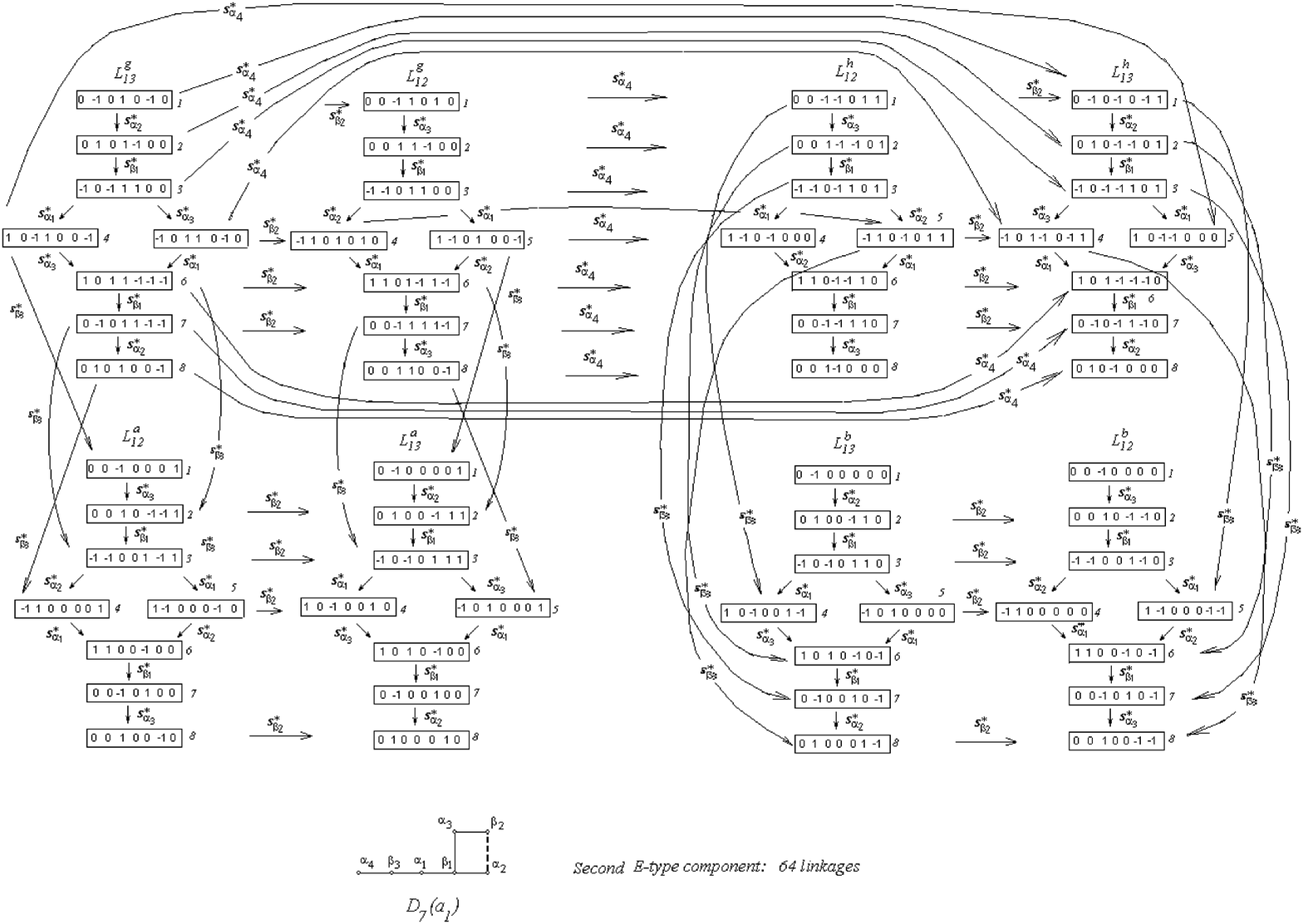}
\caption{\hspace{3mm}$D_7(a_1)$,second $E$-type component: $64$ linkages, $8$ loctets}
\label{D7a1_linkages_abgh}
\end{figure}


\begin{figure}[H]
\centering
\includegraphics[scale=1.25, angle=270]{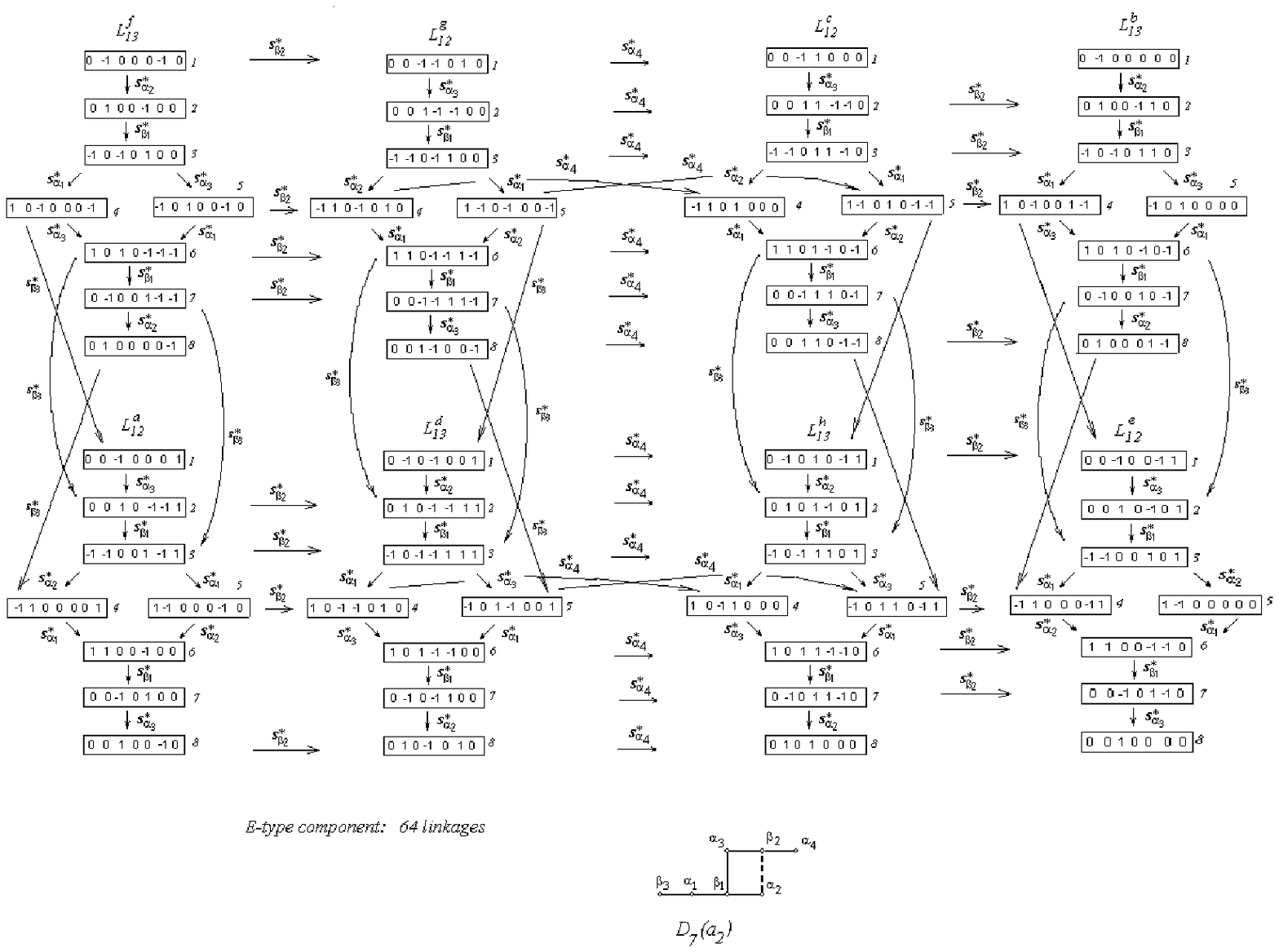}
\caption{\hspace{3mm}The linkage system $D_7(a_2)$, $1$st $E$-type component , $64$ linkage diagrams, $8$ loctets}
\label{D7a2_linkages_comp1}
\end{figure}


\begin{figure}[H]
\centering
\includegraphics[scale=0.55, angle=270]{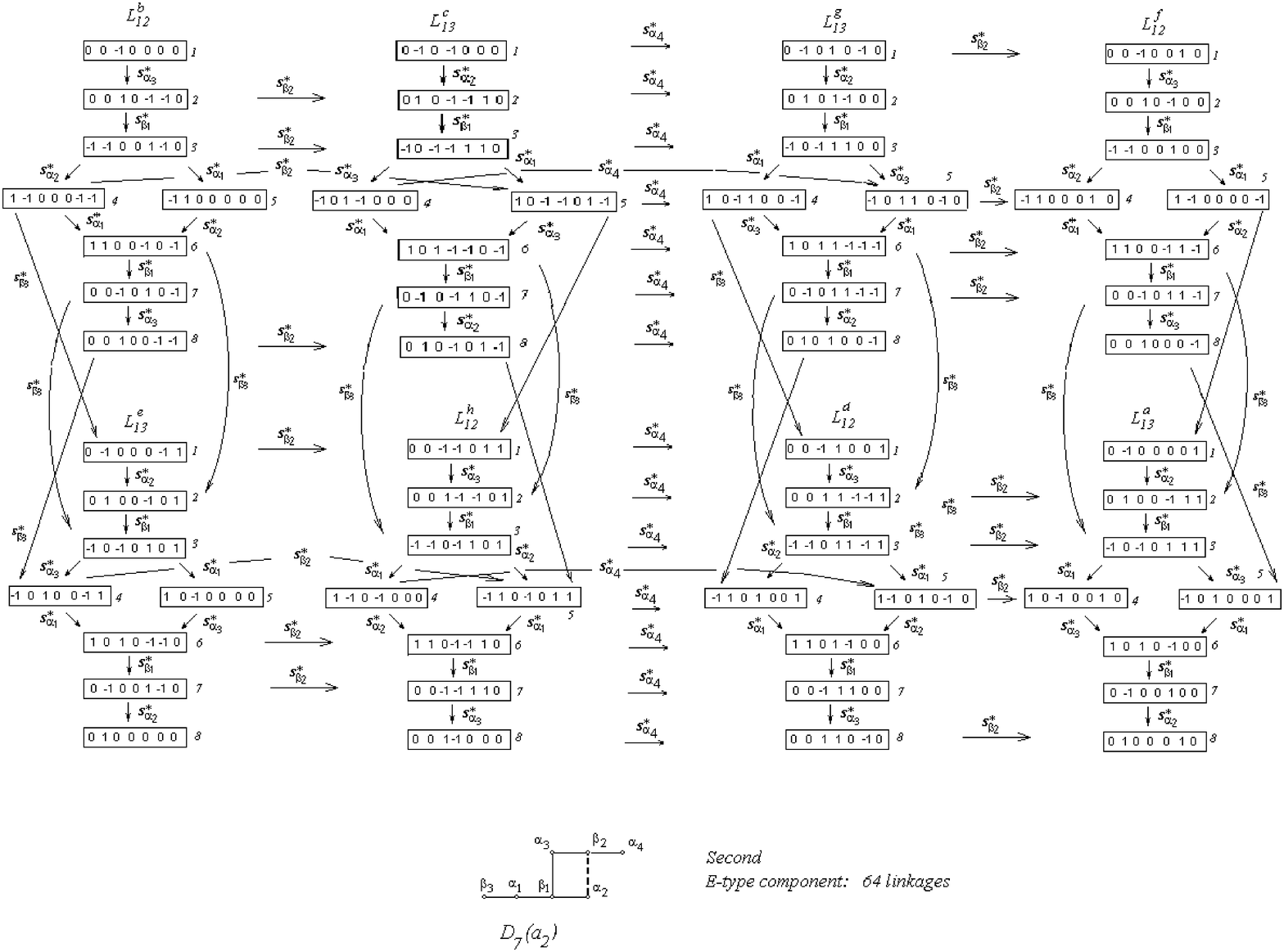}
\caption{\hspace{3mm} The linkage system $D_7(a_2)$, $2$nd $E$-type component, $64$ linkage diagrams, $8$ loctets}
\label{D7a2_linkages_comp2}
\end{figure}

\begin{figure}[H]
\centering
\includegraphics[scale=0.45, angle=270]{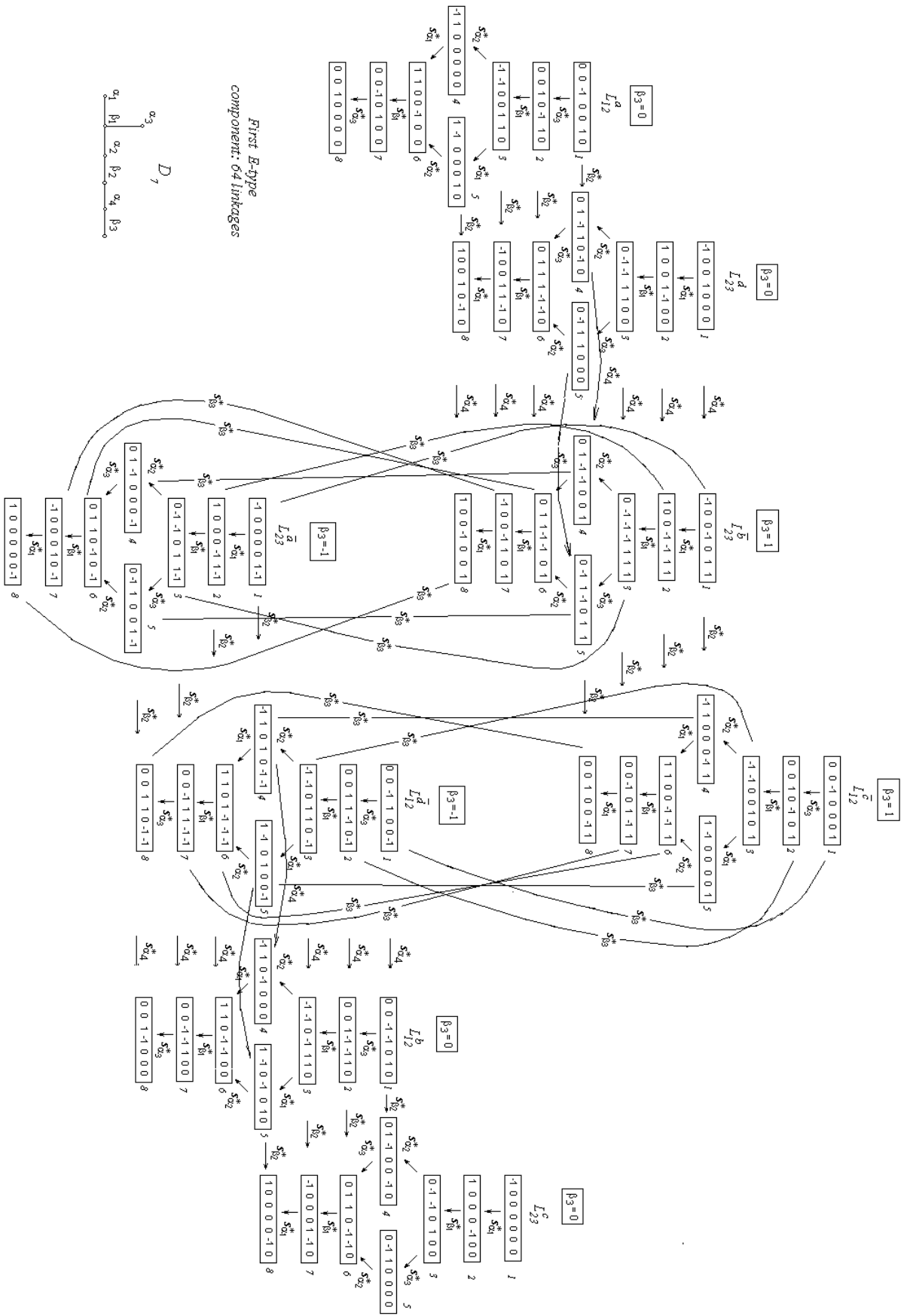}
\vspace{3mm}
\caption{\hspace{3mm}The first component of the linkage system $D_7$, 64 linkages}
\label{D7pu_loctets_comp1}
\end{figure}

\begin{figure}[H]
\centering
\includegraphics[scale=0.47, angle=270]{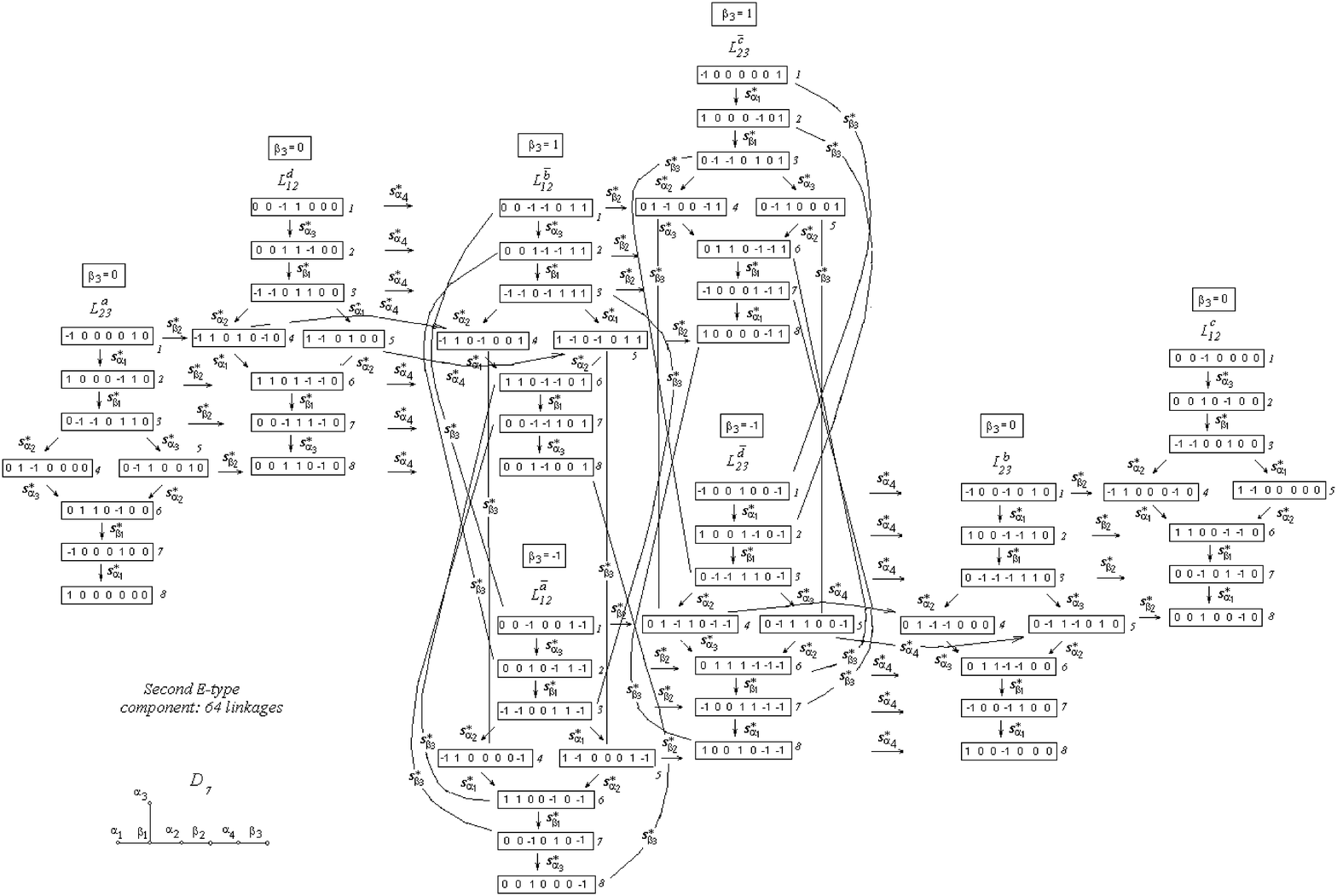}
\vspace{3mm}
\caption{\hspace{3mm}The second component of the linkage system $D_7$, 64 linkages}
\label{D7pu_loctets_comp2}
\end{figure}

\begin{figure}[H]
\centering
\includegraphics[scale=0.47, angle=270]{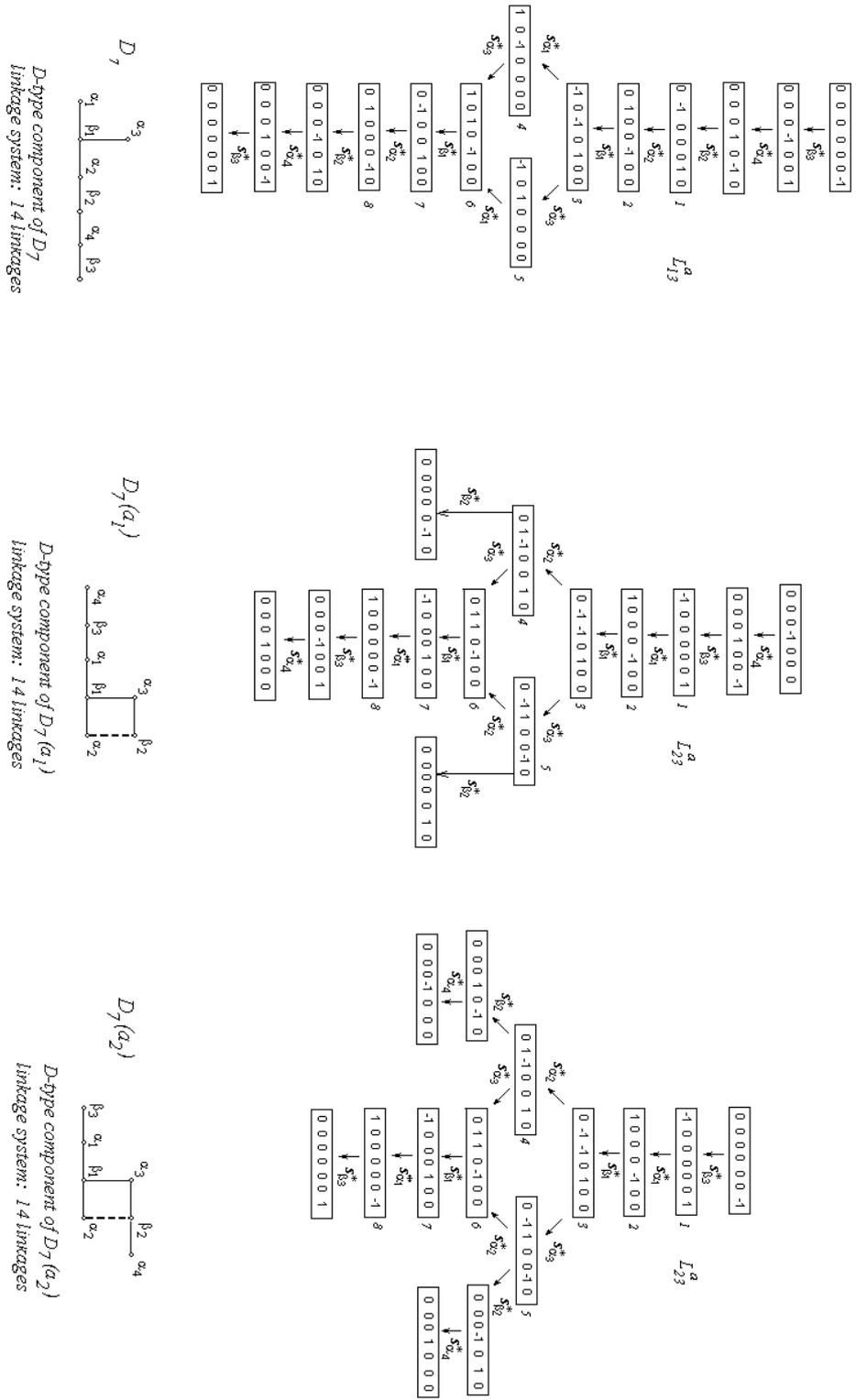}
\vspace{3mm}
\caption[\hspace{3mm}The D-type component for linkage systems
$D_7$, $D_7(a_1)$, $D_7(a_2)$; every D-type components contains 14 linkages]
{The D-type component for linkage systems $D_7$, $D_7(a_1)$, $D_7(a_2)$}
\label{D7a1_D7a2_D7pu_loctets_comp3}
\end{figure}

\newpage
\subsection{The linkage systems $D_l(a_k)$, $D_l$, $A_l$}
~\\

\begin{figure}[H]
\centering
\includegraphics[scale=0.4]{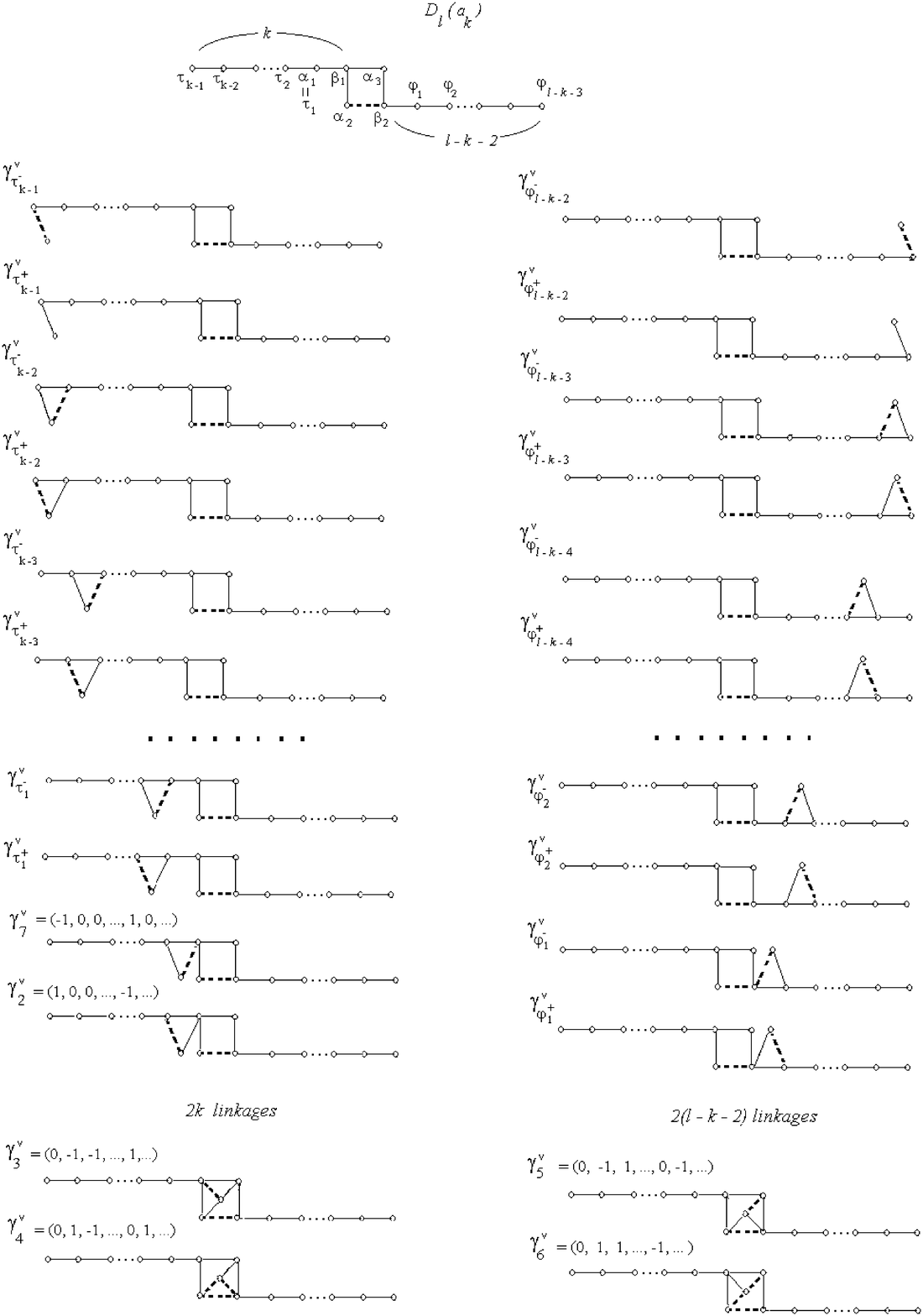}
\caption{\hspace{3mm}$D_l(a_k)$ for $l > 7$, $1$ loctet, $2l$ linkage diagrams}
\label{Dk_al_linkages}
\end{figure}

\begin{figure}[H]
\centering
\includegraphics[scale=1.3]{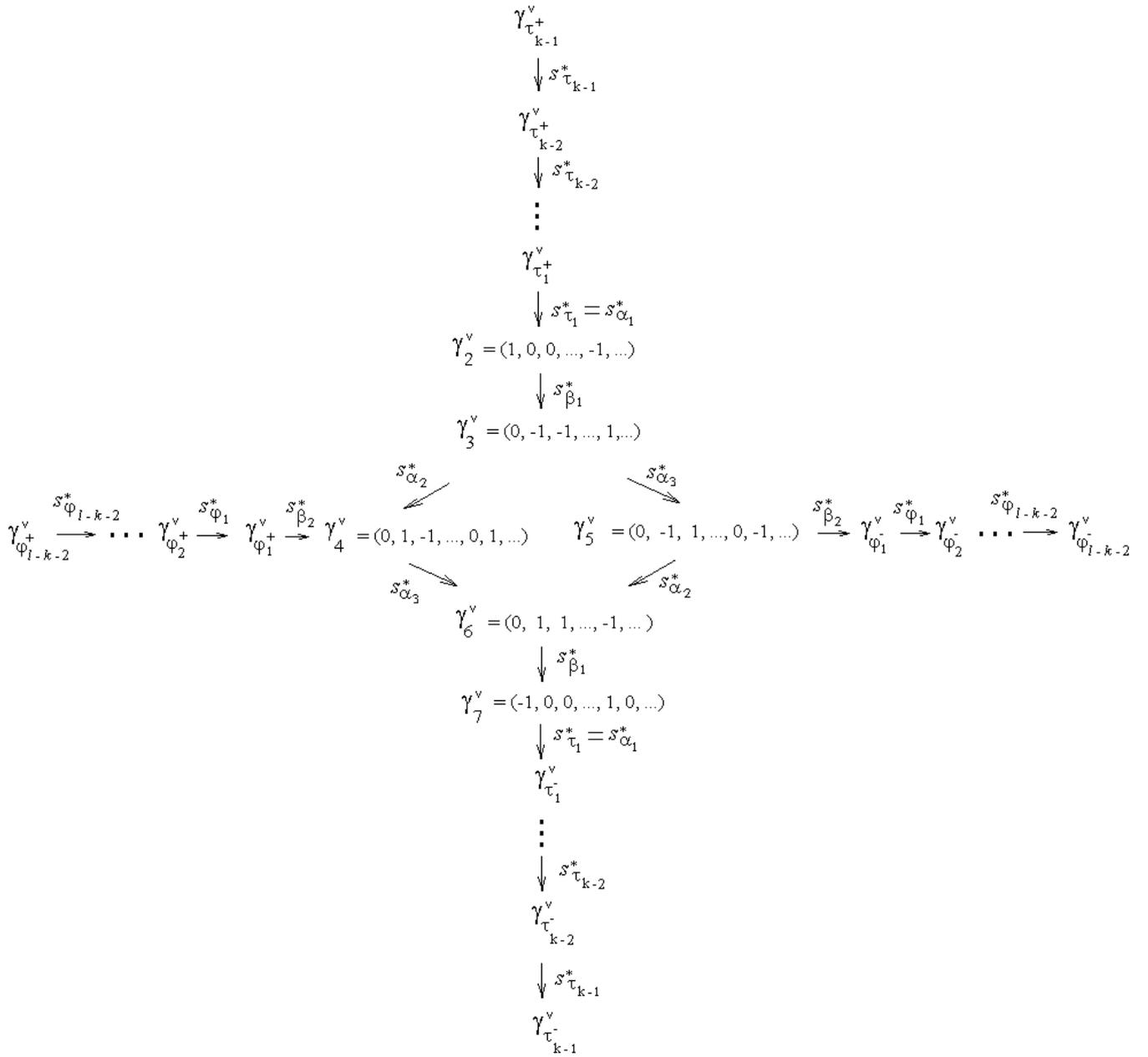}
\caption{\hspace{3mm}The linkage system $D_l(a_k)$ for $l > 7$ (wind rose of linkages)}
\label{Dk_al_wind_rose}
\end{figure}

\begin{figure}[H]
\centering
\includegraphics[scale=1.2, angle=270]{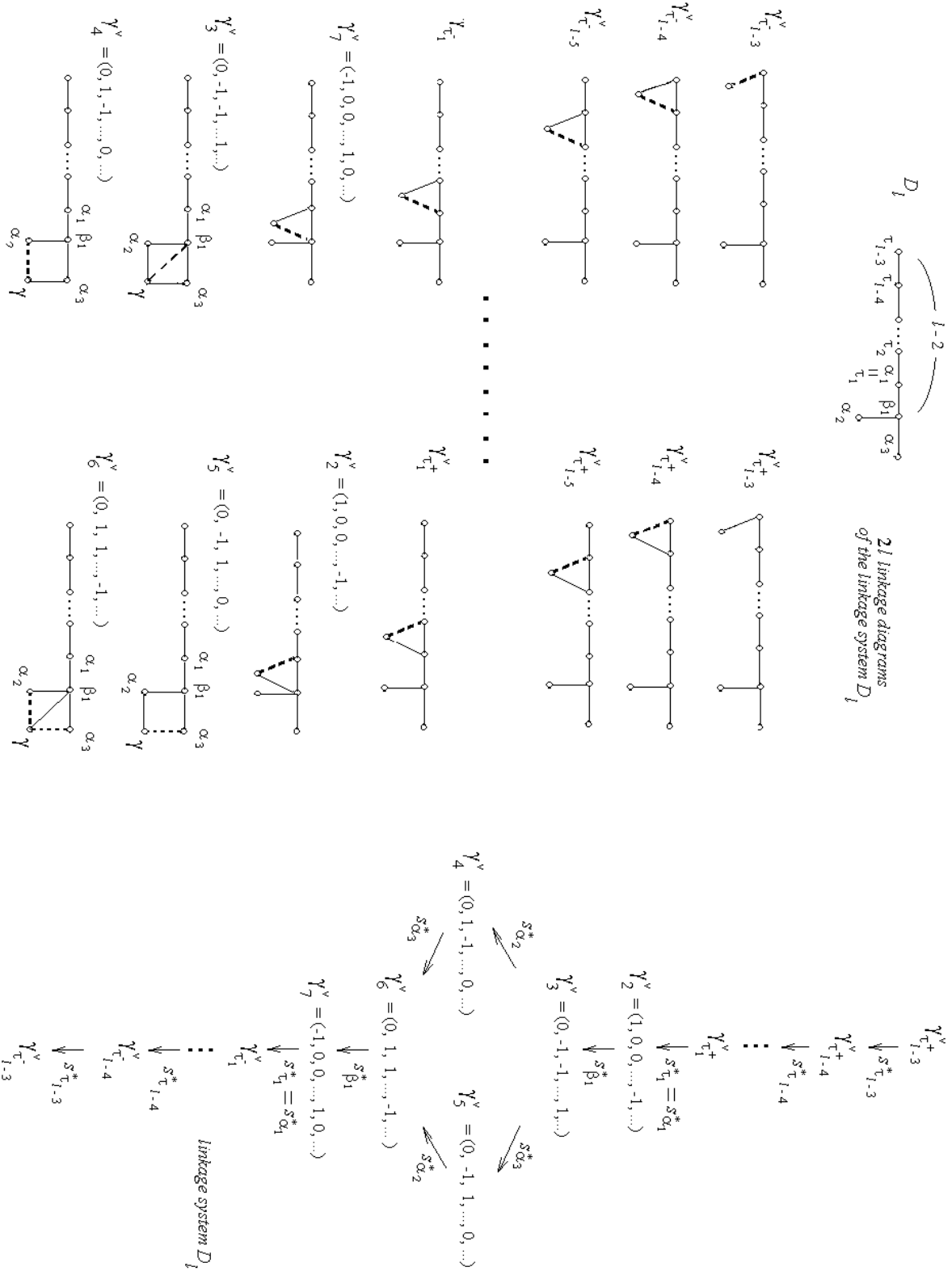}
\caption{\hspace{3mm}The linkage system $D_l$ for $l > 7$, $2l$ linkages, $1$ loctet}
\label{Dlpu_linkages}
\end{figure}

\begin{figure}[H]
\centering
\includegraphics[scale=1.8]{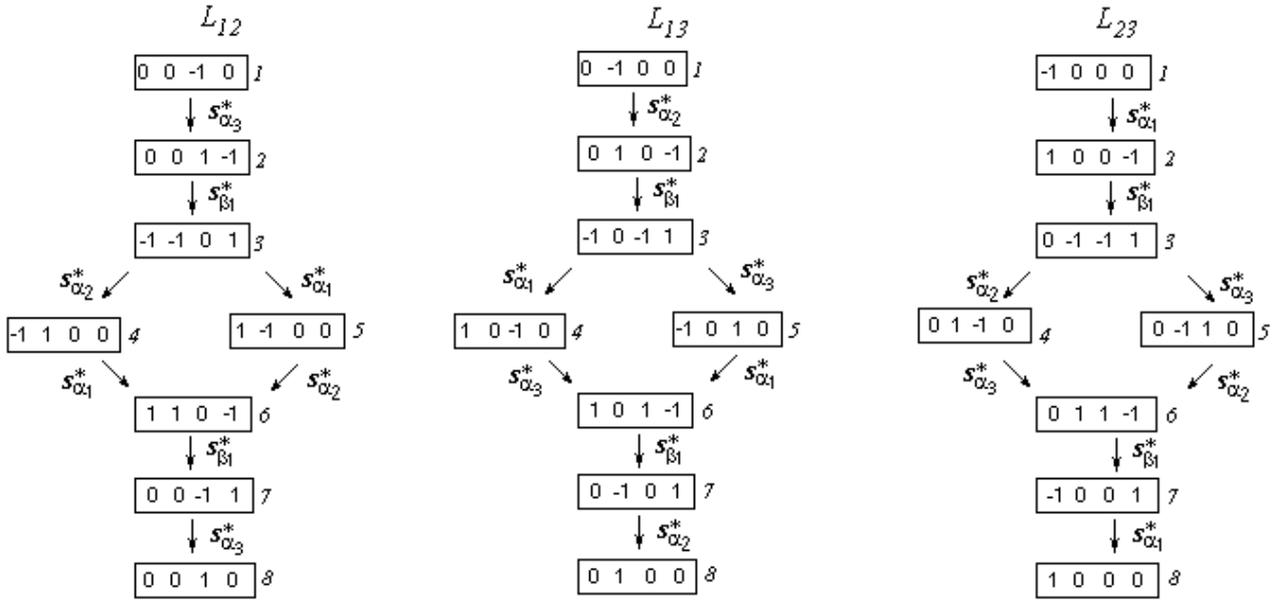}
\caption{\hspace{3mm}The linkage system $D_4$, $24$ linkages, $3$ loctets}
\label{D4_loctets}
\end{figure} 

%% file: notation.tex
\section*{List of Notations}
~\\

$\begin{array}{ll}
   W & \text{ the Weyl group associated with a Dynkin diagram } \\
   w = w_1w_2 & \text{ the bicolored decomposition of } w \in W; w_1, w_2 \text{ - involutions } \\
  \Gamma & \text{ the Carter diagram associated with the bicolored decomposition of } w \\
  \mathsf{DE4} & \text{ the class of simply-laced Dynkin diagrams containing $D_4$ as subdiagram}\\
  \mathsf{C4} & \text{ the class of simply-laced connected Carter diagrams
         containing $D_4(a_1)$ as subdiagram}\\
  \varPhi  & \text{ the root system associated with the Weyl group } W \\
  \Pi & \text{ the set of all simple roots of the root system } \varPhi \\
  E & \text{ the linear subspace spanned by the root system } \varPhi \\
  \Pi_w & \text{ the root subset (containing not necessary simple roots) associated with $w \in W$ } \\
  \Pi_w(\theta) & \text{ the root subset obtained from $\Pi_w$ by adding some $\theta \in \varPhi$ linearly independent of $\Pi_w$} \\
  L & \text{ the linear subspace spanned by the root subset } \Pi_w \\
  {\bf B} & \text{ the Cartan matrix associated with the Dynkin diagram } \Gamma \\
  B_L & \text{ the partial Cartan matrix associated with the Carter diagram } \Gamma \\
  (\cdot, \cdot) & \text{ the symmetric bilinear form associated with } {\bf B} \\
  (\cdot, \cdot)_{\botL} & \text{ the symmetric bilinear form associated with } B_L \\
  \mathscr{B} & \text{ the quadratic Tits form associated with } {\bf B} \\
  \mathscr{B}_L & \text{ the quadratic form associated with } B_L \\
  B^{-1}_L  & \text{ the inverse of for the partial Cartan matrix } B_L \\
  \mathscr{B}^{\vee}_L & \text{ the inverse quadratic form associated with } B_L \\
  W_L &  \text{ the partial Weyl group associated with the Carter diagram } \Gamma \\
  W_L(\theta) & \text{ the subgroup of the Weyl group $W$ generated by $W_L$ and $s_{\theta}$ } \\
  \mathscr{L}(\Gamma) & \text{ the linkage system associated with the Carter diagram $\Gamma$ } \\
  \varPhi_w(\theta) & \text{ the root subset is the orbit of the action of $W_L(\theta)$
      on $\Pi_w(\theta)$ } \\
  W^{\vee}_L &  \text{ the dual partial Weyl group associated with the Carter diagram } \Gamma \\
  \alpha_i, \beta_j &  \text{ the roots in the bicolored notation } \\
  \gamma & \text{ the linkage or $\gamma$-linkage } \\
  \gamma^{\vee} & \text{ the linkage diagram, or the linkage labels } \\
  \gamma^{\vee}(n) & \text{ the $n$-th linkage diagram of the loctet (in the vertical numbering $(n=1,\dots,8)$)} \\
  \gamma^{\vee}(8) & \text{ the eighth linkage diagram of the loctet} \\
  L_{ij}  & \text{ the loctet (=linkage octet) of type $ij \in \{ 12, 13, 23 \}$ }  \\
  \alpha\text{-set} & \text{ the subset of roots corresponding to $w_1$ in
            the bicolored decomposition } \\
  \beta\text{-set} & \text{ the subset of roots corresponding to $w_2$ in
            the bicolored decomposition } \\
  \alpha\text{-label} & \text{ any coordinate from $\alpha$-set of the linkage labels vector } \\
  \beta\text{-label} & \text{ any coordinate from $\beta$-set of the linkage labels vector } \\
  k & \text{ the number of coordinates in $\alpha$-set (= number of $\alpha$\text{-labels}) } \\
  h & \text{ the number of coordinates in $\beta$-set (= number of $\beta$\text{-labels}) } \\
  l & \text{ the number of vertices in the Carter diagram $\Gamma$, $l = k + h$ }
\end{array}$

%% file: biblio-2.tex
\pagebreak[4]

%% file: 1DiagrCalc-2_newd.bbl
\begin{thebibliography}{9999}

\bibitem[]{}

\bibitem[Bo02]{Bo02}
N.~Bourbaki, {\em Lie groups and Lie algebras. Chapters 4,5,6}.
Translated from the 1968 French original by Andrew Pressley. Elements of Mathematics (Berlin). Springer-Verlag, Berlin, 2002. xii+300 pp.

\bibitem[Bo05]{Bo05}
N.~Bourbaki, {\em Lie groups and Lie algebras. Chapters 7,8,9}.
Translated from the 1975 and 1982 French originals by Andrew Pressley. Elements of Mathematics (Berlin). Springer-Verlag, Berlin, 2005. xii+434 pp.

\bibitem[Ch84]{Ch84}
R.~N.~Cahn, {\em Semi-Simple Lie Algebras and Their Representations}.
Berkeley, Benjamin-Cummings publishing company, 1984.

\bibitem[Ca70]{Ca70}
R.~W.~Carter, {\em Conjugacy classes in the Weyl group}. 1970
Seminar on Algebraic Groups and Related Finite Groups (The Institute
for Advanced Study, Princeton, N.J., 1968/69) pp. 297--318 Springer,
Berlin.

\bibitem[Ca72]{Ca72}
R.~W.~Carter, {\em Conjugacy classes in the Weyl group}.  Compositio
Math.  25  (1972), 1--59

\bibitem[Co89]{Co89}
A.~J.~Coleman, {\em The greatest mathematical paper of all time}. Math. Intelligencer  11  (1989), no. 3, 29--38.

\bibitem[Dy50]{Dy50}
E.~B.~Dynkin, {\em Some properties of the system of weights of a linear representation of a semisimple Lie group}. (Russian) Doklady Akad. Nauk SSSR (N.S.)  71,  (1950). 221--224.

\bibitem[Dy52]{Dy52}
E.~B.~Dynkin, {\em Maximal subgroups of the classical groups}. (Russian) Trudy Moskov. Mat. Obsh. 1,  (1952). 39--166.

\bibitem[FOT00]{FOT00}
M.~Fukuma, T.~Oota, H.~Tanaka, {\em Weyl groups in $\rm AdS\sb 3/CFT\sb 2$}. Progr. Theoret. Phys.  103  (2000), no. 2, 447--462

\bibitem[GOV90]{GOV90}
V.~V.~Gorbatsevich, A.~L.~Onishchik, E.~B.~Vinberg, {\em Structure of Lie groups and Lie algebras}.
(Russian)  Current problems in mathematics. Fundamental directions, Vol. 41 (Russian),  5--259,
Itogi Nauki i Tekhniki, Akad. Nauk SSSR, Vsesoyuz. Inst. Nauchn. i Tekhn. Inform., Moscow, 1990.
English Translation: {\em Lie groups and Lie algebras III}, Encyclopediya of Mathematical Sciences, v. 41.

\bibitem[Kac80]{Kac80}
V.~Kac, {\em Infnite  root  systems,  representations  of  graphs  and  invariant  theory}.
Invent. Math. 56 (1980), no. 1, 57–-92.


\bibitem[KOV95]{KOV95}
F.~I.~Karpelevich, A.~L. Onishchik, and E.~B.~Vinberg, {\em On the
work of E. B. Dynkin in the theory of Lie groups in: Lie Groups and
Lie Algebras: E.~B.~Dynkin's Seminar}, 1--12,
 Amer. Math. Soc. Transl. Ser. 2, 169, Amer. Math. Soc. 1995,
202 pp, vol. 169.

\bibitem[K70]{K70}
G.~F.~Kushner, {\em On the compactification of noncompact symmetric spaces.} Dokl. Akad. Nauk, 190 (1970), 1282--1285, English transl., Sov. Math. Docl. 11 (1970), no. 1, 284-287.

\bibitem[K72]{K72}
G.~F.~Kushner, {\em The compactification of noncompact symmetric spaces}, Trudy Sem.Vektor.Tenzor.Anal. 16 (1972), 99-152.

\bibitem[K79]{K79}
G.~F.~Kushner, {\em F.I.Karpelevich's compactification is homeomorphic to a ball.} Trudy Sem. Vektor. Tenzor. Anal., 19 (1979), 96--111; English transl., Amer.Math.Soc.Translations, ser.2, 134 (1987), 119-133.

\bibitem[PSV98]{PSV98}
E.~Plotkin, A.~Semenov, N.~Vavilov, {\em Visual basic representations: an atlas}. Internat. J. Algebra Comput.  8  (1998), no. 1, 61--95.

\bibitem[Sl81]{Sl81}
R.~Slansky, {\em Group theory for unified model building}.
Phys. Rep.  79  (1981), no. 1, 1--128.

\bibitem[St08]{St08}
R.~Stekolshchik, {\em Notes on Coxeter Transformations and the McKay
Correspondence}, Springer Monographs in Mathematics,
 2008, XX, 240 p.

\bibitem[St10]{St10}
R.~Stekolshchik, {\em Root systems and diagram calculus. I. Cycles in the Carter diagrams},
arXiv:1005.2769v3.

\bibitem[St11]{St11}
R.~Stekolshchik, {\em Root systems and diagram calculus. III. Semi-Coxeter orbits of linkage diagrams and the Carter theorem}, arXiv:1105.2875v1.

\bibitem[S07]{S07}
I.~Stewart, {\em Why beauty is truth: a history of symmetry}, Basic Books, 2007.

\bibitem[Va00]{Va00}
N.~Vavilov, {\em A third look at weight diagrams}. (English summary)
Rend. Sem. Mat. Univ. Padova 104 (2000), 201--250.

\end{thebibliography}
